\documentclass[a4, 12pt]{article}
\usepackage{amssymb}
\usepackage{yfonts}
\def\pointir{\unskip . --- \ignorespaces}
\newtheorem{guess}{Conjecture}

\newtheorem{proposition}{Proposition}[section]
\newtheorem{corollary}{Corollary}[section]
\newtheorem{lemma}{Lemma}[section]
\newtheorem{theorem}{Theorem}[section]

\makeatletter

\title{\bf Multiplicity one Conjectures}
\author{\textbf {Steve Rallis\hskip 0.3cm and\hskip 0.3 cm G\'erard Schiffmann}}
\begin{document}

\def\section
{\@startsection {section}{1}
{-1pt}{-5ex \@plus -1ex \@minus -.2ex}
{2ex \@plus .2ex}
{\normalfont \large \bfseries } }
\maketitle

\vskip 0.5cm
 \begin{abstract}
In the first part, in the local non archimedean case,  we consider distributions on $GL(n+1)$ which are invariant under the adjoint action of  $GL(n)$. We conjecture that such distributions are invariant by transposition. This would imply multiplicity at most one for rerstrictions from $GL(n+1)$ to $GL(n)$. We reduce ourselves to distributions with "singular" support and then finish the proof for $n\leq 8$.\hfill\break
In the second part we show that similar Theorems for orthogonal or unitary groups follow from the case of $GL(n)$
\end{abstract}

 \footnote{The first author is partially supported by NSF Grant DMS-0500392}
 
 \centerline{\bf Introduction}
 \vskip 0.5cm
Let ${\mathbb F }$ be a local field non archimedean and of characteristic 0. Let $W$ be a vector space over ${\mathbb F}$ of finite dimension $n+1\geqslant  2$ and let $W=V\oplus U$ be a direct sum decomposition with $\dim V=n $. Then we have an imbedding of $GL(V)$ into $GL(W)$. Our goal is to study the following Conjecture:

\vskip 0.5cm
{\parindent 0pt
{\bf Conjecture $\mathbf 1$:}
 {\sl If $\pi$ (resp. $\rho$) is an irreducible admissible representation of $GL(W)$ (resp. of $GL(V)$) then
\[ \mathrm{dim} \left ({\mathrm {Hom}}_ {GL(V)}(\pi _{|GL(V)}\, ,\, \rho )\right )\leqslant 1\] }}

We choose a basis of $V$ and a non zero vector in $U$ thus getting a basis of $W$. We can identify $GL(W)$ with $GL(n+1,{\mathbb F})$ and $GL(V)$ with $GL(n,{\mathbb F})$. The transposition map is an involutive anti automorphism of $GL(n+1,{\mathbb F)}$ which leaves $GL(n,{\mathbb F})$ stable. It acts on the space of distributions on $GL(n+1,{\mathbb F})$. 

Conjecture 1 is a Corollary of :
\vskip 0.5cm
{\parindent 0pt
{\bf Conjecture $\mathbf 2$:}
{\sl  A distribution on $GL(W)$ which is invariant under the adjoint action of $GL(V)$ is invariant by transposition.}}
\vskip 0.5cm
To try to prove this conjecture we proceed by induction on $n$. We use Harish-Chandra's descent and another type of descent, of a more elementary nature, similar to the argument used to get Frobenius reciprocity. This allows us to show that an invariant distribution, skew with respect to transposition must be supported by the ''singular set'' ( strictly speaking we first linearize the problem). We have not been able to produce a complete proof that such a distribution must be 0, except for very small values of $n$ ( precisely $n\leqslant 8$). At this stage we can not rule out the possibility that Conjecture 2 is false in general, thus casting a doubt on Conjecture 1.

On the opposite we are still confident that this does not happen and that the missing part of the proof should not lie too deep.

One can raise a similar question for orthogonal and unitary groups. Let ${\mathbb D}$ be a  either ${\mathbb F}$ or a quadratic extension of ${\mathbb F}$ . If $x\in{\mathbb D}$ then $\overline x$ is the conjugate of $x$ if ${\mathbb D}\neq{\mathbb F}$ and is equal to $x$ if ${\mathbb D}={\mathbb F}$.

Let $W$ be a vector space over ${\mathbb D}$ of finite dimension $n+1\geqslant 2$.  Let $\langle  .,.\rangle  $ be a non degenerate hermitian form on $W$. This form is bi-additive and
\[ \langle  dw,d'w'\rangle  =d\,\,\overline{d'}\langle  w,w'\rangle   ,\quad \langle  w',w\rangle  =\overline{\langle  w,w'\rangle  }\]
Given a ${\mathbb D}-$linear map $u$ from $W$ into itself, its adjoint $u^*$ is defined by the usual formula
\[ \langle  u(w),w'\rangle  =\langle  w,u^*(w')\rangle  \]

Choose a vector $e$ in $W$ such that $\langle  e,e\rangle  \ne 0$; let $U={\mathbb D}e$ and $V=U^{\perp}$ the orthogonal complement. Then $V$ has dimension $n$ and the restriction of the hermitian form to $V$ is non degenerate.

Let $M$ be the unitary group of $W$ that is to say the group of all ${\mathbb D}-$linear maps $m$ of $W$ into itself which preserve the hermitian form or equivalently such that $mm^*=1$. Let $G$ be the unitary group of $V$. With the p-adic topology both groups are of type lctd ( locally compact, totally discontinuous and countable at infinity). They are reductive groups of classical type.

The group $G$ is naturally imbedded into $M$. \vskip 0.5cm

{\parindent 0pt
{\bf Conjecture $\mathbf 1'$:} 
 {\sl If $\pi $ (resp $\rho $) is an irreducible admissible representation of $M$ (resp of $G$) then
\[ \mathrm{dim}\left  ({\mathrm {Hom}}_G(\pi _{| M},\rho )\right  )\leq 1\]
}}
Choose a basis $e_1,\dots e_n$ of $V$ such that $\langle  e_i,e_j\rangle  \in{\mathbb F}$. For 
\[ w=x_0e+\sum _1^nx_ie_i\]
put
\[ \overline w=\overline xe+\sum _1^n\overline {x_i}\, e_i\] 
If $u$ is a ${\mathbb D}-$linear map from $W$ into itself, let $\overline u$ be defined by
\[ \overline u(w)=\overline{u(\overline w)}\]

Let $\sigma $ be the anti-involution $\sigma (m)=\overline m^{-1}$ of $M$;  Conjecture 1' is a consequence of \vskip 0.5cm

{\parindent 0pt
{\bf Conjecture $\mathbf 2'$:} 
{\sl A distribution on $M$ which is invariant under the adjoint action of $G$ is invariant under $\sigma $.
}}
\vskip 0.5cm
We proceed exactly as in the case of the general linear group. Note that the Levi subgroups of some parabolic subgroups of $G$ have components of type $GL$ so that we need to assume Conjecture 2. However, in this case, the singular set has a natural stratification stable by the involution and it is possible to use induction on each strata, thus finishing the proof that Conjecture 2' follows from Conjecture 2.

The above questions are not new and Conjectures 1 and 1' have been around for some time. At the beginning stage of this project the first author had conversations with J.Bernstein who was interested in the same problem.

Assuming multiplicity at most one, a more difficult question is to find when it is one. Some partial results are known.

For the orthogonal group (in fact the special orthogonal group) this question has been studied by B. Gross and D.Prasad (\cite{1},\cite{8}) who formulated a precise conjecture. An up to date account  is given by B.Gross and M.Reeder (\cite{2}). Perhaps their approach could eventually give a proof of Conjecture 1'.  In a different setup, in their work on "Shintani" functions A.Murase and T.Sugano obtained complete results for $GL(n)$ and the split orthogonal case but only for spherical representations (\cite{5},\cite{7}). Finally we should mention, Hakim's publication \cite{3}, which, at least for the discrete series, could perhaps lead to a proof of the Conjectures.

This work is divided in two parts (General Linear Group and Orthognal/Unitary Groups), each having a short introduction and an Appendix where we present a certainly well known elementary descent method based on an argument \`a la Frobenius. 
\vfill\eject
\centerline{\normalfont \Large \bfseries  Part I}
\vskip 0.5cm
\centerline{\normalfont \Large \bfseries  General Linear Group}
\vskip 1cm{}

\section {\hglue -12pt . Introduction}

Let ${\mathbb F }$ be a local field non archimedean and of characteristic 0. Let $W$ be a vector space over ${\mathbb F}$ of finite dimension $n+1\geqslant 2$ and let $W=V\oplus   U$ be a direct sum decomposition with $\dim V=n $. Then we have an imbedding of $GL(V)$ into $GL(W)$. Our ( at present unachieved) goal is to prove the following Conjecture:
\begin{guess}
 If $\pi $ (resp. $\rho $) is an irreducible admissible representation of $GL(W)$ (resp. of $GL(V)$) then
\[ \dim \left  ({\mathrm {Hom}}_ {GL(V)}(\pi _{|GL(V)}\, ,\, \rho )\right  )\leq 1\] 
\end{guess}
We choose a basis of $V$ and a non zero vector in $U$ thus getting a basis of $W$. We can identify $GL(W)$ with $GL(n+1,{\mathbb F})$ and $GL(V)$ with $GL(n,{\mathbb F})$. The transposition map is an involutive anti automorphism of $GL(n+1,{\mathbb F)}$ which leaves $GL(n,{\mathbb F})$ stable. It acts on the space of distributions on $GL(n+1,{\mathbb F})$. 

Conjecture 1 is a Corollary of :
\begin{guess}
 A distribution on $GL(W)$ which is invariant under the adjoint action of $GL(V)$ is invariant by transposition.
\end{guess}
{\bf Remarks\pointir}The transposition depends on the choice of the basis. However, changing the basis of $V$ amounts to the adjoint action of some element of $GL(V)$ and so does not change the action on $GL(V)-$invariant distributions. Changing the choice of the non zero vector in $U$ amounts to the adjoint action of some element of $GL(U)$. However, because the adjoint action of the center of $GL(W)$ is trivial  the adjoint action of $GL(U)$ reduces to the adjoint action of the center of $GL(V)$ and so is trivial on the space of $GL(V)-$invariant distributions. In other words a distribution invariant under the adjoint action of $GL(V)$ is in fact invariant under the adjoint action of $GL(V)\times GL(U)$ and so, on this space of distributions, the action of the "transposition" is independant the choice of the basis, adapted to the decomposition $W=V\oplus   U$. 

We now describe in more details the content of this part. In section 2 we prove that Conjecture 1 follows from Conjecture 2. In section 3 we perform a partial linearization. Namely we consider the action of $G=GL(V)$ on $G\times V\times V^*$ with $V^*$ the dual space of $V$; the group $G$ acts on itself by the adjoint action, on $V$ by the natural action and on $V^*$ by the inverse of the transpose. If $(e_i)$ is a basis of $V$ and $(e_i^*)$ the dual basis of $V^*$ define the linear map $u$ from $V$ to $V^*$ by $u(e_i)=e_i^*$ and consider the involution
\[ \sigma : (g,v,v^*)\mapsto (u^{-1}\,\,{}^tgu,u^{-1}(v^*),u(v))\] 
We show that Conjecture 2 is implied by
\begin{guess}
A distribution on $G\times V\times V^*$, invariant under $G$ is invariant under $\sigma $.
\end{guess}
Using the adjoint action on the Lie algebra ${\mathfrak {g}}$ of $G$ and  defining $\sigma $ as before, we then consider
\begin{guess} 
A distribution on ${\mathfrak{g}}\times V\times V^*$, invariant under $G$ is invariant under $\sigma $.
\end{guess}
In section 4 we prove Conjecture 4 for the case $n=1$, thus starting the induction. In section 5 we review Harish-Chandra's descent and first prove that Conjecture 4 implies Conjecture 3. 

Using Harish-Chandra's descent on the Lie algebra we get that an invariant distribution $T$ on ${\mathfrak { g}}\times V\times V^*$, skew-symmetric with respect to $\sigma $, has a support contained in $({\mathcal {Z}}\times {\mathcal {N}})\times V\times V^*$ where ${\mathcal {Z}}$ is the center of ${\mathfrak { g}}$ and ${\mathcal {N}}$ the cone of nilpotent elements.

The next step is to study the support on the $V\oplus V^*$ side. In section 6 we define and study the regular orbits. A triplet $(X,v,v^*)$ is regular if $(v,Xv,\dots ,x^{n-1}v)$ is a basis of $V$ and if $(v^*,^tXv^*,\dots ,^tX^{n-1}v^*)$ is a basis of $V^*$. The set of regular elements is a non empty Zariski open subset, stable by $G$. The orbits of the regular elements are the "regular" orbits. They turn out to be closed, with trivial centralizer and stable by $\sigma $. A well known result of Gelfand and Kazhdan then tells us that an invariant distribution on the set of regular elements is stable by $\sigma $.

A triplet  $(X,v,v^*)$ is "singular" if $X\in {\mathcal {N}}$ and if, for all $i$ we have $\langle  v^*,X^iv\rangle  =0$. In section 7 we show that if the triplet is non singular then the situation may be reduced to a mixture of the regular case and a case of lower dimension. We then can conclude that an invariant skew symmetric distribution must have a support contained in the singular set ( note that ${\mathcal {Z}}$ does not play any role).

In the remainder of this first Part we explain how far we can go toward proving that an invariant distribution with singular support and skew symmetric is 0. First in section 8 we collect further properties of these distributions. They are distributions on $[\mathfrak g,\mathfrak g]\times V\times V^*$ . Define the Fourier transform on $[\mathfrak g,\mathfrak g]$ using the Killing form. The   Fourier transform of the distribution has the same invariance properties. In particular its support must be singular. Corresponding to the Killing form there is a representation of the twofold covering of $SL(2)$ and the distribution is invariant.For $n$
even this implies that the distribution is 0 and for $n$ odd that it is homogeneous and equal to its Fourier transform. A similar remark is valid for the quadratic form $\langle v^*,v\rangle$ on $V\oplus V^*$. We fix a nilpotent orbit and transfer the problem to $V\oplus V^*$.

In sections 9 and 10, working on $V\oplus V^*$, we have a fix nilpotent matrix $X$ and the distribution is invariant under the centralizer $C$ of $X$. We try to work by  induction on ÒnÔ( independantly of our general induction) but "going down" we get only invariance under a subgroup of the centralizer of the new nilpotent matrix. This kills the hope of an easy proof. Still it works in a limited number of cases and we managed to finish the proof for $n\leq 8$.

Finally section 11 outlines a completely different approach. We show that it would be enough to establish a symmetry property either of $C$-invariant linear forms on the space of a finite number of  representations belonging to a degenerate principal series or of $G$ intertwinning maps between the Schwartz space of the nilpotent cone and these degenerate principal series.

\section{\hglue -12pt . Conjecture 2 implies Conjecture 1}

This is a well known argument but we recall the proof. A group of type lctd is  a locally compact, totally discontinuous group which is countable at infinity. We consider  smooth representations of such groups. If $(\pi ,E_\pi )$ is such a representation then $(\pi ^*,E_\pi ^*)$ is the smooth contragredient. Smooth induction is denoted by $Ind$ and compact induction by $ind$. For any topological space $T$ of type lctd, $\mathcal {S}(T)$ is the space of functions locally constant, complex valued,  defined on $T$ and with compact support.
\begin{proposition}
 Let $M$ be a tdlc group and $N$ a closed subgroup, both unimodular. Suppose that there exists an involutive anti-automorphism $\sigma $ of $M$ such that $\sigma (N)=N$ and such that any distribution on $M$, biinvariant under $N$, is fixed by $\sigma $. Then, for any irreducible admissible representation $\pi $ of $M$
\[ {\mathrm {dim}} \left  ({\mathrm {Hom}} _ M(ind_ N^M(1),\pi )\right  )\times {\mathrm {dim}} \left  ({\mathrm {Hom}} _ M(ind_ N^M(1),\pi ^*)\right  )\leq 1\]
\end{proposition}

 Let $E_\pi $ be the space of $\pi $ and $E^*_\pi $ its smooth dual. Let $A$ (resp. $A'$) be a linear map from ${\mathcal {S}}(M/N)$ into $E_\pi $  (resp. $E^*_\pi $) commuting with $M$. For $f_1,f_2\in {\mathcal{ S}}(M/N)$ put
 \[ B(f_1,f_2)=\langle  Af_1,A'f_2\rangle  \] 
 Choose Haar measures on $M$ and $N$. For $\varphi \in {\mathcal {S}}(M)$ define
 \[ \varphi ^\sharp(mN)=\int _N\varphi (mn)dn\] 
 This is a surjective map from ${\mathcal {S}}(M)$ onto ${\mathcal {S}}(M/N)$. Then
 \[ (\varphi _1,\varphi _2)\mapsto B(\varphi _1^\sharp,\varphi _2^\sharp)\]
 is a bilinear form on $M\times M$ left invariant under the action $(m_1,m_2)\mapsto (mm_1,mm_2)$ of $M$ onto $M\times M$ and right invariant under the action of $N\times N$. View it as a distribution $S$ on $M\times M$ and let us use integral notations
 \[ B(\varphi _1^\sharp,\varphi _2^\sharp)=\int _{M\times M}\varphi _1(m_1)\varphi _2(m_2)dS(m_1,m_2)\] 
 
 Next we consider the homeomorphism of $M\times M$ onto itself;
 \[ \Phi  : \,\,\,(m_1,m_2)\mapsto (m_2^{-1}m_1,m_2)\] 
 The distribution
 \[ \varphi \mapsto \int _{M\times M}\varphi \circ \Phi  (m_1,m_2)dS(m_1,m_2)=\int _{M\times M}\varphi (m_2^{-1}m_1,m_2)dS(m_1,m_2)\] 
 is invariant under the action of $M$ by left translations on the second variable and is biinvariant under the action of $N$ acting on the first variable by left and right multiplication. Hence there exists on $M$ a distribution $T$, biinvariant under $N$ and such that
 \[ \int _{M\times M}\varphi (m_2^{-1}m_1,m_2)dS(m_1,m_2)=\int _{M\times M}\varphi (m_1,m_2)dm_2\, dT(m_1)\]
 Suppose that $\varphi $ is decomposed $\varphi (m_1,m_2)=\varphi _1(m_1)\varphi _2(m_2)$. Let $\tau : m\mapsto \sigma (m)^{-1}$; it is an automorphism of $M$ leaving $N$ stable. Define $\varphi _i^\sigma (m)=\varphi _i(\sigma (m))$,$\,\,\check \varphi _i(m)=\varphi _i(m^{-1})$ and $\varphi _i^\tau (m)=\varphi _i(\tau (m))$. Then
 \begin{eqnarray*}
\int _{M\times M}\varphi _1(m_1)\varphi _2(m_2) dS(m_1,m_2)&=&\int _{M\times M}\varphi _1(m_2m_1)\varphi _2(m_2)dm_2\, dT(m_1)\\
 &=&\langle  T,\varphi _1*\check \varphi _2\rangle  \cr &=&\langle  T,\varphi _2^\tau *\varphi _1^\sigma \rangle  \\  
 \end{eqnarray*}
The last equality is due to the invariance of $T$ under the anti automorphism $\sigma $. Working back we get
\[ \langle  T,\varphi _2^\tau *\varphi _1^\sigma \rangle  =\int _{M\times M}\varphi _2^\tau (m_1)\varphi _1^\tau (m_2) dS(m_1,m_2)\]
As $(\varphi _i^\tau )^\sharp=(\varphi _i^\sharp )^\tau $ this implies that
\[ \langle  Af_1,A'f_2\rangle  =\langle  Af_2^\tau ,A'f_1^\tau \rangle  \] 
Suppose that we can find an $A'\ne 0$. Then $A'$ is onto and $Af_1=0$ if and only if, for all $f_2$  we have $\langle  Af_1,A'f_2\rangle  =0$. This means that $A'f_1^\tau $ is orthogonal to the image of $A$. If $A\ne 0$ this means that $A'f_1^\tau =0$. Therefore if both $A$ and $A'$ are non zeros then the kernel of $A$ is the image under $\tau $ of the kernel of $A'$. Because of Schur Lemma $A$ and $A'$ are, up to a constant factor determined by their respective kernels. So we are left with two possibilities for the dimensions ${\rm dim} \left  ({\mathrm {Hom}} _ M(ind_ N^M(1),\pi )\right  )$ and ${\rm dim} \left  ({\mathrm {Hom}} _ M(ind_ N^M(1),\pi ^*)\right  )$: either they are both equal to 1 or one of them is 0 and the other arbitrary ( possibly infinite).

{\bf Remark\pointir} There is a variant for the non unimodular case; we will not need it.

\begin{corollary}
 Let  $M$ be a tdlc group and $N$ a closed subgroup, both unimodular. Suppose that there exists an involutive anti-automorphism $\sigma $ of $M$ such that $\sigma (N)=N$ and such that any distribution on $M$, invariant under the adjoint action of  $N$, is fixed by $\sigma $. Then, for any irreducible admissible representation $\pi $ of $M$ and any irreducible admissible representation $\rho $ of $N$
\[ \dim\left  ({\mathrm {Hom}}_ N(\pi _{|N},\rho ^*)\right  )\times \dim\left  ({\mathrm {Hom}}_ N((\pi ^*)_{|N},\rho )\right  )\leq 1\]
\end{corollary}
Let $M'=M\times N$ and $N'$ the closed subgroup of $M'$ image of the homomorpism $n\mapsto (n,n)$  of $N$ into $M$. The map $(m,n)\mapsto mn^{-1}$ of $M'$ onto $M$ defines an homeomorphism of $M'/N'$ onto $M$.The inverse map is $m\mapsto (m,1)N'$ On $M'/N'$ left translations by $N'$ correspond to the adjoint action of $N$ onto $M$. We have a bijection between the space of distributions $T$ on $M$ invariant under the adjoint action of $N$ and the space of distributions $S$ on $M'$ which are biinvariant under $N'$. Explicitly
\[ \langle  S,f(m,n)\rangle  =\langle  T,\int _Nf(mn,n)dn\rangle  \]
Suppose that $T$ is invariant under $\sigma $ and consider the involutive anti-automorphism $\sigma '$ of $M'$ given by $\sigma '(m,n)=(\sigma (m),\sigma (n))$. Then
\[ \langle  S,f\circ \sigma '\rangle  =\langle  T,\int _Nf(\sigma (n)\sigma (m),\sigma (n))dn\rangle  \]
Using the invariance under $\sigma $ and for the adjoint action of $N$ we get
 \begin{eqnarray*}
 \langle  T,\int _Nf(\sigma (n)\sigma (m),\sigma (n))dn\rangle  &=&\langle  T,\int _Nf(\sigma (n)m,\sigma (n))dn\rangle  \\
&=&\langle  T,\int _Nf(mn,n)dn\rangle  \\
&=&\langle  S,f\rangle  
\end{eqnarray*}
\overfullrule=0pt
Hence $S$ is invariant under $\sigma '$. Conversely if $S$ is invariant under $\sigma '$ the same computation shows that $T$ is invariant under $\sigma $. Under the assumption of Corollary 2-1
 we can now apply Proposition 2-1 and we obtain the inequality
 \[ {\mathrm {dim}} \left  ({\mathrm {Hom}} _ {M'}(ind_ {N'}^{M'}(1),\pi \otimes \rho )\right  )\times {\mathrm {dim}} \left  ({\mathrm {Hom}} _ {M'}(ind_ {N'}^{M'}(1),\pi ^*\otimes \rho ^*)\right  )\leq 1\]
 We know that $Ind_ {N'}^{M'}(1)$ is the smooth contragredient representation of $ind_{N'}^{M'}(1)$hence
 \[ {\mathrm {Hom}} _ {M'}(ind_ {N'}^{M'}(1),\pi ^*\otimes \rho ^*)\approx {\mathrm{ Hom}} _ {M'}(\pi \otimes \rho , Ind _ {N'}^{M'}(1))\] 
 Frobenius reciprocity tells us that
 \[ {\mathrm {Hom}} _ {M'}(\pi \otimes \rho , Ind _ {N'}^{M'}(1))\approx {\mathrm {Hom}} _ {N'}(\pi \otimes \rho )_{|N'},1)\] 
Clearly
\[ {\mathrm {Hom}} _ {N'}(\pi \otimes \rho )_{|N'},1)\approx {\mathrm {Hom}} _ {N}(\rho ,(\pi _ {|N})^*)\approx {\mathrm {Hom}}_ N(\pi _{|N},\rho ^*)\]
Using again Frobenius reciprocity we get
\[ {\mathrm {Hom}} _ {N}(\rho ,(\pi _ {|N})^*)\approx {\mathrm {Hom}} _ {M}(ind_ N^M(\rho ),\pi ^*)\] 
In the above computations we may replace $\rho $ by $\rho ^*$ and $\pi $ by $\pi ^*$. Finally
\begin{eqnarray*}
{\mathrm {Hom}} _ {M'}(ind_ {N'}^{M'}(1),\pi ^*\otimes \rho ^*)&\approx& {\mathrm {Hom}} _ {N}(\rho ,(\pi _ {|N})^*)\\
&\approx& {\mathrm {Hom}}_ N(\pi _{|N},\rho ^*)\\
&\approx& {\mathrm {Hom}} _ {M}(ind_ N^M(\rho ),\pi ^*)\\
{\mathrm {Hom}} _ {M'}(ind_ {N'}^{M'}(1),\pi \otimes \rho )
&\approx&{\mathrm {Hom}} _ {N}(\rho ^*,((\pi ^*)_ {|N})^*)\\
&\approx& {\mathrm {Hom}}_ N((\pi ^*)_{|N},\rho )\\
&\approx& {\mathrm {Hom}} _ {M}(ind_ N^M(\rho ^*),\pi )
\end{eqnarray*}

Going back to the situation of   section 1 we take $M=GL(W)$ and $N=GL(V)$. let $E_\pi $ be the spaces of the representation $\pi $  and let $E_\pi ^*$ be the smooth  dual (relative to the action of $GL(W)$. Let $E_\rho $ be the space of $\rho $ and $E^*_\rho $ be the smooth dual for the action of $GL(V)$. We know that the contragredient representation $\pi ^*$ in $E_\pi ^*$ is equivalent to the representation $g\mapsto \pi (^tg^{{-1}})$ in $E_\pi $. The same is true for $\rho ^*$. Therefore an element of ${\mathrm {Hom}}_ N(\pi _{|N},\rho ^*)$ may be described as a linear map $A$ from $E_\pi $ into $E_\rho $ such that, for $g\in N$
\[ A\pi (g)=\rho (^tg^{-1})A\]
An element of ${\mathrm {Hom}}_ N((\pi ^*)_{|N},\rho )$ may be described as a linear map $A'$ from $E_\pi $ into $E_\rho $ such that, for $g\in N$
\[ A'\pi (^tg^{-1})=\rho (g)A'\] 
We have obtained the same set of linear maps:
\[ {\mathrm {Hom}}_ N((\pi ^*)_{|N},\rho )\approx {\mathrm {Hom}}_ N(\pi _{|N},\rho ^*)\] 
We are left with 2 possibilities: either both spaces have dimension 0 or they both have dimension 1 which is exactly what we want.

From now on we forget Conjecture 1 and try to prove Conjecture  2

\begin{section}
 {\hglue -12pt . A partial linearization}
\end{section}

We keep the notations of section 1. Put $G=GL(V)$. Let $V^*$ be the dual space of $V$. Choose a bijective linear map $u$ from $V$ to $V^*$ which is self-adjoint. On $G\times V\times V^*$ we consider the involution $\sigma $:
\[ \sigma (g,v,v^*)=(u^{-1}\,{}^tg\, u,u^{-1}(v^*),u(v))\]
The group $G$ acts on $G\times V\times V^*$ by
\[ g(x,v,v^*)=(gxg^{-1},gv\, ,^t\hskip -0.1cm  g^{-1}v^*)\]
We are going to show that Conjecture 2 is in turn implied by Conjecture 3

{\bf Remark\pointir}1) The involution $\sigma $ depends on the choice of $u$. However the action on invariant distributions does not so that the choice is irrelevant. If we need a precise $u$ we start with a basis  $e_1,\dots ,e_n$ of $V$ and the dual basis $e_1^*,\dots ,e_n^*$ and define $u$ by $u(e_i)=e_i^*$. 

Assume Conjecture 3. We shall apply it to $GL(W)$ acting on $GL(W)\times W\times W^*$. Let $e_1,\dots ,e_n$ be a basis of $V$ as above and choose $e_{n+1}\in U$,a non zero vector. Call $e_1^*,\dots e_{n+1}^*$ the dual basis of $W^*$. Define $u$, hence $\sigma $, using this basis. Let $Y$ be the set of all triples $(g,w,w^*)\in GL(W)\times W\times W^*$ 
 such that $\langle  w^*,w\rangle  =1$; it is a closed subset, stable under $\sigma $ and under $GL(W)$. Therefore, by Conjecture 3, any $GL(W)-$invariant distribution on $Y$ is fixed by $\sigma $. Now let $X$ be the closed subset of $Y$ defined by $w=e_{n+1}$ and $w^*=e^*_{n+1}$. Clearly $GL(W)X=Y$ and if $\xi \in X$ and $g\in GL(W)$ are such that $g\xi \in X$ then $g\in GL(V)$ and $gX=X$
and any $g\in GL(V)$ has this property. This means that we can use Frobenius descent as described in the Appendix : there is a one to one correspondence between $GL(W)-$invariant distributions on $Y$ and $GL(V)-$invariant distributions on $X$. If we identify $X$ with $GL(W)$ and if we choose an invariant measure on $GL(W)/GL(V)$ then, to a $GL(V)-$invariant
 distribution $S$ on $GL(W)$ corresponds the distribution $T$ on $Y$ given by
 \[ \langle  T,f\rangle  =\int _{GL(W)/GL(V)}\langle  S,f(gxg^{-1},ge_{n+1},^t\hskip -0.1cm g^{-1}e^*_{n+1}\rangle  dg\] 
 If we replace $S$ by $\sigma (S)$ we get the distribution $T'$:
\[ \langle  T',f\rangle  =\int _{GL(W)/GL(V)}\langle  S,f(gu^{-1}{\,\,}^t\hskip -0.1 cm x\,ug^{-1},ge_{n+1},^t\hskip -0.1cm g^{-1}e^*_{n+1}\rangle  dg\]  
Now $g\mapsto u^{-1}{\,\,}^t\hskip -0.1 cm g^{-1}\, u$ is an involutive automorphisme of $GL(W)$ which defines an involutive homeomorphism of $GL(W)/GL(V)$. Changing variables :
\[ \langle  T',f\rangle  =\int _{GL(W)/GL(V)}f(u^{-1}\,\,{}^t(gxg^{-1})u,u^{-1}(^tg^{-1}e^*_{n+1}),u(ge_{n+1}))\rangle  dg\] 
which is equal to $\langle  T,\sigma (f)\rangle  $ so that $T'=\sigma (T)$. We assumed that $\sigma (T)=T$ and the map $S\mapsto T$ is one to one so $\sigma (S)=S$.

We let $G=GL(V)$ acts on its Lie algebra by the adjoint action and consider the action on ${\mathfrak {g} }\times V\times V^*$. Clearly the involution $\sigma $ acts also on ${\mathfrak {g} }$.

We are going to show that Conjecture 4 implies Conjecture 3. The basic tool is Harish-Chandra's  descent method.   
The proof is by induction on the dimension $n$ of $V$. So let us assume that Conjectures 3 and 4 are true for all local fields ${\mathbb F}$ non archimedean of characteristic 0 and $\dim V<n$. We take $V$ of dimension $n$.

\vskip 5ex
\begin{section}{\hglue -12pt . The case $\mathbf n=1$}
\end{section}
We take $V={\mathbb F}$ and also $V^*={\mathbb F}$, the duality being the usual product. The group $GL(V) $ is ${\mathbb F}^*$. It acts trivially on itself and on $V\times V^*\approx {\mathbb F}^2$ by $(x,y)\mapsto (tx,t^{-1}y)$.The involution $\sigma $ is $(t,x,y)\mapsto (t,y,x)$. We have to show that any distribution on ${\mathbb F}^2$, invariant under the action of ${\mathbb F}^*$ is invariant by the map $(x,y)\mapsto (y,x)$. 
Let $\Omega \subset{\mathbb F}^2$ be the open subset $xy\ne 0$. The fibers of the map $(x,y)\mapsto xy$ of $\Omega $ onto ${\mathbb F}^*$ are the orbits of ${\mathbb F}^*$. Each of them is fixed under $\sigma $. It follows from Bernstein's localization principle  that a distribution on $\Omega $, invariant under ${\mathbb F}^*$ is invariant under $\sigma $. If $T$ is a distribution on ${\mathbb F}^2$, invariant under ${\mathbb F}^*$ and such that $\sigma (T)=-T$ its support must be contained in the closed subset $Z:\,\, xy=0$. Then there exist two constants $a$ and $b$ such that, on $Z\setminus \{(0,0)\}$
\[ \langle  T,f\rangle  =a\int _{{\mathbb F}^*}f(x,0)d^*x+b\int _{{\mathbb F}^*}f(0,y)d^*y\] 
One possible extension of the right hand side to a distribution on $Z$ is
\[ \langle  T',f\rangle  =a\int _{{\mathbb F}^*}(f(x,0)-f(0,0)\chi (x))d^*x+b\int _{{\mathbb F}^*}(f(0,y)-f(0,0)\chi (y))d^*y\]
where $\chi$ is the charactereistic function of the ring of integers in $\mathbb F$.
Then $T-T'$ is a multiple of the Dirac measure at the origin hence is invariant under $\sigma $ and ${\mathbb F}^*$. On the other hand  a trivial computation shows that $T'$ is invariant by ${\mathbb F}^*$ if and only if $a=b$. On $Z\setminus \{(0,0)\}$ the distributions $T$ and $T'$ are equal so $T'$ is skew with respect to $\sigma $ which means that $a=-b$. Thus we must have $a=b=0$ and , as there is no skew distributions supported by the origin we conclude that $T=0$. 
In other words the space of ${\mathbb F}^*-$invariant distributions on $Z$ has dimension 2 and a basis consists of the Dirac measure at the origin and the distribution
\[ f\mapsto \int _{{\mathbb F}^*}(f(x,0)-f(0,0)\chi (x))d^*x+\int _{{\mathbb F}^*}(f(0,y)-f(0,0)\chi (y))d^*y\] 
Invariance under ${\mathbb F}^*$ forces symmetry under $\sigma $!
\vskip 5ex
\begin{section}
{\hglue -12pt . Harish-Chandra's descent}
\end{section}
Just for this section we work in a more general setup.
Let ${\mathbb F}$ be a local field, of characteristic 0 and non archimedean. Let ${\mathfrak {g}}$ be a reductive
Lie algebra, defined over ${\mathbb F}$ and $G$ an algebraic group, defined over ${\mathbb F}$ and of Lie
algebra ${\mathfrak { g}}$. 

For $x\in {\mathfrak  {g}}$ we denote by ${\mathcal {Z}}_{\mathfrak {g} }(x)$ the centralizer of $x$ in ${\mathfrak {g} }$ and by
$Z_ G(x)$ its centralizer in $G$.

The group $G$ acts on ${\mathfrak {g} }$ by the adjoint action; we write $gx$ instead of $\mathrm {Ad} g(x)$. 

Let $a$ be a non central semi-simple element of ${\mathfrak {g} }$. Put
\[ {\mathfrak {m}   ={\mathcal Z}_{\mathfrak {g} }}(a),\qquad M=Z_ G(a)\] 
Let ${\mathfrak {q}  }$ be the image of $\mathrm {ad}      (a)$; we have ${\mathfrak {g} }={\mathfrak {m}   }\oplus   {\mathfrak {q}  }$. For $x\in {\mathfrak {m}   }$
\[ \mathrm {ad}      (x){\mathfrak {m}   }\subset{\mathfrak {m}   },\qquad \mathrm {ad}      (x){\mathfrak {q}  }\subset{\mathfrak {q}  } \] 
Let ${\mathfrak {m}   }'$ be the set of all $x\in {\mathfrak {m}   }$ such that $(\mathrm {ad}\,      x)_{|{\mathfrak {q}   }}$ is invertible.
It is a Zariski open set in ${\mathfrak {m}   }$.
Choose a complete set ${\mathfrak {h}    }_1,\dots ,{\mathfrak {h}    }_r$ of representative of conjugacy classes of Cartan
subalgebras of ${\mathfrak {m}   }$ (the group which acts is thus $M$). If $G$ acts continuously on some topological
space $T$, then a $G-$domain in $T$ is a subset of $T$ which is closed, open and invariant under $G$. The
following Lemma is due to Harish-Chandra \cite{4}, Corollary 2-3 
\begin{lemma}
 Let $\omega $ be a neighbourhood of $a$ in ${\mathfrak {m}   }$. There exists an $M-$domain $\Omega $ in ${\mathfrak {m}   }$ such that\\
 {a)} $a\in \Omega \subset M\omega $,\\
 {b)}  $\Omega {\mathfrak {h}    }_i\subset\omega $ for  $i=1,\dots ,r$,\\
 {c)}$\,\, \Omega \subset{\mathfrak {m}   }'$,\\
 {d)} for all compact sets $Q$ of ${\mathfrak {g} }$ there exists a compact set $C$ in $G$ such that $x\Omega Q\ne \emptyset  $
implies $x\in C$
\end{lemma}

Fix such an $\Omega $. Consider the map $(g,x)\mapsto gx$ from $G\times \Omega $ into ${\mathfrak {g} }$.  The differential of
this map is
\[ (X,Y)\mapsto g([X,x]+Y),\qquad X\in {\mathfrak {g} },\,\, Y\in {\mathfrak {m}   }\]
It is onto if and only if ${\mathfrak {m}   }+{\rm Im}(\mathrm {ad} \,     x)={\mathfrak {g} }$ which means that $x\in {\mathfrak {m}   }'$. As
we know that $\Omega \subset{\mathfrak {m}   }'$ the map is everywhere submersive. In particular its image $G\Omega $ is
open. It follows from  condition d) of Lemma 5-1 that $G\Omega $ is also closed \cite{4}, Corollary 2-4.

Let $E$ be a finite dimensional rational $G-$module (defined over ${\mathbb F}$). The map $(g,x,e)\mapsto (gx,ge)$ from $G\times \Omega \times E$ into ${\mathfrak {g} }\times E$ is everywhere submersive so that we may apply Harish-Chandra submersion principle. There exists a map $f\mapsto F_f$ of ${\mathcal S}(G\times \Omega \times E)$ onto ${\mathcal S}(\mathrm {Ad}  (G)\Omega \times E)$ such that for any $\varphi \in {\mathcal S}(\mathrm {Ad}  (G)\Omega \times E)$
\[ \int _{G\times \Omega \times E}f(g,x,e)\varphi (gxg^{-1},ge)dg\, dx\, de=\int _{\mathrm {Ad}  (G)\Omega \times E}F_f(X,e)\varphi (X,e) dX\, de\] 
Here $dg$ is a Haar measure of $G$, $dx$ the restriction to $\Omega $ of a Haar measure of ${\mathfrak {m}   }$, $de$ a Haar measure of $E$ and $dX$ the restriction to $\mathrm {Ad}  (G)\Omega $ of a Haar measure of ${\mathfrak {g} }$. By transposition we then get a one to one map $T\mapsto \Phi  (T)$ of ${\mathcal S}'(\mathrm {Ad}  (G)\Omega \times E)$ into ${\mathcal S}'(G\times \Omega \times E)$. We let $G$ acts on $G\times \Omega \times E$ by left translations on the first factor and on $\mathrm {Ad}  (G)\Omega \times E$ by adjoint on the first factor and the given action on the second factor so that the map $f\mapsto F_f$ commutes with the action of $G$. It follows that if $T$ is invariant then $\Phi  (T)$ is also invariant and therefore may be written in a unique way as $\Phi  (T)=dg\otimes \theta (T)$ with $\theta (T)\in {\mathcal S}'(\Omega \times E)$. If we define
\[ H_f(x,e)=\int _Gf(g,x,e)dg\] 
then $\langle  T,F_f\rangle  =\langle  \theta (T),H_f\rangle  $. The map $\theta $ is one to one. Furthermore the group $M$ acts on $\Omega $ by the adjoint action and  it acts on $E$. The distribution $\theta (T)$ is semi-invariant. Indeed let $\psi \in {\mathcal S}(\Omega \times E)$ and choose $f$ such that $\psi =H_f$. Fix $m\in M$ and define $\psi _1(x,e)=\psi (mxm^{-1},me)$. If $f_1$ is given by $f_1(g,x,e)=f(g,mxm^{-1},me)$ then $\psi _1=H_{f_1}$. To compute $F_{f_1}$ we start from
\[ \int _{G\times \Omega \times E}f(g,mxm^{-1}x,me)\varphi (gxg^{-1},ge)dg\, dx\, de=\int _{\mathrm {Ad}  (G)\Omega \times E}F_{f_1}(X,e)\varphi (X,e) dX\, de\]
and on the left-hand side we change $m$ into $m^{-1}xm$ and $e$ into $m^{-1}e$ which gives
\[ \int _{G\times \Omega \times E}f(g,x,e)\varphi (gm^{-1}xmg^{-1},gm^{-1}e)dg\, dx\, de\, |\mathrm {Det} _ Em|^{-1}\] 
and now we change $g$ into $gm$:
\[ \int _{G\times \Omega \times E}f(gm,x,e)\varphi (gxg^{-1},ge)dg\, dx\, de\, |\mathrm {Det} _ Em|^{-1}\]
Define $f_2(g,x,m)=f(gm,x,e)$. We have $F_{f_2}=|\mathrm {Det} _ E m| F_{f_1}$. Now
 \begin{eqnarray*}
 \langle  \theta (T),\psi _1\rangle  &=\langle  T,F_{f_1}\rangle  =\langle  T,F_{f_2}\rangle  \mathrm {Det} _ E m|^{-1}=\langle  \Phi  (T),f_2\rangle  |\mathrm {Det} _ E m|^{-1}\\
  &=\langle  \Phi  (T),f_1\rangle | \mathrm {Det} _ E m|^{-1}=\langle  T,\psi _1\rangle | \mathrm {Det} _ E m|^{-1}
  \end{eqnarray*}

\begin{proposition} 
 The map $\theta $ is a one to one map from the space   ${\mathcal S}(G\Omega \times E)^{'G}$  of $G-$invariant distributions on $G\Omega \times E$ into the space of distributions on $\Omega \times E$ having the above semi-invariance. 
\end{proposition}

Similar results are true for the group. Let $a$ be a semi-simple element of $G$; we do not exclude the case $a$ central. Let $M$ be the centralizer of $a$ in $G$ and ${\mathfrak {m}   }$ its centralizer in ${\mathfrak {g} }$. The endomorphism ${\rm Id}-Ad(a)$ of ${\mathfrak {g} }$ is semi-simple so that ${\mathfrak {g} }={\mathfrak {m}   }\oplus   {\mathfrak {q}  }$ where ${\mathfrak {q}  }$ is the image of ${\rm Id}-Ad(a)$. For any $m\in M$ one has ${\rm Id}- \mathrm {Ad}  m({\mathfrak {q}  }\subset{\mathfrak {q}  }$  and if $m$ is sufficiently close to  the neutral element of $M$ then the restriction to ${\mathfrak {q}  }$ of ${\rm Id}-\mathrm {Ad}  m$ is bijective.  Choose a linearization $G\subset GL_ N({\mathbb F})$ of $G$; then ${\mathfrak {g} }\subset{\mathfrak {g} }{\mathfrak  l}_ N({\mathbb F})$ and we have an exponential map defined on a subset of ${\mathfrak {g} }$.
\begin{lemma}
 There exists a neighborhood  $\omega $ of 0 in ${\mathfrak {m}   }$ with the following properties:\hfill\break
a) $\omega $ is open, closed and invariant under the action of $M$\hfill\break
b) the exponential map is defined on $\omega $ and submersive at each point of $\omega $\hfill\break
c) for any $X\in \omega $ the restriction to ${\mathfrak {q}  }$ of ${\rm Id}-\mathrm {Ad}  (ae^X)$ is bijective.
\end{lemma}
Indeed there exists in ${\mathfrak {g} }{\mathfrak  l}_ N({\mathbb F})$ a $GL_ N({\mathbb F)}-$invariant open neighborhood $U$ of 0 such that the exponential map is defined on $U$ and is an analytic isomorphism of $U$ onto $\mathrm {Exp}   (U)$. Choose an open neighborhood $\omega _0$ of 0 in ${\mathfrak {m}   }$ such that $\omega _0\subset{\mathfrak {m}   }\cap U$ and such that condition c) of the Lemma is satisfied for $X\in \omega _0$ hence also for $X\in \mathrm {Ad}  (M)\omega _0$. Then we take for $\omega $ an open, closed and $M-$invariant neighborhood of 0 such that $\omega\subset\mathrm {Ad}  (M)\omega _0$ ( see Lemma  5-1). Let $\Omega =\mathrm {Exp}   (\omega )$; it is an open invariant neighborhood of the identity in $M$ ( be careful that  $\Omega $ does not have the same meaning as before!).

Let $E$ be a finite dimensional rational $G-$module.
\begin{lemma}
 The map $(g,m,e)\mapsto (gmag^{-1},ge)$ of $G\times \Omega \times E$ into $G\times E$ is everywhere submersive.
\end{lemma}
The differential is a map from ${\mathfrak {g} }\oplus   {\mathfrak {m}   }\oplus   E$ into ${\mathfrak {g} }\oplus   E$. To prove the surjectivity we compute the partial differentials. First for $(g,e)$, let $(X,Z)\in {\mathfrak {g} }\oplus   E$ then
\begin{eqnarray*}
g\mathrm {Exp} ( {tX})(e+uZ)&=&g(e+tX+uZ+\cdots )\\
 g\mathrm {Exp}  (tX)ma\mathrm {Exp}  (-tX)g^{-1}&=&gmag^{{-1}}\mathrm {Exp}  \left  (t\mathrm {Ad} (g)\mathrm {Ad} (m^{-1}a^{-1})({\rm Id}-\mathrm {Ad} (ma))X+\cdots \right  )
 \end{eqnarray*} 
Then the image of this partial differential is 
\[ \left  (\mathrm {Ad} (gm^{-1}a^{-1}){\rm Im}({\rm Id}-\mathrm {Ad} (ma))\right  )\oplus   E\] 
The image of ${\mathrm {Id}}-\mathrm {Ad}  (ma)$ contains ${\mathfrak {q}  }$ and ${\mathfrak {q}  }$ is fixed by $\mathrm {Ad}  (m^{{-1}}a^{-1})$ so the image of he differential contains the subspace $\mathrm {Ad}  (g)({\mathfrak {q}  })\oplus   E$.
  
Next we take $Y\in {\mathfrak {m}   }$, then
\[ gm\exp(tY)ag^{-1}=gmag^{-1}\mathrm {Exp}  (t\mathrm {Ad}  (g)Y)\]
and so the image of the differential contains $\mathrm {Ad}  (g){\mathfrak {m}   }$. Thus the differential is onto, as required.

Now we apply Harish-Chandra submersion principle: there exists a surjective linear map $f\mapsto F_f$ of ${\mathcal S}(G\times \Omega \times E)$ onto ${\mathcal S}(\mathrm {Ad}  (G)(\Omega a)\times E)$ such that, for any $\varphi \in {\mathcal S}(\mathrm {Ad}  (G)(\Omega a)\times E)$
\[ \int _{G\times \Omega \times E}f(g,m,e)\varphi (gmag^{-1},ge)dg\, dm\, de=\int _{\mathrm {Ad} (G)(\Omega a)\times E}F_f(x,e)\varphi (x,e)dx\, de\] 
On the left hand side $dm $ is a Har measure on $M$, $dg$ a Haar measure on $G$ and $de$ on $E$. On the right hand side $dx$ is the restriction of $dg$ to $\mathrm {Ad}  (G)\Omega $. By transposition we have a one to one linear map $T\mapsto \Phi  (T)$ of ${\mathcal S}'(\mathrm {Ad}  (G)(\Omega a)\times E)$ into ${\mathcal S}'(G\times \Omega \times E)$. The group $G$ acts on $\mathrm {Ad}  (G)(\Omega a)\times E$ by the adjoint action on $\mathrm {Ad}  (G)(\Omega a)$ and the given action on $E$ and we let it act on $G\times \Omega \times E$ by left translations on the second component. As the map $f\mapsto F_f$ commutes with these actions, if the distribution $T$ is invariant so is the distribution $\Phi  (T)$. In that case there is a unique distribution $\theta (T)\in {\mathcal S}'(\Omega \times E)$ such that $\Phi  (T)=dg\otimes \theta (T)$. If we define
\[ H_f(m,e)=\int _Gf(g,m,e)dg\] 
than $\langle  T,F_f\rangle  =\langle  \theta (T),H_f\rangle  $.

The group $M$ acts on $\Omega $ by adjoint action and so it acts on $\Omega \times E$. If $u\in M$ then
\[ H_f(umu^{-1},ue)=\int _Gf(g,umu^{-1},ue)dg\]
Define  $f'(g,m,e)=f(g,umu^{-1},ue)$ and $f''(g,m,e)=f(gu,m,e)$. Then
 \begin{eqnarray*}
 &\int _{G\times \Omega \times E}f(g,umu^{-1},ue)\varphi (gmag^{-1},ge) dgdmde\\
 & =|\mathrm {Det} _ E u|^{-1}\int _{G\times \Omega \times E}f(gu,m,e)\varphi (gmag^{-1},ge)dgdmde
\end{eqnarray*}

which means that $F_{f'}=|\mathrm {Det} _ E u|^{-1}F_{f''}$. Therefore $\langle  T,F_{f'}\rangle  =|\mathrm {Det} _ E u|^{-1}\langle  T,F_{f''}\rangle  $ and then $\langle  \Phi  (T),f'\rangle  =|\mathrm {Det} _ E u|^{-1}\langle  \Phi  (T),f''\rangle  $. However $\Phi  (T)=dg\otimes \theta (T)$ implies that $\langle  \Phi  (T),f''\rangle  =\langle  \Phi  (T),f\rangle  $ so that $\langle  \theta (T),H_{f'}\rangle  =|\mathrm {Det} _ E u|^{-1}\langle  \theta (T),f\rangle  $: the distribution $\theta (T)$ is semi-invariant under $M$; in integral notations:
\[ \int _{\Omega \times E}\psi (umu^{-1},ue)\, d\theta (T)(m,e)=|\mathrm {Det} _ E u|^{-1}\int _{\Omega \times E}\psi (m,e)\, d\theta (T)(m,e)\] 
\begin{proposition} 
 The map $\theta $ is a one to one linear map from the space ${\mathcal S }'(\mathrm {Ad} (G)(\Omega a)\times E)^G$ of $G-$invariant distributions on $\mathrm {Ad}  (G)(\Omega a)\times E$ into the space of distributions on $\Omega \times E$ with the above semi-invariance property.
\end{proposition}

Now we go back to our setup: $G=GL(V)$ and $E=V\oplus   V^*$. In this case $\mathrm{Det}_ E g=1$
Let $a$ be a  semi-simple element of $\mathrm{End}  \,_{\mathbb F}(V)$. We want to describe the centralizer
${\mathfrak {m}   }$ of $a$ in ${\mathfrak {g} }$ and the centralizer $M$ of $a$ in the linear group $G$.  Let $P$ be its minimal polynomial; all its roots are simple. Let $P=P_1\dots P_r$ be
the decomposition of $P$ into irreducible factors, over ${\mathbb F}$. Then the $P_i$
 are two by two relatively prime. If $V_i=\mathrm{Ker}    P_i(a)$, then $V=\oplus   V_i$ and $V^*=\oplus   V_i^*$. An element $x$ of $G$
which commutes with $a$ is given by  a familly $\{x_1,\dots ,x_r\}$ where each $x_i$ is a linear map from $V_i$
to $V_i$, commuting with the restriction of $a$ to $V_i$. Now ${\mathbb F}[T]$ acts on $V_i$, by specializing $T$ to
$a_{|V_i}$ and $P_i$ acts trivially so that, if ${\mathbb F}_i={\mathbb F}[T]/(P_i)$, then $V_i$ becomes a vector space
over ${\mathbb F}_i$. The ${\mathbb F}-$linear map $x_i$ commutes with $a$ if and only if it is ${\mathbb
F}_i-$linear.\hfill\break

Fix $i$. Let $\ell$ be a non zero ${\mathbb F}-$linear form on ${\mathbb F}_i$. If $v_i\in V_i$ and $v'_i\in V_i^*$ then
$\lambda \mapsto \langle  \lambda v_i,v'_i\rangle  $ is an ${\mathbb F}-$linear form on ${\mathbb F}_i$, hence there exists a unique element
$S(v_i,v'_i)$ of ${\mathbb F}_i$ such that $\langle  \lambda v_i,v'_i\rangle  =\ell\left  (\lambda S(v_i,v'_i)\right  )$. One checks trivially that $S$ is
${\mathbb F}_i-$linear with respect to each variable and defines a non degenerate duality, over ${\mathbb F}_i$ between
$V_i$ and $V_i^*$. Here ${\mathbb F}_i$ acts on $V_i^*$ by transposition, relative to the ${\mathbb F}-$duality $\langle  .,.\rangle  $ ,of
the action on $V_i$. Finally if $x_i\in \mathrm{End}  _{{\mathbb F}_i}V_i$ its transpose, relative to the duality $S(.,.)$ is the
same as its transpose relative to the duality $\langle  .,.\rangle  $. 

Thus $M$ is a product of linear groups and the situation $(M,V,V^*)$ is a composite case, each component being
a linear case (over various extensions of ${\mathbb F}$).

Recall that the involution on ${\mathfrak {g} }\times V\times V^*$ or on $G\times V\times V^*$is not  unique  (although the action on the space of invariant distributions is).We choose a basis of $V$ that is to say an isomorphism $\varphi $ of ${\mathbb F}^n$ onto $V$
onto $V$, identify ${\mathbb F}^n$ with its dual and then put $u=^t\hskip -0.6pt \varphi \circ\varphi $. The involution is then given by
\[ (x,v,v^*)\mapsto (u^{-1}\,\, ^tXu,u^{-1}(v^*),u(v))\] 
We assume that $\varphi $ is compatible with the decomposition $V=\oplus V_i$
in the sense that the image by $\varphi $ of the canonical basis of ${\mathbb F}^n$ is the union of basis of the $V_i$.
We are going to define $\sigma _a\in G$ such that $\sigma _a\, u^{-1}au\, \sigma _a^{-1}=a$. This can be done independantly for each $V_i$.

 Fix $i$. Let $s$ be the dimension of $V_i$ over ${\mathbb F}_i$ and choose an isomorphism $\Phi  $ of ${\mathbb F}_i^s$ onto $V_i$, then identify ${\mathbb F}_i^s$ with its dual and put $U_i=^t\hskip -0.6 pt \Phi  \Phi  $.
 
Let $U$ be the direct sum of the $U_i$. Note that $Ua={}^t\hskip -1.6pt a\,\, U^{-1}$. Call ${ s}$  the involution on ${\mathfrak {m}   }\times V\times V^*$ or on $M\times V\times V^*$
product of the involutions built with the $U_i$ on each factor.

Let us look at the situation on the group. 
Let $\psi \in {\mathcal S}(\Omega \times V\times V^*)$; choose $f\in {\mathcal S}(G\times \Omega \times V\times V^*)$ such that $\psi =H_f$. Let $\psi _1$ be defined by
\[ \psi _1(m,v,v^*)=\psi (U^{-1}{\,}^tm\, U,U^{-1}(v^*),U(v))\] 
This make sense because any $M-$orbit (resp. $G-$orbit) in ${ M}$ (resp. in ${ G}$) is stable by the involution $m\mapsto U^{-1}{\,}^tm\, U$ (resp .$g\mapsto u^{-1}{\,}^tg\, u$)

Then $\psi _1=H_{f_1}$ with 
\[ f_1(g,m,v,v^*)=f(g,U^{-1}{\,}^tm\, U,U^{-1}(v^*),U(v))\]
We want to compute $F_{f_1}$. We have
\begin{eqnarray*}
&\int _{G\times \Omega \times V\times V^*}f(g,U^{-1}{\,}^tm\, U,U^{-1}(v^*),U(v))\varphi (gmag^{-1},gv,^tg^{-1}v^*)dg\, dm\, dv\, dv^*\\
=&\int _{G\times \Omega \times V\times V^*}f(g,m,v,v^*)\varphi (gU^{-1}{\,}^tm\, Uag^{-1},gu^{-1}v^*,^tg^{-1}Uv)dg\, dm\, dv\, dv^*\\
=&\int _{G\times \Omega \times V\times V^*}f(u^{-1}{\,}^tg^{-1}\, U,x,v,v^*)\varphi (u^{-1}{\,}^t(gmag^{-1})u,u^{-1}{\,}^tg^{-1}v^*,guv)dg\, dx\, dv\, dv^*
\end{eqnarray*}
The last equality corresponding to the involutive change of variables $g\mapsto u^{-1}{\,}^tg^{-1}\, U$. Then define $f_2$ by
\[ f_2(g,m,v,v^*)=f(u^{-1}{\,}^tg^{-1}\, U,m,v,v^*).\]
The above integral is equal to
\begin{eqnarray*}
&\int _{\mathrm {Ad} (G)(\Omega a)\times V\times V^*}F_{f_2}(g,v,v^*)\varphi (u^{-1}{\,}^tg\, u,u{{-1}(v^*),u(v))}dg\, dv\, dv^*\\
=&\int _{\mathrm {Ad} (G)(\Omega a)\times V\times V^*}F_{f_2}(u^{-1}{\,}^tg\, u,u^{-1}(v^*),u(v))\varphi (g,v,v^*)dg\, dv\, dv^*
\end{eqnarray*} 
 Finally
 \[ F_{f_1}(g,v,v^*)=F_{f_2}(u^{-1}{\,}^tg\, u,u^{-1}(v^*),u(v))\] 
If $T$ is an invariant distribution  on $\mathrm {Ad}  (G)(\Omega a)\times V\times V^*$, skew symmetric with respect to the involution  $\sigma $ then 
\[ \langle  T,F_{f_1}\rangle  =-\langle  T,F_{f_2}\rangle  .\] 
Therefore
\[ \langle  \theta (T),\psi _1\rangle  =-\langle  \theta (T),H_{f_2}\rangle  .\]
However $H_{f_2}=H_f=\psi $ so the distribution $\theta (T)$ is skew symmetric with respect to the involution ${ s}$.

 The exponential map is an isomorphism of $\omega $ onto $\Omega $ and commutes with the adjoint action of $M$. We may view $\theta (T)$ as an invariant distribution on $\omega \times V\times V^* $ but $\omega $ is closed in ${\mathfrak {m}   }$ so $\theta (T)$ extends to  an invariant distribution on ${\mathfrak {m}   }\times V\times V^*$. If we assume Conjecture 4 then $\theta (T)$ must be symmetric, hence 0 and $T=0$.

Now let $T$ be an invariant distribution on $G\times V\times V^*$ which is skew with respect to $\sigma $ and let $(g,v,v^*)$ be any element of $G\times V\times V^*$. Take for $a$ the semi-simple part of $g$. Then , always assuming Conjecture 4, the distribution $T$ will be 0 on $\mathrm {Ad}  (G)(\Omega a)\times V\times V^*$. However $a$ belongs to the closure of the orbit of $g$ so, for some $g'\in G$ we have $(g'gg^{'-1},g'v,^tg^{'-1}v^*)\in \mathrm {Ad}  (G)(\Omega a)\times V\times V^*$ and we conclude that $T=0$ in a neighborhood of $(g,v,v^*)$ and finally $T=0$.

{\noindent\bf Conclusion\pointir Conjecture 4 implies Conjecture 3 }

We now deal, in exactly the same way,  with the Lie algebra case. We assume that $a$ is not central so that, by the induction hypothesis, Conjecture 4 is valid for each $i$ and hence for the composite case ${\mathfrak {m}   }\oplus   V\oplus   V^*$. We use the notations of Proposition 5-2.

Let $\psi \in {\mathcal S}(\Omega \times V\times V^*)$; choose $f\in {\mathcal S}(G\times \Omega \times V\times V^*)$ such that $\psi =H_f$. Let $\psi _1$ be defined by
\[ \psi _1(x,v,v^*)=\psi (U^{-1}{\,}^tx\, U,U^{-1}(v^*),U(v))\]
This make sense because any $M-$orbit (resp. $G-$orbit) in ${\mathfrak {m}   }$ (resp. in ${\mathfrak {g} }$) is stable by the involution $x\mapsto U^{-1}{\,}^tx\, U$ (resp .$X\mapsto u^{-1}{\,}^tX\, u)$
Then $\psi _1=H_{f_1}$ with 
\[ f_1(g,x,v,v^*)=f(g,U^{-1}{\,}^tx\, U,U^{-1}(v^*),U(v))\] 
We want to compute $F_{f_1}$. We have
\begin{eqnarray*}
&\int _{G\times \Omega \times V\times V^*}f(g,U^{-1}{\,}^tx\, U,U^{-1}(v^*),U(v))\varphi (gxg^{-1},gv,^tg^{-1}v^*)dg\, dx\, dv\, dv^*\\
=&\int _{G\times \Omega \times V\times V^*}f(g,x,v,v^*)\varphi (gU^{-1}{\,}^tx\, Ug^{-1},gu^{-1}v^*,^tg^{-1}Uv)dg\, dx\, dv\, dv^*\\
=&\int _{G\times \Omega \times V\times V^*}f(u^{-1}{\,}^tg^{-1}\, U,x,v,v^*)\varphi (u^{-1}{\,}^t(gxg^{-1})u,u^{-1}{\,}^tg^{-1}v^*,guv)dg\, dx\, dv\, dv^*
\end{eqnarray*}
 
The last equality corresponding to the involutive change of variables $g\mapsto u^{-1}{\,}^tg^{-1}\, U$. Then define $f_2$ by
\[ f_2(g,x,v,v^*)=f(u^{-1}{\,}^tg^{-1}\, U,x,v,v^*).\] 
The above integral is equal to
\begin{eqnarray*}
&\int _{\mathrm {Ad} (G)\Omega \times V\times V^*}F_{f_2}(X,v,v^*)\varphi (u^{-1}{\,}^tX\, u,u{{-1}(v^*),u(v))}dX\, dv\, dv^*\\
=&\int _{\mathrm {Ad} (G)\Omega \times V\times V^*}F_{f_2}(u^{-1}{\,}^tX\, u,u^{-1}(v^*),u(v))\varphi (X,v,v^*)dX\, dv\, dv^*
\end{eqnarray*}
 Finally
 \[ F_{f_1}(X,v,v^*)=F_{f_2}(u^{-1}{\,}^tX\, u,u^{-1}(v^*),u(v))\]
If $T$ is an invariant distribution  on $\mathrm {Ad}  (G)\Omega \times V\times V^*$, skew symmetric with respect to the involution  $\sigma $ then 
\[ \langle  T,F_{f_1}\rangle  =-\langle  T,F_{f_2}\rangle  .\] 
Therefore
\[ \langle  \theta (T),\psi _1\rangle  =-\langle  \theta (T),H_{f_2}\rangle  .\]
However $H_{f_2}=H_f=\psi $ so the distribution $\theta (T)$ is skew symmetric with respect to the involution ${ s}$, but $\Omega \times V\times V^*$ is open and close so that $\theta (T)$ may be considered as a distribution on ${\mathfrak {m}   }\times V\times V^*$ and, by induction, is symmetric so it must be 0 and $T=0$.

Now $\mathrm {Ad}  (G)\Omega \times V\times V^*$ is open in ${\mathfrak {g} }\times V\times V^*$ so if $T$ is now an invariant distribution on ${\mathfrak {g} }\times V\times V^*$, skew symmetric with respect to $\sigma $ we conclude that its restriction to $\mathrm {Ad}  (G)\Omega \times V\times V^*$ is 0. If $(X,v,v^*)$ belongs to the support of such a distribution and if the semi-simple part $a$ of $X$ is not central then, with the above notations, $T$ is 0 on $\mathrm {Ad}  (G)\Omega \times V\times V^*$.  We know that $a$ belongs to the closure of the orbit of $X$, hence there exists $g\in G$ such that $(gXg^{-1},gv,{}^tg^{-1}v^*)\in \mathrm {Ad}  (G)\Omega \times V\times V^*$. It follows that the orbit of $(X,v,v^*)$ is contained in $\mathrm {Ad}  (G)\Omega \times V\times V^*$ and so $T$ is 0 in a neighborhood of $(X,v,v^*)$. Thus:

{\noindent\bf Conclusion\pointir If $T$ is a skew symmetric invariant distribution on ${\mathfrak {g} }\times V\times V^*$ then its support is contained in ${\mathfrak { z}}\times {\mathcal N}\times V\times V^*$ where ${\mathfrak  z}$ is the center of ${\mathfrak {g} }$ and ${\mathcal N}$ the c\^one of nilpotent elements in $[{\mathfrak {g} },{\mathfrak {g} }]$.}

\begin{section}
{\hglue -12pt . Regular orbits}
\end{section}
For $v\in V$,$v^*\in V^*$ and $X\in {\mathfrak {g} }$ put
\[ q_i(X,v,v^*)=\langle  v^*,X^iv\rangle  =\langle  {}^tX^iv^*,v\rangle  \]
and
\[ \mathrm {Det}  (T\,{\rm Id} -X)=T^n-\sum _0^{n-1}D_j(X)T^j\] 
The $q_i$ and $D_i$ are invariant polynomial functions. Let
\[ A(X,v,v^*)=(a_{i,j}),\quad a_{i,j}=\langle  {}^tX^{j-1}v^*,X^{i-1}v\rangle  ,\,\, i,j=1,\dots n\] 
an $n\times n$ matrix and put
\[ D(X,v,v^*)=\mathrm {Det} (A(X,v,v^*))\] 

A triplet $(X,v,v^*)$ is {\bf regular} if $\{ v,Xv,\dots X^{n-1}v\} $ is a basis of $V$ and
$\{v^*,{}^tXv^*,\dots ,{}^tX^{n-1}v^*\} $ is a basis of $V^*.$
\begin{proposition}
 The triplet $(X,v,v^*)$ is regular if and only if $D(X,v,v^*)\ne 0$. The set of
regular elements is a non empty Zariski open subset.
\end{proposition}
On $V\oplus   V^*$ the bilinear form
\[ \left  ((v,v^*),(w,w^*)\right  )\mapsto \langle  v^*,w\rangle  +\langle  w^*,v\rangle  \]
is symmetric and non degenerate. If $(X,v,v^*)$ is regular then 
\[v,\dots
,X^{n-1}v,v^*,\dots ,{}^tX^{n-1}v^*\]
 is a basis of $V\oplus   V^*$ and, relative to this basis, the
matrix of the above bilinear form is
\[ \pmatrix{0&A\cr A&0
\cr}\] 
and is non singular which means that $A$ is non singular. Conversely, any  non
trivial linear relation among the $X^iv,\,\, i=0,\dots n-1$ (resp the ${}^tX^jv^*,\,\,
j=0,\dots n-1)$ gives a non trivial linear relation among the rows (resp the columns)
of $A$ so that if the triplet is not regular then $D(X,v,v^*)=0$.

To prove the second assertion we simply have to exhibit a regular element. For
example choose a basis $e_1,\dots e_n$ of $V$ and let $e_1^*,\dots ,e_n^*$ be the dual
basis of $V^*$. Then define $X$ by
\[ X(e_i)=e_{i+1},\,\, i=1,\dots n-1,\quad X(e_n)=0\quad {\rm and}\quad
v=e_1,v^*=e^*_n \] 
Then 
\[ {}^tXe^*_i=e^*_{i-1},\,\, i=2,\dots ,n,\quad {}^tXe^*_1=0\]
it follows that $(X,v,v^*)$ is regular.
\begin{proposition}
 Two regular elements are conjugate under $G$ if and only if they
give the same values to the invariants $q_j$ and $D_j$.
\end{proposition}
The necessity is clear. Conversely let $(X,v,v^*)$ and $(Y,w,w^*)$ be two regular
elements such that
\[ q_j(X,v,v^*)=q_j(Y,w,w^*),\quad D_j(X,v,v^*)=D_j(Y,w,w^*)\quad j=1,\dots n-1\] 
In particular 
\[ A(X,v,v^*)=A(Y,w,w^*)\]
Let $g$ be the linear map from $V$ to $V$ defined by $g(X^pv)=Y^pw$ for $p=0,\dots
n-1$. We claim that $gXg^{-1}=Y$. It is enough to check that
\[ g^{-1}YgX^pv=X^{p+1}v,\quad p=0,\dots n-1\]
Now, if $p\leq n-2$
\[ g^{-1}YgX^pv=g^{-1}Y^{p+1}w=X^{p+1}v\] 
For $p=n-1$
\[ g^{-1}YgX^{n-1}v=g^{-1}Y^nw=g^{-1}\sum _0^{n-1}D_j(Y)Y^jw=\sum _{0}^{n-1}D_j(X)X^jv
=X^nv\] 
Also
\[ \langle  {}^tg^{-1}v^*,Y^jw\rangle  =\langle  v^*,g^{-1}Y^jw\rangle  =\langle  v^*,X^jv\rangle  =q_j(X,v,v^*)\] 
and
\[ \langle  w^*,Y^jw\rangle  =q_j(Y,w,w^*)\]
Hence, for all $p$ we get
\[ \langle  {}^tg^{-1}v^*,Y^jw\rangle  =\langle  w^*,Y^jw\rangle  \] 
which implies that ${}^tg^{-1}v^*=w^*$. The two triplets are conjugate by $g$.

We define the {\bf regular orbits} as the orbits of regular elements. Each
such orbit is defined by the values of the $q_j$ and the $D_j$. These values
must be such that $D$, which is a polynomial in the $q_j$ should not take the
value 0 and also the relations between the $q_j$ and
$D_j$ must be satisfied by the chosen values. Each such orbit is (Zariski)
closed (the invariants are constant  on the closure of any orbit !). If
$(X,v,v^*)$ is regular and if $g\in G$ is such that $g(X,v,v^*)=(X,v,v^*)$ which means in
particular  that $gv=v$ and $gXg^{-1}=X$, then $gX^pv=(gXg^{-1})^pgv=X^pv$ for all $p$. By
definition of regular $\{v,Xv,\dots ,X^{n-1}v\}$ is a basis of $V$ so we
conclude that $g={\rm Id}$: the isotropy subgroup of a regular element is
trivial. In fact
\begin{theorem}
  An orbit  is regular if and only if it is closed and if the centralizer in $G$ of
an element of the orbit is trivial..
\end{theorem}
Let $(X,v,v^*)$ be a non regular triplet. Suppose for example that the subspace $E$ of
$V$ generated by $v,Xv,\dots ,X^{n-1}v$ is a proper subspace. Choose a subspace $F$
such that $V=E\oplus   F$. With respect to this decomposition we can write $X$ as
\[ X=\pmatrix{A&B\cr 0&C
\cr}\]
For t$\in {\mathbb F}^*$ define
\[ a_t=\pmatrix{{\rm Id}_ E&0\cr 0&t^{-1}{\rm Id}_ F
\cr}\] 
Then $v$ is fixed by $a_t$ and
\[ a_tXa_t^{-1}=\pmatrix{A&tB\cr 0&C
\cr}\] 
so that
\[ \lim_{t\rightarrow 0}a_tXa_t^{-1}=\pmatrix{A&0\cr 0&C
\cr}\] 
Furthermore we have $V^*=E^*\oplus   F^*$; decompose $v^*=v_1^*+v_2^*$. Then
\[ {}^ta_t^{{-1}}v^*=v_1^*+tv_2^*\] 
so that
\[ \lim_{t\rightarrow 0}{}^ta_t^{{-1}}v^*=v_1^*\] 
Note that the limit point  (when $t $ goes to 0 ) of $(X,v,v^*)$ is fixed by $a_t$. Thus either the orbit is not closed or, if it is closed the centralizer of one of its point is not trivial.

We can also give a characterization of the closed orbits. as before let $E$ be the
subspace generated by  $v,Xv,X^2v,\dots$ and let $E'$ be the subspace generated by
$v^*,{}^tXv^*,{}^tX^2v^*,\dots$. Choose a subspace $F$ of $V$ such that $V=E\oplus   F$.
Then, relative to this decomposition, we may write \[ X=\pmatrix{A&B\cr0&C
\cr}\]
 
If we view $C$ as a linear map from $V/E$ into itself then it does not depend upon the
choice of $F$.
\begin{theorem}
 The orbit of $(X,v,v^*)$ is closed if and only if $V^*=E'\oplus   E^\perp $ and  $C$ is
semi-simple.
\end{theorem}
The condition  $V^*=E'\oplus   E^\perp $ is equivalent to the condition $V=E\oplus   E^{'\perp}$ and
means that $E'$ (resp. $E$) may be identified with the dual space $E^*$  (resp $E^{'*}$) of
$E$ (resp $E'$). Suppose that this condition is not satisfied and that, for example $\dim
E'\geqslant \dim E$ Define $a_t$ as before. Then
\[ \lim _{t\rightarrow 0}(X,v,v^*)=(Y, v, e^*)\] 
where $e^*$ is obtained as follows: we identified $V^*$ with $E^*\oplus   F^*$ and write
$v^*=e^*+f^*$ and
\[ Y=\pmatrix{A&0\cr 0&C
\cr}\]
The subspace $E'(Y,e^*)$ generated by $e^*,{}^tYe^*,\dots$ is contained in $E^*$ so that
\[ \dim E'(Y,e^*)\leq \dim E^*=\dim E\leq \dim E'\] 
If $(X,v,v^*)$ is conjugate to $(Y,v,e^*)$ we must have $\dim E'(Y,e^*)=\dim E'$ so that
$E'(Y,e^*)=E^*$. But the subspace $E(Y,v)$ generated by $v,Yv,\dots$ is $E$ so that
$(Y,v,e^*)$ satisifes the first condition of the theorem, hence can not be conjugate to
$(X,v,v^*)$.

Thus, looking for closed orbits, we may assume that $V^*=E'\oplus   E^\perp$ and
$V=E\oplus   E^{'\perp}$. So we choose $F=E^{'\perp}$ and identify $F^*$ with $E^\perp$ and
$E^*$ with $E'$. Now
\[ {}^tX=\pmatrix{
{}^tA&0\cr {}^tB&{}^tC\cr}\] 
However $E'=E^*$ is stable under ${}^tX$ so that $B=0$.

Suppose first that the orbit of $(X,v,v^*)$ is closed. Let $D$ be some element in the
closure of the orbit ${\rm GL} ({F})C$. There exists a sequence $(u_n)$ of elements of
${\rm GL}(F)$ such that $u_nCu_n^{-1}\rightarrow D$. Put
\[ g_n=\pmatrix{{\rm Id}_ E&0\cr 0&u_n
\cr}\]
Then
\[ g_nv=v,\quad ^tg_n^{-1}v^*=v^*,\quad g_nXg_n^{{-1}}\rightarrow \pmatrix{A&0\cr
0&D
\cr}\] 
As the orbit is closed there exists $g\in G$ such that
\[ gv=v,\quad ^tg^{-1}v^*=v^*,\quad gXg^{-1}=\pmatrix{A&0\cr
0&D
\cr}\] 
Now $X^qv=A^qv$ and $g(X^qv)=(gXg^{-1})^qg(v)=A^qv$ so that $g$ is the identity on $E$
and, by the same argument, ${}^tg^{-1}$ is the identity on $E^*$. Finally $g$ may
be written as
\[ g=\pmatrix{{\rm Id}_ E&0\cr 0&\delta 
\cr}\] 
and $D=\delta C\delta ^{-1}$ which shows that $D$ belongs to the orbit of $C$. Hence this orbit is
closed which is  equivalent to the semi-simplicity of $C$.

Conversely assume that $C$ is semi-simple. Let $(g_n)$ be a sequence of elements of
$G$ such that $g_n(X,v,v^*)$ has a limit $(Y,u,u^*)$. Consider the subspace $E(Y,u)$
generated by $u, Yu,Y^2u,\dots$. For any $q$ the sequence $(g_n(X^qv))$ converges to
$Y^qu$. Hence any linear relation between the vectors $X^qv$ remains valid for the
$Y^qu$. In particular if $E$ is of dimension $p$ then $E(Y,u)$ is of dimension at most $p$,
the set $\{v,Xv.\,\dots ,X^{p-1}v\}$ is a basis of $E$ and the set $\{u,Yu,\dots
,Y^{p-1}u\}$ is a set of generators of $E(Y,u)$. The matrix $\left 
(\langle  {}^tX^{i-1}v^*,X^{j-1}v)\right  ),\,\, 0\leq i,j\leq p-1$ is invariant under $G$ and, as
$E'\approx E^*$ is non singular. Hence the matrix  $\left 
(\langle  {}^tY^{i-1}u^*,Y^{j-1}u)\right  ),\,\, 0\leq i,j\leq p-1$ is also non singular. This implies
that the vectors $u,Yu,\dots ,Y^{p-1}u$ are linearly independant. Similarly we prove
that $E'({}^tY,u)$ the subspace of $V^*$ generated by $u^*,{}^tYu^*,\dots$ has
dimension $p$ and admits $\{u^*,{}^tYu^*,\dots ,{}^tY^{p-1}u^*\}$ as basis. Also the
pairing between $E(Y,u)$ and $E'({}^tY,u^*)$ is non singular. 

Choose $g\in G$ such that $g(F)=E'({^tY,u^*)^\perp} $ and $g(X^qv)=Y^qu,\,\,
0\leq q\leq p-1$. Note that ${}^tg^{-1}{}^tX^qv^*\in E'({}^tY,u^*)$.We have
\[ \langle  {}^tg^{-1}{}^tX^qv^*, gX^rv\rangle  =\langle  {}^tX^qv^*, X^rv\rangle  =\langle  {}^tY^qu^*,Y^ru\rangle   \]
and
\[ \langle  {}^tg^{-1}{}^tX^qv^*, gX^rv\rangle  =\langle  {}^tg^{-1}{}^tX^qv^*, Y^ru\rangle  \] 
This being true for $0\leq r\leq p-1$ implies that ${}^tg^{-1}{}^tX^qv^*={}^tY^qu^*$.

If we replace  the sequence $(g_n)$ by the sequence $g^{{-1}}g_n$ we see that we can
assume that, for all $e\in E$  and $e^*\in E^*$ one has
\[ \lim g_n(e)=e,\quad \lim {}^tg_n^{-1}e^*=e^*\]
In particular $X^rv=Y^ru$ for all $r$.
Next consider the decomposition $V=E\oplus   E^{*{\perp}}$ and put
\[ g_n=\pmatrix{\alpha  _n&\beta _n\cr \gamma _n&\delta _n
\cr}\quad g_n^{-1}=\pmatrix{a_n&b_n\cr c_n&d_n
\cr}\]
We get
\[ \lim \alpha  _n={\rm Id}_ E,\,\, \,\,
\lim \gamma _n=0\]
 Also
\[ X=\pmatrix{A&0\cr 0&C
\cr}\,\,\quad Y=\pmatrix{A&0\cr 0&D
\cr}\] 
From $g_ng_n^{-1}=g_n^{-1}g_n={\rm Id}$ we get
\begin{eqnarray*}
\alpha  _nb_n+\beta _nd_n&=&0\\
\gamma _nb_n+\delta _nd_n&=&1 
\end{eqnarray*}
For $n$ large enough $\alpha  _n$ is invertible so that we may write
\[ b_n=-\alpha  _n^{-1}\beta _nd_n\]
and
\[ -\gamma _n\alpha  _n^{-1}\beta _nd_n+\delta _nd_n=1\]
which shows that $d_n$is invertible with inverse $d_n^{-1}=-\gamma _n\alpha  _n^{-1}\beta _n+\delta _n$.

From $\lim g_nXg_n^{-1}=Y$ we get
\begin{eqnarray*}
\alpha  _nAb_n+\beta _nCd_n&\rightarrow& 0\\
\gamma _nAb_n+\delta _nCd_n&\rightarrow& D
\end{eqnarray*}
Put
\[\beta _nCd_n=-\alpha  _nAb_n+\varepsilon _n,\quad \lim \varepsilon _n=0\] 
Then
\begin{eqnarray*} D=\lim(\gamma _nAb_n+\delta _nCd_n)&=&\lim
(\gamma _nAb_n+d_n^{-1}Cd_n+\gamma _n\alpha  _n^{-1}\beta _nCd_n)\\
&=&\lim \left  (\gamma _nAb_n+d_n^{-1}Cd_n+\gamma _n\alpha  _n^{-1}(-\alpha  _nAb_n+\varepsilon _n)\right  )\\
&=&\lim (d_n^{-1}Cd_n+\gamma _n\alpha  _n^{-1}\varepsilon _n)\\
&=&\lim (d_n^{-1}Cd_n)
\end{eqnarray*} 
As $C$ is semi-simple its orbit under ${\rm GL}(E^{*\perp})$ is closed. Hence, for some
$u\in {\rm GL}(E^{*\perp})$ we have $D=uCu^{{-1}}$. Let
\[ g=\pmatrix{{\rm Id}_ E&0\cr 0&u
\cr}\] 
We get
\[ g(X,v,v^*)=(Y,u,u^*)\]
The orbit of $(X,v,v^*)$ is closed.

Let us go back to the regular orbits. The invariant polynomials $q_i$ and $D_j$ are also invariant under the involution $\sigma $. Then  Proposition 6-2 implies that the regular orbits are also stable by $\sigma $. For the p-adic topology each of them is homeomorphic to $G$ and carries an invariant measure which is symmetric with respect to $\sigma $. It follows from a classical result of Gelfand and Kazdhan  that a distribution on the regular set, invariant under $G$ is symmetric with respect to $\sigma $. Therefore an invariant distribution on ${\mathfrak {g} }\times V\times V^*$ which is skew symmetric with respect to $\sigma $ is supported in the complement of the regular set.

We want to prove  similar results for $G$ acting on $G\times V\times V^*$. 

For $v\in V$,$v^*\in V^*$ and $x\in G$ put
\[ q_i(x,v,v^*)=\langle  v^*,x^iv\rangle  =\langle  {}^tx^iv^*,v\rangle  \] 
and
\[ \mathrm {Det}  (T\,{\rm Id} -x)=T^n-\sum _0^{n-1}D_j(x)T^j\] 
The $q_i$ and $D_i$ are invariant polynomial functions. let
\[ A(x,v,v^*)=(a_{i,j}),\quad a_{i,j}=\langle  {}^tx^{j-1}v^*,x^{i-1}v\rangle  ,\,\, i,j=1,\dots n\] 
an $n\times n$ matrix and put
\[ D(x,v,v^*)=\mathrm {Det} (A(x,v,v^*)\] 

A triplet $(x,v,v^*)$ is {\bf regular} if $\{ v,xv,\dots x^{n-1}v\}$ is a basis of $V$ and
$\{ v^*,{}^txv^*,\dots ,{}^tx^{n-1}v^*\}$ is a basis of $V^*.$
\begin{proposition}
 The triplet $(x,v,v^*)$ is regular if and only if $D(x,v,v^*)\ne 0$. The set of
regular elements is a non empty Zariski open subset.
\end{proposition}
The proof is the same as the proof of Proposition 6-1 except that to exhibit a regular element we replace the principal nilpotent element $X$ by the unipotent element $x=1+X$.

\begin{proposition}
 Two regular elements are conjugate under $G$ if and only if they
give the same values to the invariants $q_j$ and $D_j$.
\end{proposition}
Exactly the same proof as Proposition 6-2. In fact if we use the  inclusion $G\subset{\mathfrak {g} }$ it is even a particular case.

We define the {\bf regular orbits} as the orbits of regular elements. Each
such orbit is defined by the values of the $q_j$ and the $D_j$. These values
must be such that $D$, which is a polynomial in the $q_j$ should not take the
value 0 and also the relations between the $q_j$ and
$D_j$ must be satisfied by the chosen values; also as we want $x$ to be invertible, we require that $D_0\ne 0$. Each such orbit is (Zariski)
closed (the invariants are constant  on the closure of any orbit !). If
$(x,v,v^*)$ is regular and if $g\in G$ is such that $g(x,v,v^*)=(X,v,v^*)$ which means in
particular  that $gv=v$ and $gxg^{-1}=x$, then $gx^pv=(gxg^{-1})^pgv=x^pv$ for all $p$. By
definition of regular $\{v,Xv,\dots ,X^{n-1}v\}$ is a basis of $V$ so we
conclude that $g={\rm Id}$: the isotropy subgroup of a regular element is
trivial. In fact
\begin{theorem}
  An orbit  is regular if and only if it is closed and if the centralizer in $G$ of
an element of the orbit is trivial..
\end{theorem}
Same proof as for Theorem 6-1. If we use the fact that $G\subset{\mathfrak {g} }$ it is a particular case because the determinant being constant on the closure of an orbit, the closure in ${\mathfrak {g} }$ is the same as the closure in $G$.

Finally we can also give a characterization of the closed orbits.  Let $E$ be the
subspace generated by  $v,xv,x^2v,\dots$ and let $E'$ be the subspace generated by
$v^*,{}^txv^*,{}^tx^2v^*,\dots$. Choose a subspace $F$ of $V$ such that $V=E\oplus   F$.
Then, relative to this decomposition, we may write \[ x=\pmatrix{A&B\cr0&C
\cr}\] 
If we view $C$ as a linear map from $V/E$ into itself then it does not depend upon the
choice of $F$.
\begin{theorem}
 The orbit of $(x,v,v^*)$ is closed if and only if $V^*=E'\oplus   E^\perp $ and  $C$ is
semi-simple.
\end{theorem}
We already remarked that the closure in ${\mathfrak {g} }$ is equal to the closure in $G$ so it is a particular case of Theorem 6-2.

\begin{section}
 {\hglue -12pt . A Frobenius type descent.}
\end{section}
With the same notations as in the preceeding section we want to show that an invariant distribution on ${\mathfrak {g} }\times V\times V^*$, skew with respect to $\sigma $ must be supported by the set of all $(X,v,v^*)$ such that $\langle  v^*,X^iv\rangle  =0$ for all $i$.

 For $1\leq r\leq n$ define
\[ E_r(X,v)=\sum _{j=0}^{r-1}{\mathbb F}X^jv,\quad E'_r(X,v^*)=\sum _{j=0}^{r-1}{\mathbb F}\,\,^tX^jv^* \] 
For $1\leq r\leq n$ let $\Sigma(r)$ be the set of all $(x,v,v^*)$ such that $E_r(X,v)$ and $E_r'(X,v^*)$ are both of
dimension $r$ and that $\langle  .,.\rangle  $ defines a non degenerate duality between these two spaces.As before put
\[ q_j=\langle  v^*,X^jv\rangle  ,\quad a_{i,j}=\langle  ^tX^{j-1}v^*,X^{{i-1}}v\rangle  =q_{i+j-2}\]
and let $A_r(X,v,v^*)$ be the $r\times r$ matrix with coefficients $a_{i,j}$ for $1\leq i,j\leq r$. Then $(X,v,v^*)\in \Sigma (r)$
if and only if the matrix $A_r$ is non singular. It follows that $\Sigma (r)$ is open; also $\Sigma (n)$ is the
subset of regular elements. 
\begin{lemma}
  Each $\Sigma (r)$ is an open subset. The complement of the union of all $\Sigma (r)$ is the
set of all $(X,v,v^*)$ such that  $E'_n(X,v^*)$ is contained in $E_n(X,v)^\perp$.
\end{lemma}
Indeed  $(X,v,v^*)$ does not belong to $\Sigma (1)$ if and only if $q_0=0$, does not belong to $\Sigma
(1)\cup \Sigma (2)$ if and only if $q_0=q_1=0$ and so on. It does not belong to $\Sigma (1)\cup \cdots \cup \Sigma (n)$ if
and only if $q_0=q_1=\cdots =q_{n-1}=0$. This last condition implies that all $q_r=0$ because $X^n$ is a linear
combination of the $X^j$ for $j=0,\dots ,n-1$
\hfill\break
 
Fix $r\leq n$. Let $E$ be a subspace of $V$ of dimension $r$ and choose a subspace $F$ such that $V=F\oplus   E$. identify
$V^*$ with $F^*\oplus   E^*$ so that, in particular $E^\perp=F^*$ and $F=E^{*\perp}$. Let $(X,v,v^*)\in \Sigma (r)$. Then
$V=E(X,v)\oplus   E'(X,v^*)^\perp$. Hence there exists $g\in G$ such that $gE(X,v)=E$ and $gE'(X,v^*)^{\perp}=F$ which
implies that $^tg^{-1}E'(X,v^*)=E^*$. Let $\Xi (r)$ be the set of all $(X,v,v^*)\in \Sigma (r)$ such that $E(X,v)=E$ and
$E'(X,v^*)=E^*$; it is a closed subset and we just saw that $G\Xi (r)=\Sigma(r)$. Furthermore let
$(X,v,v^*)\in \Xi (r)$ and suppose that for some $g\in G$, we have $g(X,v,v^*)\in \Xi (r)$. put $g(X,v,v^*)=(X' ,v',v^{'*})$ We
know that $\{ v,Xv,\dots ,X^{r-1}v\}$ and $\{ v',X'v',\dots ,X'{}^{r-1}v'\}$ are two basis of $E$ and that $g(X^jv)=X'{}^jv'$.
It follows that $g(E)=E$. The same argument shows that $^tg(E^*)=E^*$. If we use the decompositions $V=F\oplus   E$
and $V^*=F^*\oplus   E^*$ to write $g$ and $^tg$ as two by two matrices  then
\[ g=\pmatrix{g_ F&0\cr 0&g_ E
\cr}\] 
with $g_ F\in {\rm GL}(F)$ and $g_ E\in {\rm GL}(E)$ so that $g\Xi (r)=\Xi (r)$. and we are in a situation to use the results of
the Appendix. We call $H$ the stabilizer of $\Xi (r)$  in $G$. In matrix form, as above, it is the subgroup of
diagonal matrices.

We have to look in more details at $\Xi (r)$. Let $(X,v,v^*)\in \Xi (r)$ and put
\[ X=\pmatrix{x_{1,1}&x_{1,2}\cr x_{2,1}&x_{2,2}
\cr}\] 
We have $x_{1,2}X^jv=0$ for $j=0,1,\dots ,r-2$. Let $w^*\in E^*$ be defined by
\[ \langle  w^*,\sum _0^{r-1}\lambda _jX^jv\rangle  =\lambda _{r-1}\] 
and put 
\[ u=x_{1,2}X^{r-1}v\] 
Then 
\[ x_{1,2}e=\langle  w^*,e\rangle  u,\quad e\in E\] 
In a similar way. we define $w\in E$ by
\[ \langle  \sum _0^{r-1}\mu _j\,\,^tX^jv^*,w\rangle  =\mu _{r-1}\] 
and
\[ u^*={}^tx_{2,1}\,\, ^tX^{r-1}v^*\] 
to get
\[ {}^tx_{2,1}e^*=\langle  e^*,w\rangle  u^*,\quad e^*\in E^*\] 
Note that $(x_{1,1},u,u^*)\in {\mathfrak {g} }_ F\times F\times F^*$ and $(x_{2,2},v,v^*)\in {\mathfrak {g} }_ E\times E\times E^*$. and is a
regular element. Conversely if we start with two such elements then we can recover $(X,v,v^*)\in \Xi (r)$ by the
above formulas. Thus if we denote by $({\mathfrak {g} }_ E\times E\times E^*)_{\rm reg}$ the open set of regular elements
then we have an homeomorphism
\[ \Xi (r)\rightarrow \biggl ({\mathfrak {g} }_ F\times F\times F^*\biggr )\times \biggl (({\mathfrak {g} }_ E\times E\times E^*)_{\rm reg}\biggr )\]
Next let us look at the action of $H$; take
\[ h=\pmatrix{h_ F&0\cr 0&h_ E
\cr}\] 
Then $h$ acts by\overfullrule=0pt
\[ \displaylines{u\mapsto h_ F u,\,\,\, u^*\mapsto {}^th_ F^{-1}u^*,\cr
 v\mapsto h_ Ev,\,\,\,
v^*\mapsto {}^th_ E^{-1}v^*,\cr
 x_{1,1}\mapsto h_ Ex_{2,2}h_ E^{-1},\,\,\, x_{1,1}\mapsto
h_ Fx_{1,1}h_ F^{-1}} \]
The action of $H$ on $\Xi (r)$ is the product of the action on the two components.
\begin{proposition}
 $(X,v,v^*)$ is regular if and only if $(x_{1,1},u,u^*)$ is regular.
\end{proposition}
Indeed $v,Xv,\dots ,X^{r-1}v$ is a basis of $E$ and ,by a trivial induction, there exists constants $\lambda _{i,j}$ such that
\[ X^{r+s}v\equiv x_{1,1}^su+\lambda _{s,1}x_{1,1}^{s-1}u+\dots +\lambda _{s,s}u\pmod E \] 
so that $v,Xv,\dots ,X^{n-1}v$ is a basis of $V$ if and only if $u,x_{1,1}u,\dots ,x_{1,1}^{n-r+1}u$ is a basis of $F$. The same argument is valid for $v^*$ and $u^*$, hence the result.

To define the involution $\sigma $ relative to $V$ we start from a basis of $E$ and a basis of $F$. Then we have two maps $u_ E$ and $u_ F$ from $E$ and $F$ to their respective dual space, sending the basis to the dual basis. we define $u=(u_ F, u_ E)$. Recall that
\[ \sigma (X,v,v^*)=(u^{-1}{}^tXu,\, u^{-1}(v^*),\, u(v))\] 
Then $\Sigma (r)$ is stable. We use Proposition 1
 of the Appendix: for any distribution $T$ on $\Sigma (r)$, invariant under $G$, there exists a unique distribution $S$ on $\Xi (r)$, invariant under $H$ and such that, for $f\in {\mathcal S}(\Sigma (r))$
 \[ \langle  T,f\rangle  =\int_{G/H}\langle  S,f(gx)d\mu (g)\]
 In this case $\mu $ is an invariant measure on $G/H$. Now $\Xi (r)$ is homeomorphic to $\biggl ({\mathfrak {g} }_ F\times F\times F^*\biggr )\times \biggl (({\mathfrak {g} }_ E\times E\times E^*)_{\rm reg}\biggr )$. The involution $\sigma $ fixes $\Xi (r)$ and its restriction is the product of the partial involutions  $\sigma _ F$ and $\sigma _ E$.

\begin{proposition}
 Let $(X,v,v^*)\in \Xi (r)$. The orbit $G(X,v,v^*)$ is fixed by $\sigma $ if and only if the orbit
${\rm GL}(F)(x_{1,1},u,u^*)$ is fixed by $\sigma _ F$. If  any $G_ F-$invariant distribution on ${\mathfrak {g} }_ F\times F\times F^*$ is invariant under $\sigma _ F$ then
 then any $G-$invariant distribution on $\Sigma (r)$ is invariant by $\sigma $.
\end{proposition}
The first assertion follows from the above remarks and the fact that any regular orbit is fixed by the
involution.. For the second we have to check that the map $S\mapsto T$ given by  the Appendix is
compatible with $\sigma $. Recall that $\sigma \circ g=g'\circ \sigma $ with $g'=u^{-1}\circ {}^tg^{-1}\circ u$ where $g\in G$ and,
as before
\[ u=\pmatrix{u_ F&0\cr 0&u_ E
\cr}\] 
 As $\sigma $ is an involution so is the map
$g\mapsto g'$, hence the Haar measure of $G$ is invariant under this map and the same is true for the stable
subgroup $H$ and so also for the quotient measure on $G/H$. If we replace $S$ by $\sigma (S)$ we obtain
\[ \int _{G/H}\langle  S,f(g\sigma (x)\rangle  d\mu (g)=\int _{G/H}\langle  S,f(\sigma g'x)\rangle  d\mu (g)=\int _{G/H}\langle  S,f(\sigma gx)\rangle  d\mu (g)\] 
so that $T$ is replaced by $\sigma (T)$.

As we proceed by induction on $n$ we conclude that any invariant distribution on ${\mathfrak {g} }\times V\times V^*$, skew with respect to $\sigma $, is 0 on $\Sigma(r)$. This is true for any $r$ ( the case  $r=n$ is the case of regular orbits). 

Let $\Sigma \subset [{\mathfrak {g} },{\mathfrak {g} }]\times V\times V^*$ be the set of all $(X,v,v^*)$ such that $X\in {\mathcal N}$, the nilpotent cone and that, for all $i$, $\langle  v^*,X^iv\rangle  =0$. Recall that ${\mathfrak z }$ is the center of ${\mathfrak {g} }$.
\begin{lemma}
 If any invariant distribution on $[{\mathfrak {g} },{\mathfrak {g} }]\times V\times V^*$
 is symmetric with respect to $\sigma $ then any invariant distribution on ${\mathfrak {g} }\times V\times V^*$ is symmetric with respect to $\sigma $.
 \end{lemma} 
This is clear as $G$ and $\sigma $ act trivially on ${\mathfrak z }$.

{\noindent\bf Conclusion\pointir The support of an invariant distribution on ${\mathfrak {g} }\times V\times V^*$, skew with respect to $\sigma $ is contained in $\Sigma$}

\section {\hglue -12pt 
 Invariant distributions with singular support}

Our first goal is to find extra conditions satisfied by these distributions.

\subsection
 {Some recollections}

Let $E$ be a vector space of dimension $n$ over ${\mathbb  F}$. Let $Q$ be a non
degenerate quadratic form on $E$ and let
\[ B(X,Y)=Q(X+Y)-Q(X)-Q(Y)\] 
We define the Fourier transform by
\[ \hat f(Y)=\int _Ef(X)\tau (B(X,Y))dX\] 
where $\tau $ is a non trivial additive character of ${\mathbb  F}$ and $dX$ the Haar measure
on $E$ such that the inversion formula is
\[ f(X)=\int _E\hat f(Y)\tau (-B(X,Y))dY\] 
We denote by $(a|b)$ he Hilbert symbol of ${\mathbb  F}$. The metaplectic group ${\mathrm 
Mp}_2({\mathbb  F})$ is the twofold covering of ${\mathrm  SL}_2({\mathbb  F})$ defined as
follows. If
\[ g=\pmatrix{a&b\cr c&d
\cr}\in {\mathrm  SL}_2({\mathbb  F})\] 
then we let $j(g)=c$ if $c\ne 0$ and $j(g)=d$ if $c=0$. For $g,g'\in {\mathrm  SL}_2({\mathbb  F})$ and
$g''=gg'$ put
\[ \varepsilon (g,g')=(j(g)j(g'')|j(g')j(g'')\] 
Then an element of the metaplectic group is written as $(g,\varepsilon )$ with $g\in {\mathrm 
SL}_2({\mathbb  F})$ and $\varepsilon =\pm 1$ and the group law is
\[ (g,\varepsilon )(g',\varepsilon ')=(gg',\varepsilon \varepsilon '\varepsilon (g,g'))\] 
To the quadratic form $Q$ ( and the fixed character $\tau $) is attached an eight root
of unity $\gamma (Q)$ and there exists a (unique) representation $\pi _Q$ of ${\mathrm  Mp}_2({\mathbb 
F})$ into ${\mathcal  S}(E)$ such that
\begin{eqnarray*}
 \pi _Q\Biggl [\pmatrix{1&u\cr 0&1
\cr},\varepsilon \Biggr ]f(X)&=&\varepsilon ^n\tau (uQ(X))f(X)\\
\pi _Q\Biggl [\pmatrix{t&0\cr 0&t^{-1}
\cr},\varepsilon \Biggr ]f(X)&=&\varepsilon ^n\frac{\gamma (Q)}{\gamma (tQ)}|t|^{n/2}f(tX)\\
\pi _Q\Biggl [\pmatrix{\hfill 0&1\cr -1&0
\cr},\varepsilon \Biggr ]f(X)&=&\varepsilon ^n\gamma (Q)\widehat f(X)\\
\end{eqnarray*}

Let $\Gamma _0$ be the cone of all $X\in E$ such that $Q(X)=0$.
\begin{lemma}
 Let $T\in {\mathcal  S}'(E)$ be such that both $T$ and $\widehat T$ have support
contained in $\Gamma _0$. Then $T$ is homogeneous:
\[ \langle T,f(tX)\rangle =|t|^{-n/2}\frac{\gamma (tQ)}{\gamma (Q)}\langle T,f\rangle \] 
In particular if $n$ is odd then $T=0$.
\end{lemma}
The condition on the support of $T$ may be written as
\[ \langle T,\pi _Q\Biggl [\pmatrix{1&u\cr 0&1
\cr},\varepsilon \Biggr ]f\rangle =\varepsilon ^n\langle T,f\rangle \] 
The condition on the support of $\widehat T$ may be written as
\[ \langle \widehat T,\pi _Q\Biggl [\pmatrix{1&u\cr 0&1
\cr},\varepsilon \Biggr ]f\rangle =\varepsilon ^n\langle \widehat T,f\rangle \] 
or
\overfullrule=0pt
\[ \langle T,\gamma (Q)^{-1}\pi _Q\Biggl (\Biggl [\pmatrix{\hfill 0&1\cr -1&0
\cr},1\Biggr ]\Biggl [\pmatrix{1&u\cr 0&1
\cr},\varepsilon \Biggr ]\Biggr )f\rangle =\varepsilon ^n\langle T,\gamma (Q)^{-1}\Biggl [\pmatrix{\hfill 0&1\cr -1&0
\cr},1\Biggr ]f\rangle \] 
or
\[ \langle T,\pi _Q\Biggl (\Biggl [\pmatrix{\hfill 0&1\cr -1&0
\cr},1\Biggr ]\Biggl [\pmatrix{1&u\cr 0&1
\cr},\varepsilon \Biggr ]\Biggl [\pmatrix{\hfill 0&1\cr -1&0
\cr},1\Biggr ]^{-1}\Biggr )f\rangle =\varepsilon ^n\langle T,f\rangle \] 
Computing the product we get
\[ \langle T,\pi _Q\Biggl [\pmatrix{1&0\cr v&1
\cr},\varepsilon \Biggr ]f\rangle =\varepsilon ^n\langle T,f\rangle \] 
with $v=-u$.

But
\begin{eqnarray*}& &\Biggl [\pmatrix{t&0\cr 0&t^{-1}
\cr},\varepsilon \Biggr ]=\\ \\
\Biggl [\pmatrix{1&-t\cr 0&1
\cr},1\Biggr ]\Biggl [\pmatrix{1&0\cr t^{-1}&1
\cr},\varepsilon \Biggr ]&&\Biggl [\pmatrix{1&-t\cr 0&1
\cr},1\Biggr ]\Biggl [\pmatrix{1&1\cr 0&1
\cr},1\Biggr ]\Biggl [\pmatrix{\hfill 1&0\cr -1&1
\cr},1\Biggr ]\Biggl [\pmatrix{1&1\cr 0&1
\cr},1\Biggr ]\\
\end{eqnarray*}
Hence
\[ \langle T,\pi _Q\Biggl [\pmatrix{t&0\cr 0&t^{-1}
\cr},\varepsilon 
\Biggr ]f\rangle =\varepsilon ^n\langle T,f\rangle \] 
and we obtain the homogeneity of $T$. 

If $T\ne 0$ then $\gamma (tQ)/\gamma (Q)$ must be a multiplicative character which is true if and
only if $n$ is even.

Note also that
\[ \Biggl [\pmatrix{\hfill 0&1\cr -1&0
\cr},\varepsilon \Biggr ]=\Biggl [\pmatrix{1&-1\cr 0&1
\cr},\varepsilon \Biggr ]\Biggl [\pmatrix{1&0\cr 1&0
\cr},1\Biggr ]\Biggl [\pmatrix{1&-1\cr 0&\hfill 1
\cr},1\Biggr ]\Biggl [\pmatrix{-1&0\cr \hfill 0&1
\cr},1\Biggr ]\] 
which implies
\[ \langle T,\pi _Q\Biggl [\pmatrix{\hfill 0&1\cr -1&0
\cr},\varepsilon \Biggr ]f\rangle =\varepsilon ^n\langle T,f\rangle \] 
or, explicitly
\[ \gamma (Q)\widehat T=T\] 
Finally, if $n$ is even, let $D$ be the discriminant of $Q$; then
\[ \frac{\gamma (tQ)}{\gamma (Q)}=\bigl ((-1)^{n/2}D|\, t\bigr )\] 
\subsection
{Partial Fourier transform on $[{\mathfrak  g},{\mathfrak  g}]$}
Now we go back to the ${\mathrm  GL}-$case. Let  
 ${\mathfrak  g}_1=[{\mathfrak  g},{\mathfrak  g}]$ be the derived algebra.The center  will play no role hence, in this section, we systematically work with ${\mathfrak  g}_1$. Let $B$ be the
Killing form of ${\mathfrak  g}_1$ and consider the quadratic form
$R(X_1)=B(X_1,X_1)/2$. We choose, once for all, a non trivial  additive character $\tau $ of ${\mathbb  F}$.
 In this section
Fourier transform means partial Fourier transform with respect to ${\mathfrak  g}_1$
and is denoted by $f\mapsto {\mathcal  F}f$:
\[ {\mathcal  F}f(Y,v,v^*)=\int _{{\mathfrak  g}_1}f(X,v,v^*)\tau (B(X,Y))dX\] 
As $B$ is invariant under $G$, this transformation commutes with the action of $G$.
It also commutes with the involution $\sigma $. Indeed if $(e_i)$ is a basis of $V$ and
$(e_i^*)$ the dual basis of $V^*$ and if $u: V\rightarrow V^*$ is given by
$e_i^*=u(e_i)$ then we may normalize $\sigma $ by
\[ \sigma (X_1,v,v^*)=(u^{-1}\,^t\hskip -2pt X_1u,u^{-1}(v^*),u(v))\] 
Hence
\[ {\mathcal  F}(f\circ \sigma )(Y,v,v^*)=\int _{{\mathfrak  g}_1}f(u^ {-1}\,^t\hskip -2pt
Xu,u^{-1}(v^*),u(v))\tau (B(X,Y))\, dX\] 
The measure $dX$ is invariant by the involution $X\mapsto u^{-1}\, ^t\hskip -2pt
Xu$. For some $c_n>0$ we have $B(X,Y)=c_n{\mathrm Tr } (XY)$ so
\[ B(u^{-1}{}^t\hskip -2pt Xu,Y)=c_n{\mathrm Tr }(u^{-1}{}^t\hskip -2pt XuY)=c_n{\mathrm Tr }
(^t\hskip -2pt XuYu^{-1})=c_n{\mathrm Tr } (u^{-1}{}^t\hskip -1pt YuX)\] 
which means that
\[ B(u^{-1}{}^t\hskip -2pt Xu,Y)=B(X,u^{-1}{}^t\hskip -1pt Yu)\] 
and implies that
\[ {\mathcal  F}(f\circ \sigma )(Y,v,v^*)={\mathcal  F}f(u^{{-1}{}}\,^t\hskip -1pt
Yu,u^{-1}(v^*),u(v))\] 
It follows that the Fourier transform of a distribution which is invariant under
$G$ and skew relative to $\sigma $ is again invariant and skew. Suppose that $T$ is such a
distribution. We know that the support of $T$ is contained in the singular set  and the
same is true for ${\mathcal  F}T$. Let $\alpha \in {\mathcal  S}(V\oplus V^*)$; for $\varphi \in {\mathcal  S}({\mathfrak  g}_1)$ put $T_\alpha (\varphi )=T(\alpha \otimes \varphi )$. Then ${\mathcal  F}(T_\alpha )=({\mathcal  F}T)_\alpha $ so that both $T_\alpha $ and ${\mathcal  F}(T_\alpha )$ have their support contaimed in ${\mathcal  N}$. As the Killing form is 0 on ${\mathcal  N}$ we can apply Lemma 8-1. It tells us that $T$ is 0 if the dimension of ${\mathfrak 
g}_1$ is odd that is to say if $n^2-1$ is odd which means  $n$ is even.  {\textbf  From
now on we may assume that $n$ is odd}. We also get, on $[{\mathfrak  g},{\mathfrak  g}]$
\begin{eqnarray}
 \langle T,f(tX_1,v,v^*)\rangle =|t|^{-(n^2-1)/2}\langle T,f(X_1,v,v^*)\rangle 
 \end{eqnarray} 
and $T={\mathcal  F}T$ (in this case $\Delta =1$).
\begin{subsection}
{Fourier transform on $V\oplus V^*$}
\end{subsection}
Let $dv$ be a Haar measure on $V$ and let $dv^*$ be the dual Haar measure on $V^*$ so that if
\[ \psi (v^*)=\int _Vf(v)\tau (-\langle v^*,v\rangle )dv\] 
then
\[ f(v)=\int _{V^*}\psi (v^*)\tau (\langle v^*,v\rangle )dv^*\] 
On $V\oplus V^*$ we define the Fourier transform by
\[ \widehat f(w,w^*)=\int _{V\oplus V^*}f(v,v^*)\tau (-\langle v^*,w\rangle -\langle w^*,v\rangle )dv\, dv^*\] 
 and the inversion formula is
 \[ f(w,w^*)=\int _{V\oplus V^*}\widehat f(v,v^*)\tau (+\langle v^*,w\rangle +\langle w^*,v\rangle )dv\, dv^*\] 
 For $f\in {{\mathcal  S}}(V\oplus V^*)$ let
 \[ ({\mathcal  F}f)(v_1,v_2)=\int _{V^*}f(v_1,v^*)\tau (-\langle v^*,v_2\rangle )dv^*\] 
 be the partial Fourier transform.
If $g\in G$, put $f_g(v,v^*)=f(gv,\,^t\hskip -0.8pt g^{-1}v^*)$. Then
\[ {\mathcal  F}(f_g)(v_1,v_2)=|{{\mathrm {Det}}   }  g|({\mathcal  F}f)(gv_1,gv_2)\] 
Moreover
\[ ({\mathcal  F}\widehat f)(v_1,v_2)=({\mathcal  F}f)(-v_2,v_1)\] 
Let $P\in {\mathbb  F}[X]$ and 
\[ \varphi (v,v^*)=\varphi (v,v^*)\tau (-\langle v^*,P(X)v\rangle )\] 
Then
\[ ({\mathcal  F}\varphi )(v_1,v_2)=({\mathcal  F}f)(v_1,v_2+P(X)v_1)\] 
Let T
 be a distribution on ${\mathfrak  g}\times V\times V^*$, invariant under $G$, with support contained in the singular set and let $S={\mathcal  F}T$ be its partial Fourier transform on $V^*$. The support of $S$
 is contained in ${\mathcal  N}\times V\times V$. For any polynomial $P$ the function $\tau (\langle v^*,P(X)v\rangle )-1$ is equal to zero on the support of $T$, hence $(\tau (\langle v^*,P(X)v\rangle )-1)T=0$. This implies that
 \[ \langle S,\varphi (X,v_1,v_2+P(X)v_1)\rangle =\langle S,\varphi \rangle \] 
 Suppose that the distribution T is skew with respect to the involutionn$\sigma $ and consider $\widehat T$ the partial Fourier transform of $T$ with respect to $V\oplus V^*$. Then $\widehat T$ is invariant under $G$ and skew with respect to $\sigma $, hence its support is contained in the singular set. In particular
 \[ \tau (t\langle v^*,v\rangle )T=T,\quad \tau (t\langle v^*,v\rangle )\widehat T=T,\quad t\in {\mathbb  F}\] 
 On $V\oplus V^*$ the quadratic form $Q(v,v^*)=\langle v^*,v\rangle $ is non degenerate and $\gamma (Q)=1$. It follows that $T=\widehat T$. This in turn implies that
 \[ \langle S,\varphi (X,-v_2,v_1)\rangle =\langle S,\varphi \rangle \] 
 Recall that if $P$ is a polynomial and   $X$ a nilpotent endomorphism of $V$ then $P(X)$ is invertible if and only if $P(0)\ne 0$ in which case the inverse is $R(X)$ for some polynomial $R$. Consider the group $SL(2,{\mathbb  F}[X])$ of 2 by 2 matrices with coefficients polynomials in  such a $X$ and with determinant the identity endomorphism of $V$. The usual identities in $SL(2)$ imply that this group is generated by the matrices
 \[ \pmatrix{0&-1\cr 1&0\cr},\quad \pmatrix{1&0\cr P(X)&1\cr}\] 
As a distribution on ${\mathcal  N}\times V\times V$ the distribution $S$ is invariant under the transformations
\[ (X,v_1,v_2)\mapsto (X,-v_2,v_1),\quad (X,v_1,v_2)\mapsto (X,v_1,v_2+P(X)v_1)\] 
 Therefore it is invariant under the transformations
 \[ (X,v_1+v_2)\mapsto (X,A(X)v_1+B(X)v_2,C(X)v_1+D(X)v_2)\] 
 where $A,B,C,D$ are polynomials such that $A(X)D(X)-B(X)C(X)=1$. In particular, if $P$ is a polynomial such that $P(0)\ne 0$, it is invariant under
 \[ (X,v_1,v_2)\mapsto (X,P(X)v_1,P(X)^{-1}v_2)\] 
This means that the distribution $T$ satisfies the homogeneity condition
\begin{eqnarray}
  \langle T,f(X,P(X)v,{}^tP(X)v^*)=|{{\mathrm {Det}}   }  P(X)|^{-1}\langle T,f\rangle
  \end{eqnarray}
 We have to prove that a distribution $T$, having all the above properties, and $\sigma -$skew, is 0. 
 \subsection{Analysis on a fixed nilpotent orbit}
The nilpotent cone has a stratification given by the nilpotent orbits; each nilpotent orbit is stable by transposition. It would be enough to prove that, given a nilpotent orbit $\Gamma $, an invariant distribution $S$ on $\Gamma \times V\times V^*$, with singular support and $\sigma -$skew is 0. We may assume that this distribution satisfies the homogeneity conditions (1) and (2)    and is equal to its partial Fourier transform relative to $V\oplus V^*$.

If we fix $X\in \Gamma $then ,using the Appendix, we can transfer the problem to $V\oplus V^*$. Let $C$ be the centralizer of $X$ in $G$; it is known to be unimodular. Given a $T$ as above there exists a distribution $S$ on $V\oplus V^*$ such that, for $f\in {\mathcal  S}(\Gamma \times V\times V^*)$,
\[ \langle T,f\rangle =\int _{G/C}\langle S,f(gXg^{-1},gv,^{t}g^{-1}v^*)\rangle dg\] 
The distribution $S$ is invariant under $C$. We have to transfer to $S$ the properties of $T$. 

Let $\varphi \in {\mathcal  S}(\gamma )$; the function $g\mapsto \varphi (gXg^{-1})$ as a function on $G/C$ belongs to ${\mathcal  S}(G/C)$ and this is a bijection of ${\mathcal  S(\gamma )}$ onto ${\mathcal  S}(G/C)$. Let $f\in {\mathcal  S(V\oplus V^*)}$. The homogeneity condition $(2)$ implies that for any polynomial $P$ with non zero constant term
 \begin{eqnarray*}
 \int _{G/C}\varphi (gXg^{-1})&{}&\hskip -0.7cm \langle S,f(P(X)gv,^tP(X)^tg^{-1}v^*)\rangle dg\\
&{}&=|{{\mathrm {Det}}   }  P(X)|^{-1}\int _{G/C}\varphi (gXg^{-1})\langle S,f(gv,^tg^{-1}v^*)\rangle dg\\
\end{eqnarray*}
This is valid for any $\varphi $ therefore
\[ \langle S,f(P(X)gv,^tP(X)^tg^{-1}v^*)\rangle =|{{\mathrm {Det}}   }  P(X)|^{-1}\langle S,f(gv,^tg^{-1}v^*)\rangle \] 
In a similar way one shows that $S$ is equal to its Fourier transform $\widehat S$ on $V\oplus V^*$.

To transfer the homogeneity condition $(1)$ we first remark that there exists a one parameter group $\delta (t)$ in $G$ such that $\delta (t)X\delta (t)^{-1}=tX$. Indeed by Jacobson Morozov Theorem one can find in ${\mathfrak  g}$ two elements $H$ and $Y$ such that 
\[ [H,X]=2X,\,\, [H,Y]=-2Y,\,\, [X,Y]=H\] 
 The one parameter group with infinitesimal generator  $H/2$ will do. We shall be more explicit later on. Consider the adjoint representation of this TDS in ${\mathfrak  g}$. The subspace generated by the vectors having dominant weight is the Lie algebra ${\mathfrak  c}$ of the centralizer $C$ of $X$. Also note that if $\ell$ is an irreducible component with dominant vector $e$ then $[X,e]={{\mathrm {dim}}  }  \ell -1 $. It follows that the trace of the restriction of ${{\mathrm {ad}}  }  H$ to ${\mathfrak  c}$ is equal to the codimension of ${\mathfrak  c} $ in ${\mathfrak  g}$ that is to say to the dimension of $G/C$.
 
As ${{\mathrm {ad}}  }  \delta (t)X=tX $ this implies that the determinant of the restriction of ${{\mathrm {ad}}  }  \delta (t)$ to ${\mathfrak  c}$ is 
\[ {{\mathrm {Det}}   } \bigl ({{\mathrm {ad}}  }  \delta (t)_ {|\mathfrak  c}\bigr )=t^{\frac {1}{2}{{\mathrm {dim}}  }  G/C}\] 
The group $C$ is unimodular and normalized by ${{\mathrm {ad}}  }  \delta (t)$; let $dc$ be a Haar measure. We have
\[ d\bigl (\delta (t)c\delta (t)^{-1}\bigr )=|t|^{\frac {1}{2}{{\mathrm {dim}}  }  G/C}dc\] 
The Haar measure of $G$ is invariant under ${{\mathrm {ad}}  }  \delta (t)$ therefore if $\mu $ is an invariant measure on $G/C$ we get
\[ d\mu \bigl (g\delta (t)^{-1}\bigr )=|t|^{-\frac {1}{2}{{\mathrm {dim}}  }  G/C}d\mu (g)\] 
For $f$ and $\varphi $ as above 
 \begin{eqnarray*}
  \int _{G/C}\varphi (tgXg^{-1})&&\langle S,f(gv,^tg^{-1}v^*)\rangle d\mu (g)\\
&&=\int _{G/C}\varphi (g\delta (t)X\delta (t)^{-1}g^{-1})\langle S,f(gv,^tg^{-1}v^*)\rangle d\mu (g)\\
&&=|t|^{-\frac {1}{2}{{\mathrm {dim}}  }  G/C}\int _{G/C}\varphi (gX\delta g^{-1})\langle S,f(g\delta (t)^{-1}v,^tg^{-1}\,{}^t\delta (t)v^*)\rangle d\mu (g)\\
\end{eqnarray*} 
From the homogeneity condition $(1)$ we thus get
\[ \langle S,f(\delta (t)v,^t\delta (t)^{-1}v^*)\rangle =|t|^{\frac {1}{2}(n^2-1)-\frac {1}{2}{{\mathrm {dim}}  }  G/C}\langle S,f\rangle \] 
Note that
\[ \frac {1}{2}(n^2-1)-\frac {1}{2}{{\mathrm {dim}}  }  G/C=\frac {1}{2}({{\mathrm {dim}}  }  C-1)\] 
Finally we have to take care of the involution $\sigma $. We know that there exists a linear map $s$ of $V$ into $V^*$ such that $^tX=sXs^{-1}$. If $T$ is $\sigma -$skew then a trivial computation gives for $S$ the condition
\[ \langle S,f(s^{-1}v^*,sv)\rangle =-\langle S,f\rangle \] 
To summarize we have to show that, for any fixed nilpotent matrix $X$, a distribution 
$S$ on $V\oplus V^*$ which satisfies all the following conditions is 0.
\parindent 0pt

(0) The distribution $S$ is invariant under $C$,

(1) For $t\in {\mathbb F}^*$
\[   \langle S,f(\delta (t)v,^t\delta (t)^{-1}v^*)\rangle =|t|^{\frac {1}{2}({{\mathrm {dim}}  }C-1)}  \langle S,f\rangle \] 
(2) If $P$ is a polynomial with non zero constant term then
\[ \langle S,f(P(X)gv,{}^tP(X){\,}^tg^{-1}v^*)\rangle =|{{\mathrm {Det}}   }  P(X)|^{-1}\langle S,f(gv,^tg^{-1}v^*)\rangle \] 
(3) The distribution $S$ is equal to its Fourier transform

(4) The support of $S$ is contained in the singular set:
\[ \Sigma_ X=\{(v,v^*)\, |\, {\mathrm  for\, all\, p}\,\, \langle v^*,X^pv\rangle =0\}\] 
(5) The distribution $S$ is "skew".

We proceed by induction on $n$. We may assume that $n$ is odd; the case $n=1$ is done in section 4. In this case the condition (0) is enough to conclude.
\section{\hglue -12pt A descent method}
We keep the same notations and we want first to collect some informations on the centralizer $C$ of $X$.
\subsection{ Jordan types and the structure of $C$}
In the nilpotent case a Jordan block of length $r$ is the matrix of size $r$:
\[ n_r=\pmatrix{0&1&0&\dots&0&0\cr
0&0&1&\dots&0&0\cr
0&0&0&\dots&0&0\cr
0&0&0&\dots&1&0\cr
0&0&0&\dots&0&1\cr
0&0&0&\dots&0&0\cr}\] 
In a suitable basis $X$ has a matrix which is a "direct sum" of Jordan blocks. The basis is not unique but the set of length of the blocks is. This defines the Jordan type of $X$. For example we shall say that $X$ has Jordan type $(1,2,2)$ if in some basis the matrix of $X$ is the "direct sum" of one block of length 1 and two of length 2.

In general suppose that  $V=V_1\oplus V_2$ with $X(V_i)\subset V_i$. Call $X_i$ the restriction of $X$ to $V_i$. Assume also that $X_2$   is the sum of $m$ Jordan blocks of the same size $r$. Then we may identify $V_2$ with some vector space $W^r$ with $W$ a vector space of dimension $m$ in such a way that if
$v_2=(w_1,\dots ,w_r)$ is some element of $V_2$, then $X(v_2)=(w_2,\dots ,w_r,0)$. Thus the matrix of $X_2$ is $n_r$ where the "1"
 now represent the identity map of $W$. Using the decomposition $V=V_1\oplus V_2$ we have
 \[ X=\pmatrix{X_1&0\cr 0&n_r\cr}\] 
 Let 
 \[ c=\pmatrix{\alpha &\beta \cr \gamma &\delta \cr}\] 
 be an endomorphism of $V$. It will commute with $X$ if and only if
 \[ \alpha X_1=X_1\alpha ,\,\, \delta n_r=n_r\delta ,\,\, X_1\beta =\beta n_r,\,\, \gamma X_1=n_r\gamma  \] 
 The second condition means that $\delta $ may be written as
 \[ \delta =\pmatrix{d_1&d_2&\dots&d_r\cr
 0&d_1&\dots&d_{r-1}\cr
 \dots&\dots&\dots&\dots\cr
 0&0&0&d_1\cr}\] 
 with the $d_i$ endomorphisms of $W$.
 
Let us analyse the condition on $\beta $. This map is a linear map from $W^{r}$ into $V_1$; call $u_1,\dots ,u_{r}$ its columns which are elements of $V_1$. The condition is
\[ X_1u_j=u_{j-1},\,\,\,\, {\mathrm  for}\,\,\, 2\leq j\leq r,\quad X_1u_1=0 \] 
To get such a $\beta $ we choose $u\in V_1$ and put $u_j=X_1^{r-j}u$. We must have $X_1u_1=0$ that is to say $X_1^{r}u=0$.

The situation for $\gamma $ is similar except that we consider the rows and not the columns.

From now on suppose $X_1^{r-1}=0$. In other words $r$ is the largest possible size for a Jordan block of $X$ and $X_2$ is the sum of all the Jordan blocks of size $r$. In this case we have $X_1^{r}u=0$ for any choice of $u$ and in fact $u_1=0$. Similarly the first row of $\gamma $ is arbitrary and the last row is always 0.
 \begin{proposition}
 The determinant of $c$ is
\[ {{\mathrm {Det}}   }  c={{\mathrm {Det}}   }  (\alpha ){{\mathrm {Det}}   } (d_1)^r\] 
\end{proposition}
The proof is left to the reader. One can for example introduce a full Jordan basis for $X$, regroup the blocks of the same size and then proceed by induction on the number of sizes appearing in the Jordan type of $X$. In particular $c$ is invertible, hence belongs to $C$, if and only if $d_1$ and $\alpha $ are invertible.

An element  of $V$ is written as $(v_1,w_1,\dots ,w_r)$
\begin{proposition}
The open set $\{v\, |\, w_r\ne 0\}$ is an open orbit  $\Omega  $ of $C$ in $V$.
\end{proposition}
This is clear from the description of the matrix of $c$.

The involution $\sigma $ is defined using a basis $(e_i)$ of $V$, the dual basis $(e^*_i)$ and the linear map $u$ from $V$ to $V^*$ such that $u(e_i)=e^*_i$. Then
\[ \sigma (X,v,v^*)=(u^{-1}{}^tXu,u^{-1}v^*,u(v))\] 
Let $w_ X$ be a bijective linear map of $V$ onto $V^*$ such that $^tX=w_ X\, X\, w_ X^{-1}$ and put $s=uw_ X^{-1}u$. To be precise we choose a Jordan basis for $X$ and if $\varepsilon _1,\dots ,\varepsilon _p$ is the part of the basis corresponding to one block and $\varepsilon _1^*,\dots ,\varepsilon _p^*$ the dual basis we take $u(\varepsilon _i)=\varepsilon _i^*$ and $w_ X(\varepsilon _i)=\varepsilon _{p+1-i}^*$. Thus we may assume $s=w_ X$ and $s=^ts$. A trivial computation shows that $T$ is skew if and only if
\[ \langle S,f(s^{-1}v^*,sv)\rangle =-\langle S,f\rangle \] 

\subsection{ Restriction to the open orbit}
 Let ${\mathcal  O}$ be an orbit of $C$ in $V$. Let $S$ be a distribution on $\bigl ({\mathcal  O}\times V^*\bigr )$ which is invariant under $C$ and with support contained in $\Sigma_X$. Note that the orbit ${\mathcal  O}$ is stable by $P(X)$ for any polynomial $P$ with non zero constant term because $P(X)\in C$. Thus it make sense to assume that $S$ satisfies condition $(2)$. We are going to prove that for some ${\mathcal  O}$, including the open one $\Omega  $, such a distribution must be 0.

We fix a point $e\in {\mathcal  O}$ and let $D$ be the centralizer of $e$ in $C$. Put
\[ \Sigma_{X,e}=\{v^*|(X,e,v^*)\in \Sigma\}\] 
Then $\Sigma_{X,e}$ is the subspace of $V^*$ orthogonal to the subspace $E(X,e)$ of $V$
generated by $X$ and $e$. There exists a distribution $U$ on $\Sigma_{X.e}$ such that
\[ \langle U,\alpha (^t\hskip -1pt d^{-1}v^*)\rangle =\Delta _ D(d)\langle U,\alpha (v^*)\rangle \] 
and
\[ \langle S,\psi (v,v^*)\rangle =\oint _{C/D}\langle U,\psi (ce,^t\hskip -1pt c^{-1}v^*)\rangle d\mu (c)\] 
We first use the homogeneity condition. As $S$ is $C-$invariant we have, for $P$ a polynomial with non zero constant term 
\[ \langle S,\psi (P(X)^{-1}v,P(^tX)v^*\rangle =\langle S,\psi (v,v^*)\rangle \] 
Hence
\[ \langle S,\psi (v,P(^tX)^2v^*)\rangle =\langle \langle S,\psi (P(X)v,P(^tX)v^*)\rangle =|{{\mathrm {Det}}   }  P(X)|^{-1}\langle R,\psi (v,v^*)\rangle \] 
which implies
\[ \langle U,\varphi (P(^tX)^2v^*)\rangle =|{\mathrm {Det} }  P(X)|^{-1}\langle U,\varphi \rangle ,\quad \varphi \in \mathcal  S(S_{X,e})\] 
Note also that ${\mathrm {Det} }  P(X)=P(0)^n$.For our present purpose it will be enough to take for $P$ the constant polynomial $P(X)=t$.

Next we study the invariance condition under  $D$.  

We shall need another assumption. we suppose that there exists a subspace $W$
of $V$, which is invariant by $X$ and such that $V=W\oplus E(X,e)$. Put $E=E(X,e)$ and
consider the basis $\varepsilon _1,\dots ,\varepsilon _s$ of $E$ defined by $\varepsilon _{s-j}=X^je$, with $s={{\mathrm {dim}}  } 
E$. Then we may represent $X$ as \[ X=\pmatrix{Y&0\cr 0&n_s
\cr}\] 

If $d\in D$ then
\[ d=\pmatrix{a&0\cr b&1
\cr}\] 
with $aY=Ya$ and $n_sb=bY$. In particular
\[ d_t=\pmatrix{t{{\mathrm {Id}}  }_ W&0\cr 0&1
\cr}\in D\] 
hence
\[ \langle U,\varphi (t^{-1}v^*)\rangle =\Delta _ D(d_t)\langle U,\varphi \rangle \] 
We conclude that
\[ \Delta _ D(d_t)^2=n\] 
The equation $n_sb-bY=0$ is easy to solve. The condition is a condition on the
rows of $b$, relative to the basis $\varepsilon _j^*$ of $E^*$, dual to the basis $\varepsilon _j$. We
choose $w^*\in W^*$ such that $^tY^sw^*=0$ and define $b$ by
\[ \langle \varepsilon _j^*,b(w)\rangle =\langle ^tY^{j-1}(w^*),w\rangle \] 
so that we obtain a commutative group isomorphic to ${\mathrm {Ker}}  (^tY^s)$.

 To compute $\Delta _ D$ we recall that 
\[ \Delta _ D(d)=|{{\mathrm {Det}}   }  ({{\mathrm {ad}}  }  (d))|^{-1}\] 
with ${{\mathrm {ad}}  } $ the adjoint action on the Lie algebra ${\mathfrak  d}$ of $D$. Taking $d=d_t$
we deduce that
\[ \Delta _ D(d_t)=|t|^{{{\mathrm {dim}}  } ({\mathrm {Ker}}  ^tY^s})\] 
Finally, if the distribution $S$ is non zero then 
\[ n=2{{\mathrm {dim}}  } ({\mathrm {Ker}}  ^tY^s)\] 
is even contrary to one of our hypothesis.

The condition that $E=E(X,e)$ is a direct factor for $V$ viewed an ${\mathbb 
F}[X]-$module is rather strong but is satisfied at least for the open orbit of $C$ in
$V$.

 Indeed start with a Jordan decomposition of $X$:
\[ X=\pmatrix{n_{r_1}&&&&\cr&n_{r_2}&&\cr
&&n_{r_3}&&\cr
&&&\dots&\cr
&&&&n_{r_s}
\cr},\quad r_1\leq r_2\leq \cdots \leq r_s\] 
Put $r_s=r$ and take for $e$ the last vector of the basis. Then $E(X,e)$
corresponds to the last block and $W$ to the others. The orbit of $e$ is
$V\setminus {\mathrm {Ker}}  X^{r-1}$, the open one.

The conclusion is that we may add to our assumption on $S$ that its support is
contained in ${\mathrm {Ker}} (X^{r-1})\oplus V^*$. However we can reverse the role of $V$ and
$V^*$ and conclude that it is also contained in $V\oplus {\mathrm {Ker}}  (^tX^{r-1})$ so that it is
finally contained in ${\mathrm {Ker}}  (X^{r-1})\oplus {\mathrm {Ker}}  (^tX^{r-1})$
Furthermore  $S=\widehat S$ (Fourier transform on $V\oplus V^*$)
Thus the distribution $S$ is invariant by
translations by elements of the subspace orthogonal to  ${\mathrm {Ker}}  (X^{r-1})\oplus {\mathrm {Ker}} 
(^tX^{r-1})$, that is to say  ${{\mathrm {Im}}   }(X^{r-1})\oplus {{\mathrm {Im}}  } (^tX^{r-1})$. Note that,
for $r\geq 2$
 \[ {{\mathrm {Im}}   }(X^{r-1})\oplus {{\mathrm {Im}}  } (^tX^{r-1})\subset {\mathrm {Ker}}  (X^{r-1})\oplus {\mathrm {Ker}} 
(^tX^{r-1})\] 

{\textbf  Remark\pointir} If $r=1$, so $X=0$ and $C=GL(V)$, then the distribution $S$ has
its support contained in $\{0\}\times \{0\}$ hence is a multiple of the Dirac
measure.and must be equal to its Fourier transform, so it is 0. We may assume
that $r\geq 2$.

\subsection{The descent}
Assume $X\ne 0$ . Let $V'={\mathrm {Ker}}  X^{r-1}/{{\mathrm {Im}}  }X^{r-1}$. Then the dual space $V^*$ is ${\mathrm {Ker}} \, ^tX^{r-1}/{{\mathrm {Im}}  }\,^tX^{r-1}$. The operator $X$ is 0 on ${{\mathrm {Im}}  }X^{r-1}$ and $X({\mathrm {Ker}} \, X^{r-1}\subset {\mathrm {Ker}} \, X^{r-1}$ so that $X$ defines a nilpotent endomorphism $X'$ of $V'$. Of course $^tX$ will give the transpose of $X'$. An element $c$ of $C$ will also define an endomorphism of $V'$. This gives an homomorphism $\theta $ of $C$ into the commutant $C'$ of $X'$. However $\theta $ is {\textbf  not onto}.

 To be more explicit let us use the setup of 9-1. So $V=V_1\oplus W^r$ and an element of $V$ is written as $(v_1,w_1,\dots ,w_r)$. The kernel of $X^{r-1}$ is given by $w_r=0$ and the image of $X^{r-1}$ by $v_1=0,w_2=w_3=\dots =w_r=0$. We identify $V'$ with the subspace $w_1=w_r=0$. Then one goes from $V$ to $V'$ by forgetting the components relative to $w_1$ and $w_r$. For example
 \[ X=\pmatrix{X_1&0\cr 0&n_r\cr},\quad X'=\pmatrix{X_1&0\cr 0&n_{r-2}\cr}\] 
 If
 \[ c=\pmatrix{\alpha &\beta \cr \gamma &\delta \cr}\quad {\mathrm  and}\,\, c'=\theta (c)=\pmatrix{\alpha '&\beta '\cr \gamma '&\delta '\cr}\] 
then $\alpha =\alpha '$ and $\delta '$ is obtained from $\delta $ by deleting the first and lasr row and columns so that if $\delta $ is given by $D_1,\dots ,d_r$ then $\delta '$ is given by $d_1,\dots ,d_{r-2}$. 

The columns $u_j$ of $\beta $, which are elements of $V_1$ are given by $u_j=X_1^{r-j}u$ with $u$ arbitrary in $V_1$. Then $\beta '=(X_1^{r-2}u,\dots ,X_1u)$. For $C'$ a typical element $\beta '$ is $\beta '=(X_1^{r-3}u',\dots ,u')$ where $u'$ must be such that $X_1^{r-2}u'=0$ ( this is not automatic because, for some $X$, we may have $X_1^{r-2}\ne 0$). Therefore if we want $\theta $ to be onto a necessary condition is that ${\mathrm {Ker}}  X_1^{r-2}$ is contained in the image of $X_1$. We know that $X_1^{r-1}=0$ so that any Jordan block of $X_1$ has a size at most equal to $r-1$. For a block of size strictly smaller than $r-1$ the contribution to the kernel of $X^{r-2}$ is the full corresponding subspace which is not contained into the image. Thus the condition is satisfied if and only if $X_1$ contains only blocks of size $r-1$. We have to cases: either $X$ has only blocks of size $r$ or  the only possible sizes are $r$ and $r-1$.

We must also check $\gamma $; this gives the same condition.

 We now go back to $S$. From the result of the last subsection there exists a distribution $S'$ on $V'$ such that, with obvious notations
 \[ S=S'\otimes dw_1\otimes \delta _{w_1}\otimes \delta _{w^*_1}\otimes dw_r^*\] 
 The Haar measures are normalized as the Fourier transform of the Dirac measures.
 
 We try to transfer the conditions on $S$ to conditions on $S'$ and to use induction.
It is obvious that if $S$ satisfies conditions (3),(4) or (5) then the same is true for $S'$. The same is true for condition (2). This can be seen by a direct computation (note that ${{\mathrm {Det}}   }  P(X)=P(0)^n$ or by remarking that it is a consequence of conditions (3) and (4).

Unfortunately in general $S'$ is not invariant under the full commutant of $X'$ but only under $\theta (C)$.
 In the first exceptionnal case where all the Jordan blocks have the same size, which must be odd because $n$ is odd, we go down from size $r$ to the size $r-2$, keeping the full invariance and we can iterate until we arrive at size 0 which means $X=0$ and we can conclude that $S=0$.
 
 In the second case, starting from blocks of size $r$ and $r-1$ we go down to $r-1$ and $r-2$ and go on down to $X=0$ which again gives $S=0$.

  In particular suppose that $n=3$. Then the possible Jordan types for $X$ are
  \[ (1,1,1),\,\, (1,2),\,\, (3)\] 
  and all three cases are settled. Therefore
 \begin{proposition}
 Conjecture 2 is true for $n\leq 4$
 \end{proposition}
We will try to do a little better. 

Let us first deal with the behavior of condition (1) under descent.
 We have to choose $\delta (t)$. On $V_2$ we define it by
\[ \delta (t)(w_1,\dots ,w_r)=(t^{r-1}w_1,t^{r-2}w_2,\dots w_r)\] 
and make any acceptable choice on $V_1$.
The restriction to $V'$ is a possible choice for a one parameter group for $X'$. Condition $(1) $ for $S$ and the relation between $S'$ and $S$ immediately implies that
\[ \langle S',f(\delta (t)u,^t\delta (t)u^*)\rangle =|t|^{m(r-1)+\frac {1}{2}({{\mathrm {dim}}  }  C-1)}\langle R,f\rangle \] 
The degree of homogeneity is not the one we would want for $S'$, that is to say $\frac {1}{2}({{\mathrm {dim}}  }  C' -1)$. There is a shift and this has to be kept in mind when using condition (1). 
\section{\hglue -12pt Some particular cases}
To use the descent we must first show that the distribution is 0 on ${\mathrm {Ker}}  X^{r-1}\times V^*$. If we have invariance under the full group $C$ then ${\mathrm {Ker}}  X^{r-1}$ is an orbit in $C$. However after going down once we will have invariance only under a subgroup $C_1$ of $C$ and ${\mathrm {Ker}}  X^{r-1}$ is the union of a finite number of orbits. For each such orbit $\gamma $ we must prove that a distribution on $\gamma \times V^*$, invariant under $C_1$, having singular support and satisfying the homogeneities conditions is 0. Also, in condition (1)
 we have to take the shift into account. This can be done in a limited number of low dimensional cases. 
 
 We use the setup of section 2-1. We call $C_0$ the subgroup of $C$ given by the conditions $\beta =0,\, \gamma =0$ and we shall always assume that we have invariance under a subgroup $C_1$ of $C$ containing $C_0$. The one parameter group $\delta (t)$ will always be choosen inside $C_0$ so that the precise choice will be irrelevant. Finally we shall take condition $(1)$ as
 \[ \langle S,f(\delta (t)v,^\delta (t)^{-1}v^*\rangle =|t|^\mu \langle S,f\rangle \] 
 leaving the parameter  $\mu $ free. Then for a given $C_1$ we try to prove that $S$ must be 0 outside ${\mathrm {Ker}}  X^{r-1}\times V^*$. If we use condition (1) then some values of $\mu $ will have to be excluded.

   \subsection{A simple remark for $X=0$}
If $n=1$ then any distribution on $V\oplus V^*$ which is invariant under $GL(V)$ is symmetric with respect to $\sigma $. Explicitly we identify $V$ and $V^*$ to ${\mathbb  F}$ and $GL(V)$ to ${\mathbb  F}^*$. The action is
\[ \bigl (t,(x,y)\bigr )\mapsto (tx,t^{-1}y)\] 
and $((x,y))=(y,x)$.  

Take any $n$ and $X=0$. Fix a basis $(e_1,\dots e_n)$ of $V$; call $(e^*_1,\dots ,e^*_n)$ the dual basis. Let $(x_1,\dots ,x_n)$ be the coordinates in $V$ and $(y_1,\dots ,y_n)$ the \hfill\break
coordinates in $V^*$.
\begin{lemma}
 Let $S$ be a distribution on $V\oplus V^*$, invariant under the action of $({\mathbb  F}^*)^n$ acting by
\[ \bigl ((t_1,\dots ,t_n),(x_1,\dots ,x_n,y_1,\dots ,y_n)\bigr )\mapsto \bigl (t_1x_1,\dots ,t_nx_n,t_1^{-1}y_1,\dots,t_n^{-1}y_n\bigr )\] 
Then $S$ is invariant by the involution
\[ \sigma :\,\,\, (x_1,\dots ,x_n,y_1,\dots,y_n)\mapsto (y_1,\dots ,y_n,x_1,\dots ,x_n)\] 
\end{lemma}

We have to prove that for any $f\in {\mathcal  S}(V\oplus V^*)$ we have $\langle S,\sigma (f)\rangle =\langle S,f\rangle $. It is enough to do it for a product
\[ f(x,y)=\prod _1^nf_i(x_i,y_i)\] 
If we freeze all the $f_i$ except one, say $f_j$ then we get a distribution on ${\mathbb  F}^2$, invariant under the action of ${\mathbb  F}^*$. Therefore we can permute $x_j$ and $y_j$. This being true for all $j$ we conclude that $S$ is invariant by $\sigma $.

 Take $X$ of Jordan type $(1,\dots ,1,3\dots ,3)$ ( with $p$ times 1 and $q$ times 3). Apply the descent method : $X'=0$. The distribution $S'$ is not invariant under the full $GL(p+q)$ but is invariant under the subgroup $GL(p)\times GL(q)$. We may apply the Lemma to conclude that $S'=0$. 
In particular, for $n=5$ this takes care of the case $(1,1,3)$ and for $n=7$ of the case $(1,1,1,1,3)$.
\subsection{The $(1,\dots ,1,r)$ case}
We are now going to use  condition (1). Suppose that $X$ is principal and choose a basis of $V$ such that, in matrix form $X=n_r$. Let
\[ \delta (t)=\pmatrix{t^{r-1}&0&\dots &0\cr
0&t^{r-2}&\dots &0\cr
\dots&\dots&\dots&\dots\cr
0&0&\dots&1\cr}\] 
let $x_1,\dots ,x_r$ be the coordinates in $V$.
\begin{lemma}
Let $S$ be a distribution on $V\oplus V^*$, invariant under $C$,  supported by the singular set and which satisfies  condition (1). If $\mu \ne -r(r-1)/2$ then the support of $S$ is contained in $\{v|x_r=0\}\times V^*$.
\end{lemma}
Indeed $\{v\in V|x_r\ne 0$ is the open orbit of $C$ in $V$. We may restrict the distribution $S$ to
$\{v|x_r=0\}\times V^*$. The centralizer of $e_r$ being trivial, there exists a distribution $R$ on $V^*$ such that
\[ \langle S,f\rangle =\int _C\langle R,f(ce_r,^tc^{-1}v^*\rangle dc\] 
However the subspace generated by the vectors $X^qe_r$ is $V$ and the support of $S$ is contained in the singular set so that $R$ must be some multiple of the Dirac measure $\varphi \mapsto \varphi (0)$ on $V^*$. Therefore, for some scalar $\lambda \in {\mathbb  C}$, the restriction of $S$ is
\[ \langle S,f\rangle =\int _Cf(ce_r,0)dc=\lambda \int _{\mathbb  F^r}f(\sum x_ie_i,0)\,\, |x_r|^{-r}dx_1\dots dx_r\] 
and we get
\begin{eqnarray*} \langle S,f(\delta (t)v,^t\delta (t)^{-1}v^*\rangle &=&\lambda \int _{\mathbb  F^r}f(\sum t^{r-i}x_ie_i,0)\,\, |x_r|^{-r}dx_1\dots dx_r\\
&=&|t|^{-(1+2+\cdots +r-1)}\langle R,f\rangle =|t|^{-r(r-1)/2}\langle R,f\rangle \\
\end{eqnarray*}
This implies the Lemma.

{\textbf  Remark\pointir}The proof remains valid for a distribution defined only on $x_r\ne 0$. The conclusion is then that this distribution is 0.

Now we take $X$ of Jordan type $(1,\dots ,1,r)$ with $r\geq 2$ and $p$ times 1:
\[ X=\pmatrix{0&0\cr 0&n_r\cr}\] 
Let $C_0$ be the "diagonal" subgroup of $C$:
\[ \pmatrix{a&0\cr 0&b\cr},\,\,\, a\in GL(p,{\mathbb  F}),\,\, bn_r=n_rb\] 
Define
\[ \delta (t)=\pmatrix{1&0&\dots&0&0&\dots&\dots&0\cr
0&1&\dots&0&0&\dots&\dots&0\cr
\dots&\dots&\dots&\dots&\dots&\dots&\dots&\dots\cr
0&0&\dots&1&0&\dots&\dots&0\cr
0&0&0&0& t^{r-1}&\dots&\dots &0\cr
0&0&0&0&0&t^{r-2}&\dots &0\cr
\dots&\dots&\dots&\dots&\dots&\dots&\dots&\dots\cr
0&0&0&0&0&\dots&\dots&1\cr}\] 
Call $e_1,\dots ,e_{p+r}$ the basis of $V$ and $e^*_1,\dots ,e^*_{p+r}$ the dual basis. Let $x_1,\dots ,x_{p+r}$ be the coordinates in $V$ and $y_1,\dots ,y_{p+r}$ the coordinates in $V^*$.
\begin{lemma}
 Let $S$ be a distribution on $V\oplus V^*$, invariant under $C_0$,  supported by the singular set and which satisfies condition (1). If $\mu \ne -r(r-1)/2$ then the support of $S$ is contained in $\{v|x_{p+r}=0\}\times V^*$.
\end{lemma}
Put $V=V_i\oplus V_2$ with $v_1=\oplus _1^p{\mathbb  F}e_i$ and $V_2=\oplus _{p+1}^{p+r}{\mathbb  F}e_i$.
The $C_0-$orbit of $e_1+e_{p+r}$ is the open set in $V$ defined by $v_1\ne 0$ and $x_{p+r}\ne 0$. The subspace generated be $X$ and $e_1+e_{p+r}$ is ${\mathbb  F}(e_{p+r}+e_1)\oplus _{p+1}^{p+r-1}{\mathbb  F}e_j$. Let $W$ be the subspace of $V^*$ orthogonal to this subspace of $V$; then
\[ W=\{(y_1,y_2,\dots ,y_p,0,0\dots ,-y_1)\}\] 
Let $M$ be the centralizer of $e_1+e_{p+r}$. The group $M$ is isomorphic to the centralizer of $e_1$ in $GL(p,{\mathbb  F})$. In particular it is not a unimodular group. On the other hand $C_0$ is unimodular. We have an integration formula
\[ \int _{C_0}f(c)dc=\oint _{C_0/M}f(cm)d_rm\, d\mu (c)\] 
Here $d_rm$ is a right Haar measure of $M$ and $d\mu $ is a $C_0-$invariant linear form on the space of functions $\varphi $ defined on $C_0$ and such that $\varphi (cm)=\Delta _ M(m)\varphi (c)$.

Consider the restriction of $S$ to $C_0(e_1+e_{p+r})\times V^*$. There exists a distribution $R$ on $W$ such that
\[ \langle R,\psi (^tm^{-1}w)\rangle =\Delta _ M(m)\langle R,\psi \rangle \] 
and
\[ \langle S,f\rangle =\oint_{C_0/M}\langle R,f(c(e_1+e_{p+r}),^tc^{-1}w)\rangle d\mu (c)\] 
We have $\delta (t)(e_1+e_{p+r})=e_1+e_{p+r}$ and $^t\delta (t)$ is the identity on $W$ and commutes with $M$, therefore
 \begin{eqnarray*}
  \langle S,f(\delta (t)v,^t\delta (t)^{-1}v^*)\rangle &=&\oint_{C_0/M}\langle R,f\bigl (\delta (t)c\delta (t)^{-1}(e_1+e_{p+r}),^t(\delta (t)c\delta (t)^{-1}w\bigr )\rangle d\mu (c)\\
&=&\frac {d\mu (\delta (t)^{-1}c\delta (t))}{d\mu (c)}\langle R,f\rangle \\
\end{eqnarray*}
and
\[ \frac {d\mu (\delta (t)^{-1}c\delta (t))}{d\mu (c)}=\frac {d(\delta (t)^{-1}c\delta (t))}{dc}\] 
A typical element of $C_0$ is
\[ \pmatrix{a&0&0&\dots&0\cr
0&b_1&b_2&\dots&b_r\cr
0&0&b_1&\dots&b_{r-1}\cr
\dots&\dots&\dots&\dots&\dots\cr
0&0&0&0&b_1\cr}\] 
The Haar measure is
\[ dc=|b_1|^{-r}dadb_1\dots db_r\] 
and the maps $c\mapsto \delta (t)^{-1}c\delta (t)$ corresponds to
\[ (a,b_1,\dots ,b_r)\mapsto (a,b_1,t^{-1}b_2,\dots ,t^{-r+1}b_r)\] 
Therefore
\[ \frac {d\mu (\delta (t)^{-1}c\delta (t))}{d\mu (c)}=|t|^{-r(r-1)/2}\] 
It follows that the restriction of $S$ is 0. The support of $S$ is thus contained into the closed subset
\[ {\mathcal  F}=\{v\in V |\, x_{p+r}=0\}à\{v=v_1+v_2\in V\, | \, v_1=0\}\] 
Let us restrict $S$ to the open subset $x_{p+r}\ne 0$ of ${\mathcal  F}$. On this subset we have $v_1=0$. Let $f\in {\mathcal  S}(v_i)$ and $f^*_i\in {\mathcal  S}(V_i^*)$. Then, on $V_1^*$ the distribution
\[ f_1^*\mapsto \langle S,f(v)f_1^*(v_1^*)f_2^*(v_2^*)\rangle \] 
is invariant under the action of the subgroup $GL(V_1)\approx GL(V_1^*)$ of $C_0$.This distribution must be a multiple of the Dirac measure at the origin and the support of the restriction of $S$  is contained into $V_2\oplus V_2^*$. If we view this restriction as a distribution on the open subset $x_{p+r}\ne 0$ of $V_2\oplus V_2^*$ we are exactly in the situation of the remark following the proof of Lemma 10-2. We conclude that the restriction of $S$ to the subset $x_{p+r}\ne 0$ of ${\mathcal  F}$ is 0 which is exactly what we claimed.

Now let us start again with $X$ of type $(1,\dots ,1,r)$ and keep the same notations. The starting value of $\mu $ is
\[ \mu =\frac {1}{2}\bigl ((p+1)^2+r-2\bigr )\] 
The support of the distribution $S$ is contained into the set $x_{p+r}=0$.We now go down one step arriving at the type $(1,\dots ,1,r-2)$.

The new value of $\mu $ is 
 \[ \mu '= \mu +r-1=\frac {1}{2}((p+1)^2+r-2)+r-1=\frac {1}{2}((p+1)^2+r-4)+r\] 
 We have invariance under the group $C_0$. As $\mu '$ is 0 it is certainly not equal to the exceptionnal value and we can go down one more time. The value of $\mu $ will keep increasing, thus staying positive and we will keep invariance under $C_0$.

 Eventually we will get the 0 matrix, of size $p$ or $p+1$, keeping all the way invariance under at least the diagonal group and skew invariance . This implies that the distribution is 0.
  
 In particular it works for the type $(1,4)$ which was the missing case for $n=5$.
\begin{proposition}
 The general multiplicity one conjecture for $GL(n)$ is true for $n\leq 6$.
 \end{proposition}
 (this means up to restriction from $GL(7)$ to  $GL(6)$).
 
\subsection{The case $n=7$}
The remaining types are:
\[ (1,1,2,3),\,\,\, ((1,2,4)\,\,\, (2,5)\] 

We shall fix a Jordan basis for $X$, say $(e_1,\dots ,e_7)$, the Jordan blocks being taken by increasing size. The dual basis is $(e_1^*,\dots ,e_7^*)$. The coordinates in $V$ are $(x_1,\dots ,x_7)$ and in $V^*$ they are $(y_1,\dots ,y_7)$. The involution is built with the map $s: e_i\mapsto e_(8-i)^*$.

Let us start with the $(1,1,2,3)$ case. The distribution $S$ satisfies all our usual conditions with the original value of $\mu $ and invariance under the full centralizer $C$. Apply the descent method, in the original form, that is to say for the orbit of $e_7$. Then we go down to the case $(1,1,2,1)$. The distribution $S_1$ that we get is invariant under the subgroup $C_1$ of $C$ described by the picture
\[ \pmatrix{*&*&0&*&0\cr
*&*&0&*&0\cr
*&*&*&*&*\cr
0&0&0&*&0\cr
0&0&0&*&*\cr}\] 
Keeping the original numbering of the basis vector we see that $C_1e_4=\{x_4\ne 0\}$ is an open orbit in the space $W={\mathbb  F}e_1\oplus {\mathbb  F}e_2\oplus {\mathbb  F}e_3\oplus {\mathbb  F}e_4\oplus {\mathbb  F}e_6$. We restrict $S_1$ to $C_1e_4\times W^*$. we have $Xe_4=e_3,\,\, X^2e_4=0$. Let $M$ be the isotropy subgroup of $e_4$. Then a generic element of $C_1$ is
\[ c_1=\pmatrix{ a_1&a_2&0&b_1&0\cr
a_3&a_4&0&b_2&0\cr
c_1&c_2&d_1&d_2&b_3\cr
0&0&0&d_1&0\cr
0&0&0&c_3&d_3\cr}\] 
with $a_1a_4-a_2a_3\ne 0,\,\, d_1\ne 0,\,\, d_3\ne 0$.
The reductive part of $C_1$ is obtained by taking $b_1=b_2=c_1=c_2=d_2=b_3=c_3=0$ and the unipotent radical by $a_1=a_2=d_1=d_3=1,\,\, a_2=a_3=0$. Let
\[ a=\pmatrix{a_1&a_2\cr a_3&a_4\cr}\] 
The action of the reductive part onto the unipotent radical is given by
 \begin{eqnarray*}
  (c_1,c_2)&\mapsto d_1(c_1,c_2)a^{-1}\\
\pmatrix{b_1\cr b_2\cr}&\mapsto a\pmatrix{b_1\cr b_2\cr}d_1^{-1}\\
d_2&\mapsto d_2\\
c_3&\mapsto c_3d_3/d_1\\
b_3&\mapsto b_3d_1/d_3\\
\end{eqnarray*}
The Haar measure of the unipotent radical is $db_1db_2db_3dc_1dc_2dc_3d(d_2)$ and is invariant under the action of the reductive part. It follows that $C_1$ is unimodular.

Now we look at $M$ so we take $b_1=b_2=d_2=c_3=0$ and $d_1=1$. The reductive part is obtained by adding the conditions $b_3=c_1=c_2=0$ while the unipotent radical corresponds to $a=0,\, d_3=1$. The adjoint action of the reductive part onto the unipotent radical is
\[ b_3\mapsto b_3/d_3,\,\, (c_1,c_2)\mapsto (c_1,c_2)a^{-1}\] 
The Haar measure of this radical is $du=db_3dc_1dc_2$ hence, for $h$ in the reductive part $H$
\[ \frac {d(huh^{-1})}{du}=|d_3|^{-1}|{{\mathrm {Det}}   }  a|^{-1}\] 
Therefore the right Haar measure of $M$ is given by
\[ \int _Mf(m)d_rm=\int _{H\times U}f(au)|d_3|^{-1}|{{\mathrm {Det}}   }  a|^{-1}dhdu\] 
and the modulus of $M$ is
$\Delta _ M(hu)=|d_3||{{\mathrm {Det}}   }  a|$

We restrict our distribution $S_1$ to $C_0e_4\times W^*$. There exists a distribution $R$ on $W^*$ such that
\[ \langle R,\psi (^tm^{-1}w^*)\rangle =\Delta _ M(m)\langle R,\psi \rangle \] 
and
\[ \langle S_1,f\rangle =\oint_{C_1/M}\langle R,f(ce_4,^tc^{-1}w^*)\rangle d\mu (c)\] 
The support of $S_1$ is singular so that the support of $R$ must be contained into the subset $y_3=y_4=0$. The semi-invariance under $M$ means that
\[ \langle R,\psi \Bigl ((y_1,y_2)a,d_3y_6\Bigr )\rangle =|d_3||{{\mathrm {Det}}   }  a|\langle R,\psi \rangle \] 
In particular freezing $y_6$ we obtain a distribution $R_1$ on $(y_1,y_2)$ such that $\langle R_1,\varphi ((y_1,y_2)a)\rangle =|{{\mathrm {Det}}   }  a|\langle R_1,\varphi \rangle $. The only such distribution is 0 so that $R=0$ and finally the restriction of $S_1$ is 0.

We could also use the homogeneity condition
 \[ \langle R,\varphi (t^2w^*)\rangle =|t|^{-5}\langle R,\psi \rangle  \] 
 Now in $M$ we choose $a=t{{\mathrm {Id}}  }$ and $d_3=t$. By the semi-invariance
 \[ \langle R\psi \Bigl (t(y_1,y_2),ty_6\Bigr )\rangle =|t|^3\langle R,\psi \rangle \] 
 {\textbf  Remark\pointir}The point here is that there exists an element of the centralizer $M$ of the base point $e$ in $V$ which is multiplication by $t$ on $E(X,v)^\perp$. This was already the case in our original descent method. It is unfortunately not true in general but could be true in a number of cases.
 
 The support of $S_1$ is contained into the subset $x_4=0$. Then using again the Fourier transform and $\sigma $ we can go down to the Jordan type $(1,1,1)$ and an invariance groupe $C_0$ isomorphic to $GL(2)\times GL(1)$. This is enough to conclude that our new distribution $S_2$ is 0 and consequently that our original distribution $S$ is 0.
   
   The other 2 cases may be settled in a similar way We omit the details.
   
   Then:

  \begin{theorem}
 The general multiplicity one conjecture for $GL(n)$ is true for $n\leq 8$.
 \end{theorem}
 (this means up to restriction from $GL(9)$ to  $GL(8)$).

Using the same kind of elementary computations  more cases can be done. For example if there is only two possible sizes for the Jordan blocks ( with or without multiplicity) and if we assume invariance under $C_0$, then except for some values of $\mu $ the descent will be justified. If we take $n=9$ this takes care of all the types with two possible sizes and we are left with 7 cases like $(1,2,6),(1,2,2,4)\dots)$. 

Note that we are only using invariance under some subgroup $C_1$, the singularity of the support and the homogeneity conditions. Then, in each case, if $\Gamma $ is an orbit of $C_1$ in $V$ contained in $V\setminus {\mathrm {Ker}}  X^{r-1}$ we proved that on $\Gamma \times V^*$ the distribution must be 0. However it is not difficult to show that in general such $\Gamma \times V^*$ do carry non zero distributions with these properties (this already happens with a 3 block situation and $C_1=C_0$). The problem is then to prove that such distributions can not extend to invariant skew distribution on $V\oplus V^*$. In some sense the case $n=1$ is already typical.

\section{\hglue -12pt Another approach}
On $V\oplus V^*$ we have the non degenerate quadratic form $Q(v,v^*)=\langle v^*,v\rangle $ and therefore the dual pair $O(Q)\times SL(2)$ acting on ${\mathcal  S}(V\oplus V^*)$. The distribution $S$ is invariant under $SL(2)$. Therefore we may consider it as a linear form on the covariant space which is the space of the minimal representation of the orthogonal group. The group $G$ acts on $V\oplus V^*$ and this is an imbedding of $G$ into $O(Q)$. The restriction of the minimal representation to $G$ turns out to be a finite sum of inequivalent irreducible members of some degenerate  principal serie representations. This allows us to reformulate our problem in a more representation theoretic way. 

The quadratic form $Q$ does not depend on $X$ so we  go back to our distribution $T$ which we view as a linear form on  $\mathcal  S({\mathcal  N}\oplus V\oplus V^*)$. 

Consider the partial Fourier transform relative to $V^*$:
\[ f\mapsto \int _{V^*}f(X,v_1,v^*)\tau (-\langle v^*,v_2)\, dv^*\] 
As before $\tau $ is a non trivial additive character of ${\mathbb  F}$ and the Haar measures $dv$ and $dv^*$ are dual measures.

Let $R$ be the partial Fourier transform of $T$. It has the following properties.

(0) It is Invariant under $G$:
\[ \langle R,f(gXg^{-1},gv_1,gv_2)\rangle =|{{\mathrm {Det}}   }  g|^{-1}\langle R,f\rangle \] 

(1) As a distribution on $[{\mathfrak  g},{\mathfrak  g}]\oplus V\oplus V$ it is equal to its partial Fourier transform relative to $[{\mathfrak  g},{\mathfrak  g]}$  and
\[ \langle R,f(tX,v_1,v_2)\rangle =|t|^{-(n-1)/2}\langle R,f\rangle \] 

(2) If  $A,B,C,D$ are 4 polynomials in one variable, such that $AD-BC=1$ then
\[ \langle R,f(X,A(X)v_1+B(X)v_2,C(X)v_1+D(X)v_2)\rangle =\langle R,f\rangle \] 
(see section 1-3).

For the involution we choose a basis $(e_i)$ of $V$ with dual basis $(e_I^*)$. Let $(x_i)$ be the coordinates in $V$ and $(x_i^*)$ the coordinates in $V^*$. Let $u$, a linear map from $V$ onto $V^*$ be defined by $u(e_i)=e_i^*$. Then we may take
\[ \sigma (X,v,v^*)=(u^{-1}\, {}^tXu,u^{-1}(v^*),u(v))\] 
On ${\mathbb  F}$ let $dx$ be the measure self dual with respect to $\tau $. We take $dv=dx_1\dots dx_n$ and $dv^*=dx_1^*\dots dx_n^*$. These measures are self dual and the image of $dv$ by $u$ is $dv^*$. If
\[ {\mathcal  F}(X,v^*_1,v^*_2)=\int _{V\times V}f(X,v_1,v_2)\, \tau (-\langle v_1^*,v_1\rangle -\langle v_2^*,v_2\rangle )dv_1dv_2\] 
then, after partial Fourier transform on $V^*$ the involution $\sigma $ becomes
\[ f(X,v_1,v_2)\mapsto {\mathcal  F}f(X,u(v_2),-u(v_1))\] 
Therefore $T$ is skew if and only if $R$ is skew:
\[ \langle R,{\mathcal  F}f(X,u(v_2),-u(v_1))\rangle =-\langle R,f(X,v_1,v_2)\rangle \] 
Note that (2) implies in particular that $\langle R,f(X,-v_2,v_1)\rangle =\langle R,f(X,v_1,v_2)\rangle $ so that the above condition may be rewritten as
\[ \langle R,{\mathcal  F}f(X,u(v_1),u(v_2))\rangle =-\langle R,f(X,v_1,v_2)\rangle \]

In particular the distribution $R$ is invariant under $SL(2,{\mathbb  F})$ acting on $V\oplus V$ through
\[ \pmatrix{v_1\cr v_2\cr}\mapsto \pmatrix{a&b\cr c&d\cr}\pmatrix{v_1\cr v_2}=\pmatrix {av_1+bv_2\cr cv_1+dv_2\cr}\] t
and acting trivially on ${\mathfrak  g}$.

The center $Z$ of $G$ acts trivially on ${\mathfrak  g}$ and $R$ is semi-invariant under $Z$:
\[ \langle R,f(X,\lambda v_1,\lambda v_2)\rangle =|\lambda |^{-n}\langle R,f\rangle \] 
We shall first find all linear forms $L$ on ${\mathcal  S}(V\oplus V)$ which have the above invariance.

 We stratify $V\oplus V$ using the rank of
$(v_1,v_2)$; the subset of $(v_1,v_2)$ with rank $r$ (resp at most $r$) is denoted by $(V\oplus V)_{=r}$
(resp $V\oplus V)_{\leq r}$. They are all locally closed. Then for $r=1$ and $r=2$ we have short exact
sequences
\[ 0\longrightarrow {\mathcal  S}\bigl ((V\oplus V)_{=r}\bigr )\longrightarrow {\mathcal  S}\bigl
((V\oplus V)_{\leq r}\bigr )\longrightarrow {\mathcal  S}\bigl ((V\oplus V)_{\leq r-1}\bigr )\longrightarrow 0\] 

Now the stratification is stable for the action of the group $SL(2,{{\mathbb  F} })\times Z$, so working in
the category of smooth modules we have long exact homology sequence
 \begin{eqnarray*} 
 \longrightarrow H_1\bigl ({\mathcal  S}((V\oplus V)_{\leq r-1})\bigr )&\longrightarrow H_0\bigl
({\mathcal  S}((V\oplus V)_{=r})\bigr )\longrightarrow H_0\bigl ({\mathcal  S}((V\oplus V)_{\leq r})\bigr
)\\
&\longrightarrow H_0\bigl ({\mathcal  S}((V\oplus V)_{=r-1})\bigr )\longrightarrow 0 \\
\end{eqnarray*}
The group $G$ commutes with $SL(2)\times Z$ hence the above exact sequences are exact sequences
of $G-$modules.

Let us consider first $(V\oplus V)_{=0}=\{(0,0)\}$. Then ${\mathcal  S}((V\oplus V)_{=0})={\mathbb  C}$. The group
$SL(2)$ acts trivially and the group $Z$ acts as
\[ f(0,0)\mapsto f(z^{-1}(0,0))|{{\mathrm {Det}}   }  z|^{-1}\] 
that is to say by the character $|{{\mathrm {Det}}   }  z|^{-1}$. The action of $G$ is also given by this character.

There can be no linear form with the required semi-invariance  so $H_0\bigl ((V\oplus V)_{=0}\bigr )=(0)$.

Let us next consider  $(V\oplus V)_{=1}$. If $(v_1,v_2)\in (V\oplus V)_{=1}$ then the subspace
generated by $v_1$ and $v_2$ is an element of the projective space ${\mathbb  P}(V)$. hence we have a
map $p$ from $(V\oplus V)_{=1}$ onto the projective space. We use Bernstein localization principle .
The group is $SL(2)\times Z$; each fiber of $p$ is stable. We are looking for distributions on
$(V\oplus V)_{=1}$ which under the action of $SL(2)\times Z$ transform by
\[ \langle L,f(\lambda (av_1+bv_2),\lambda (cv_1+dv_2))\rangle =|\mu |^{-n}\langle L,f\rangle \] 
where
\[ ad-bc=1,\quad {\mathrm  and}\quad z=\mu \,\,{{\mathrm {Id}}  }\] 
This is clearly a closed set of distributions, stable under multiplication by locally constant
functions which factor through $p$. 

Then in order to prove that the only such distribution is 0 it is enough to prove it on each
fiber. So let us choose a point in ${\mathbb  P}(V)$ and one of its representative $v$ in $V$ so that
$(v,0)$ is a point in the fiber. The group $SL(2)\times Z$, acting on $V\oplus V$ is the same thing as the
group $GL(2)^+$, the exponent $+$ meaning that the determinant is a square. Then
\[ \pmatrix{a&b\cr c&d\cr}\pmatrix{v\cr 0\cr}=\pmatrix{av\cr cv\cr}\] 
so that the group is transitive on the orbit.
The isotropy subgroup is the group of matrices
\[ \pmatrix{1&b\cr 0&d\cr},\quad d\in {\mathbb  F}^{*2}\] 
The map
\[ (x,y)\mapsto (xv,yv)\] 
is an homeomorphism of ${\mathbb  F}^2\setminus (0,0)$ onto the orbit. The distribution must be
invariant under $SL(2)$ hence proportionnal to $dxdy$. Under the action of $Z$ the measure $dxdy$
satisfies
\[ \int _{{\mathbb  F}^2}f(\mu x,\mu y)dxdy=|\mu |^{-2}\int _{{\mathbb  F}^2}f(x,y)dxdy\] 
We may assume that $n\geq 3$ and even that $n$ is odd so this is not compatible with the required semi-invariance, except if the distribution is 0
which means  that $H_0\bigl ({\mathcal  S}(V\oplus V)_{=1})\bigr )=(0)$. The sequence
\[ H_0\bigl
({\mathcal  S}((V\oplus V)_{=1})\bigr )\longrightarrow H_0\bigl ({\mathcal  S}((V\oplus V)_{\leq 1})\bigr
)\longrightarrow H_0\bigl ({\mathcal  S}((V\oplus V)_{=0})\bigr)\] 
is exact, the first and third terms are 0, hence the second term is also 0; there is no
distribution with the required invariance supported on the ``singular set'' $rank(v_1,v_2)<2$.

The next thing to do is to find the distributions on ${\mathcal  S}(V\oplus V)_{=2}$. We first find the
covariant space for the action of $SL(2)$. Consider the map  $(v_1,v_2)\mapsto v_1ãv_2$ from
$V\oplus V$ into $ä^2V$. it is immediate to see that $SL(2) $ is simply transitive on the fibers. 
Each fiber tis homeomorphic to $SL(2)$. Fix a basis $e_1,\dots ,e_n$ of $V$. Call $H_2$ the
isotropy subgroup in $G=GL(n)$ of $e_1ãe_2$. Then $G$ acts transitively on  the quotient space
$SL(2)\setminus (V\oplus V)_{=2}$ and this space is homeorphic to $G/H_2$. if, for $f\in {\mathcal  S(\bigl
(V\oplus V)_{=2}\bigr )}$ we put 
\[ \varphi (g)=\int _{SL(2)}f\biggl (s\pmatrix{gv_1\cr gv_2\cr}\biggr )ds\] 
then $f\mapsto \varphi $ is a map onto ${\mathcal  S(G/H_2)}$ and the transpose is a one to one map from
${\mathcal  S}'(G/H)$ onto the space of distributions invariants under $SL(2)$. This means that the
space of covariants for the group $SL(2)$ is precisely ${\mathcal  S}(G/H_2)$. The group $G$ acts by
\[ g\varphi (x)=|{{\mathrm {Det}}   }  g|^{-1}\varphi (g^{{-1}}x)\] 
It will be convenient to define
\[ \psi (x)=\varphi (x)|{{\mathrm {Det}}   }  x|\] 
The function $\psi $ is such that $\psi (xh_2)=|{{\mathrm {Det}}   }  h_2|\psi (x)$ and, on $\psi $, the group $G $ acts by left
translations. 

 The group $H_2$ is almost a maximal parabolic subgroup. In fact an element
$h_2\in H_2$ is written as
\[ h_2=\pmatrix{m_1&u\cr 0&m_2\cr}\] 
with $m_1\in SL(2), m_2\in GL(n-2)$ and $u$ an arbitrary matrix with 2 rows and $n-2$ columns.
If we simply impose $m_1\in GL(2)$ we obtain a maximal parabolic subgroup $P_2$.

For any multiplicative character $\chi $ of ${\mathbb  F}^*$ the representation $\pi _\chi $ of $G$ is
defined as follows. The space $E_\chi $ of $\pi _\chi $ is the space of functions $w$ from $G$ to ${\mathbb  C}$,
locally constant and such that
\[ w(gp_2)=w(g)|{{\mathrm {Det}}   }  m_2||{{\mathrm {Det}}   }  m_1|^{-(n-2)/2}\chi ^{{-1}}({{\mathrm {Det}}   }  m_1)\] 
The group acts by left translations. This representation belongs to the degenerate principal
serie associated to $P_2$ and to the unitary degenerate principal serie if the character $\chi $ is
unitary. The central character of $\pi _\chi $ is
\[ \lambda {{\mathrm {Id}}  }\mapsto \chi ^{-2}(t)\] 
if we want this character to be trivial we must have $\chi ^2=1$ which gives only a finite number
of possibilities 

Let $d_t$ be the diagonal matrix with  diagonal coefficients $(t,1,1,\dots, 1)$.
Then for $\psi $ as above we put
\[ \psi _\chi (g)=\int _{{\mathbb  F}^*}\psi \bigl (g d_t\bigr )\chi (t)|t|^{(n-2)/2}d^*t\] 
Then $\psi _\chi \in E_\chi $.

The distribution $L$ on $(V\oplus V)=2$ is given by a distribution, which we still call $L$ on $G/H_2$. We are left with the condition
\[ \langle L,\psi (zg)\rangle =\langle L,\psi \rangle \quad z\in Z\] 

It is easy to prove that if, for all unitary characters $\chi $ such that $\chi ^2=1$, we have $\psi _\chi =0$ then $\langle L,\psi \rangle =0$.

 It follows that 
\[ H_0\Bigl ({\mathcal  S}\bigl ((V\oplus V)_{=2}\bigr )\Bigr )=\bigoplus_{\chi ^2=1}E_\chi \] 

\begin{proposition}
  The representation $\pi _\chi $ are irreducible and two by two inequivalent
\end{proposition}

Rubenthaler \cite{10}
 Theorem 5-4 page 481. The non equivalence is Theorem 5-8 of the same reference.

We will study the projection map from ${\mathcal  S}(V\oplus V)_{=2}$  onto $E_\chi $. This map is
given by
\[ f_\chi (g)=|{{\mathrm {Det}}   }  (g)|\int _{GL(2)}f\Bigl (x\pmatrix{ge_1\cr ge_2\cr}\Bigr )\chi (x)|{{\mathrm {Det}}   }  (x)|^{n/2}dx\] 
More generally for any complex number $s$ put
 \[ f_{\chi ,s}(g)=|{{\mathrm {Det}}   }  (g)|\int _{GL(2)}f\Bigl (x\pmatrix{ge_1\cr ge_2\cr}\Bigr )\chi (x)|{{\mathrm {Det}}   } 
(x)|^{s+n/2}dx\] 
Put
\[ x=\pmatrix{x_{11}&x_{1,2}\cr x_{21}&x_{22}\cr}\] 
then
\begin{eqnarray*}
 f_{\chi ,s}(g)=|{{\mathrm {Det}}   }  (g)|\int _{{\mathbb  F}^4}f\bigl
(g(x_{11}e_1+x_{12}e_2),g(x_{21}e_1+x_{22}e_2)\bigr
)\chi (x_{11}x_{22}-x_{12}x_{21})\\
|x_{11}x_{22}-x_{12}x_{21}|^{s+n/2}\frac{dx_{11}dx_{12}dx_{21}dx_{22}
}{|x_{11}x_{22}-x_{12}x_{21}|^2}\\
\end{eqnarray*}
A priori we take $f$ with support in $x_{11}x_{22}-x_{12}x_{21}\ne 0$ so that the integral is
always convergent. If we take any $f\in {\mathcal  S}(V\oplus V)$ then the integral converges for $\Re e
(s)+n/2>1$. In particular, if $n\geq 3$ the integral will be convergent for $s=0$. For such $f$ the
mapping $f\mapsto f_{\chi }$ has the same invariance properties as before so that this mapping
may be considered as a mapping from the covariant space $H_0({\mathcal  S}(V\oplus V))$ onto $E_\chi $
extending the map from $H_0({\mathcal  S}(V\oplus V)_{=2})$ onto $E_\chi $. This shows that:
\begin{proposition}
 For the action of $SL(2,{\mathbb  F}\times Z$ the covariant space of ${\mathcal  S}(V\oplus V)$is
\[ H_0({\mathcal  S}(V\oplus V))\approx H_0({\mathcal  S}(V\oplus V)_{=2})\approx \bigoplus E_\chi \] 
\end{proposition}
At this stage our problem is  to study the behaviour under the involution $\sigma $ of $G-$invariant linear forms on ${\mathcal  S}({\mathcal  N)}\otimes \bigl (\oplus E_\chi \bigr )$.

On $V\oplus V$ the involution $\sigma $ is essentially the Fourier transform. We now compute the corresponding involution on the $E_\chi $

To define the involution $\sigma $ we choose the map $u: V\rightarrow V^*$ defined by $u(e_i)=e^*_i$.
If $f\in {\mathcal  S}(V\oplus V)$ then
\[ \sigma (f)(v_1,v_2)=\widehat f(u(v_2),-u(v_1))\] 
To simplify the notations we will, until the end of the section, identify $V$ and $V^*$ with ${\mathbb  F}^n$ so that $u$ becomes
the identity operator\dots
In $V\oplus V$ we call the coordinates $x_1,\dots ,x_n$ for the first copy and $y_1,\dots ,y_n$ for
the second one.
On ${\mathbb  F}^2\oplus {\mathbb  F}^2$ consider the quadratic form $x_1y_2-x_2y_1$ and the corresponding
bilinear form
\[ x_1´_2+Å_1y_2-x_2´_1-Å_2y_1\] 
For $\varphi \in {\mathcal  S (\mathbb  F}^4)$ let
\[ \widehat{\widehat \varphi }(Å_1,Å_2,´_1,´_2)=\int _{{\mathbb 
F}^4}\varphi (x_1,x_2,y_1,y_2)\tau (x_1´_2+Å_1y_2-x_2´_1-Å_2y_1)dx_1dx_2dy_1dy_2\] 
If $\widehat \varphi $ is the Fourier transform defined with the usual bilinear form
($x_1Å_1+x_2Å_2+y_1´_1+y_2´_2$) then
\[ \widehat{\widehat \varphi }(Å_1,Å_2,´_1,´_2)=\widehat \varphi (´_2,-´_1,-Å_2,Å_1)\] 
The additive character $\tau $ is fixed . We may assume that it
is of order 0.

For $\varphi \in {\mathcal  S}({\mathbb  F}^4)$ and $\chi $ unitary  put
\[ Z_\varphi (\chi ,s)=\int _{{\mathbb 
F}^4}\varphi (x_1,x_2,y_1,y_2)\chi (x_1y_2-x_2y_1)|x_1y_2-x_2y_1|^s\frac{dx_1dx_2dy_1dy_2}{|x_1y_2-x_2y_1|^2}\] 
In \cite{9}we proved that
this integral is convergent for $\Re e(s)>1$, extends to a meromorphic function of $s$ and satisfies a functionnal equation
\[ Z_\varphi (\chi ,s)=\chi (-1)\rho   (\chi ,s)\rho  (\chi ,s-1)Z_ {\widehat{\widehat \varphi }}(\chi ^{-1},2-s)\] 
The factor $\rho  (\chi ,s)$ is defined in\cite {9}
We take $\chi ^2=1$. Then Proposition 3-1 of \cite{9}  gives explicitly the singularities of $Z_\varphi $. In particular the
only  possible real poles are 0 and 1. We want to take $s=n/2$ and $s=2-n/2$ with $n$ odd so we will have  no problem
with the poles. Finally note that $Z_{\widehat \varphi }=Z_{\widehat{\widehat \varphi }}$.

\overfullrule =0pt

Let $f_g(v_1,v_2)=f(gv_1,gv_2)$. Then
 \begin{eqnarray*}
  f_{\chi ,s}(g)=|{{\mathrm {Det}}   }  (g)|\int _{{\mathbb  F}^4}f_g\bigl (
x_1,x_2,0,\dots ,0,y_1,y_2,0,\dots ,0\bigr
)\chi (x_1y_2-x_2y_1)\\
|x_1y_2-x_2y_1|^{s+n/2}\frac{dx_{1}dx{_2}dy_{1}dy_{2}
}{|x_1y_2-x_2y_1|^2}\\
\end{eqnarray*}
Using the functional equation
\begin{eqnarray*}
 f_{\chi ,s}(g)=|{{\mathrm {Det}}   }  (g)|^{-1}\chi (-1)\rho  (\chi ,s+n/2)\rho  (\chi ,s+n/2-1)\\
\int _{{\mathbb  F}^4}\int _{{\mathbb  F}^{2n-4}}\bigl (\widehat f\bigr )_{^tg^{-1}}(x_1,x_2,x_3,\dots ,x_n,y_1,y_2,y_3,\dots
,y_n)\\
\chi (x_1y_2-x_2y_1)
|x_1y_2-x_2y_1|^{2-s-n/2}\bigl (\prod _3^ndx_jdy_j\bigr )\frac{dx_{1}dx{_2}dy_{1}dy_{2}
}{|x_1y_2-x_2y_1|^2}\\
\end{eqnarray*}
In this formula the integral is convergent for $\Re e(s)<1-n/2$ and $f_{\chi ,s}$ is the analytic continuation.
Let $\overline N$ be the subgroup of matrices
\[ \overline n=\pmatrix{I_2&0\cr c&I_{n-2}\cr}\] 
We put
\[ c=\pmatrix{Å_3&´_3\cr Å_4&´_4\cr \dots&\dots\cr \dots&\dots\cr Å_n&´_n}\] 
In the above double integral we may integrate over the Zariski open set $x_1y_2-x_2y_1\ne 0$ and then make the change
of variables
\begin{eqnarray*}
x_j&=&x_1Å_j+x_2´_j\quad j=3,\dots ,n\\
 y_j&=&y_1Å_j+y_2´_j\quad j=3,\dots n\\
 \end{eqnarray*}
 
Note that
\[ \overline n(x_1e_1+x_2e_2)=x_1e_1+x_2e_2+\sum _3^nx_1Å_j+x_2´_j\] 
and a similar formula for $\overline n(y_1e_1+y_2e_2)$. Thus
\begin{eqnarray*} 
f_{\chi ,s}(g)=|{{\mathrm {Det}}   }  (g)|^{-1}\chi (-1)\rho  (\chi ,s+n/2)\rho  (\chi ,s+n/2-1)\\
\int _{{\mathbb  F}^4}\int _{\overline N}\bigl (\widehat f\bigr )_{^tg^{-1}\overline n}(x_1,x_2,0,\dots
,0,y_1,y_2,0,\dots ,0)\\
\chi (x_1y_2-x_2y_1)
|x_1y_2-x_2y_1|^{-s+n/2}\bigl (d\overline n\bigr )\frac{dx_{1}dx{_2}dy_{1}dy_{2}
}{|x_1y_2-x_2y_1|^2}\\
\end{eqnarray*}
This integral is convergent for $\Re e(s)<1-n/2$ and is then equal to
\[ f_{\chi ,s}(g)=\chi (-1)\rho  (\chi ,s+n/2)\rho  (\chi ,s+n/2-1)\int _{\overline N}\widehat f_{\chi ,-s}(^tg^{-1}\overline n)d\overline
n\]  We want to take $s=0$ and without surprise the integral is divergent. Also we may replace $\widehat f$ by $\sigma (f)$
(they are equal up to the action of the usual Weyl group element  of  $SL(2)$ and this action dies in the integral so
\[ f_{\chi ,s}(g)=\chi (-1)\rho  (\chi ,s+n/2)\rho  (\chi ,s+n/2-1)\int _{\overline N}\sigma ( f)_{\chi ,-s}(^tg^{-1}\overline n)d\overline n\] 

Finally replace $f$ by $\sigma (f)
$
\[ \sigma (f)_{\chi ,s}(g)=\chi (-1)\rho  (\chi ,s+n/2)\rho  (\chi ,s+n/2-1)\int _{\overline N} f_{\chi ,-s}(^tg^{-1}\overline n)d\overline n\] 
If $s=0$ and $\chi ^2=1$ we have obtained a linear involution $\sigma _\chi $ on $E_\chi $ such that $\sigma _\chi \pi _\chi (g)=\pi _\chi (^tg^{-1})\sigma _\chi $. The representation $\pi _\chi $ is equivalent to its contragredient and $\sigma _\chi $ is one of the two involutive intertwinning maps.

Our distribution $R$ may be viewed as a  $G-$invariant linear form on ${\mathcal  S}({\mathcal   N})\otimes \bigl (\oplus E_\chi \bigr )$. For each $\chi $ we restrict it to ${\mathcal  S}({\mathcal   N})\otimes E_\chi $. It is skew symmetric:
\[ \langle R,f(-^tX)\oplus \sigma _\chi (e)\rangle =-\langle R,f(X)\oplus e\rangle \] 
Note that because of the homogeneity with respect to $X$ we can replace $^tX$ by $-^tX$. We have to prove that such $R$ are 0.

This can be reformulated in various ways. First one can fix a nilpotent orbit in ${\mathfrak  g}$ and some element $X$ of this orbit. Define $w_ X$ as in section ??? and let $C$ be the centralizer of $X$ in $G$. Then it would be enough to prove that a linear form $L$ on $E_\chi $ which is invariant under $C$ and skew:
\[ L(\sigma _\chi (e))=-L(\pi _\chi (w_ X)e)\] 
is 0.

Now consider the outer automorphism of $G$ defined by $g\mapsto {}^tg^{-1}$. Let
$\overline G=  {\mathbb  Z}_2\ltimes G$ be the corresponding semi-direct product. Thus if ${\mathbb 
Z_2}=\{1,\kappa \}$ then $^tg^{-1}=\kappa g\kappa ^{-1}$ and the group law is given by $(1,g)(\kappa ,g')=(\kappa ,^tg^{-1}g'),\,\,
(\kappa ,g)(1,g')=(\kappa ,gg')$.  A representation $\overline \pi $ of $\overline G$ in some complex vector
space $\overline E$ is specified by its restriction to $G$ and by a linear involution $\overline \pi (\kappa )$
such that $\overline \pi (\kappa )\pi (g)=\pi (^tg^{-1})\overline \pi (\kappa )$.

We obtain an extension $\overline \pi _\chi $ of $\pi _\chi $ to $\overline G$ by defining $\overline \pi _\chi (\kappa )=\sigma _\chi $ and an extension $\overline\nu$ of the action $\nu $ of $G$ in ${\mathcal  S}({\mathcal  N})$ by letting $\kappa $ acts through $f(X)\mapsto f(-{}^tX)$. Then  our distribution $R$ gives a $G-$ intertwinning  map $A$ between ${\mathcal  S}({\mathcal   N})$ and the smooth dual $E_\chi ^{*}$ of $E_\chi $. The skew symmetry now means that
\[ A\overline \nu (\chi )=-\overline \pi ^*_\chi (\kappa )A\] 

\vfill\eject

\centerline{\normalfont \Large \bfseries  Part II}
\vskip 0.5cm
\centerline{\normalfont \Large \bfseries  Orthogonal and Unitary Groups}
\vskip 1cm{}

\section{\hglue -12pt  . Introduction}

Let ${\mathbb F}$ be a local field non archimedean, of caracteristic 0. Let ${\mathbb D}$ be a  either ${\mathbb F}$ or a quadratic extension of ${\mathbb F}$ . If $x\in {\mathbb D}$ then $\overline x$ is the conjugate of $x$ if ${\mathbb D}\ne {\mathbb F}$ and is equal to $x$ if ${\mathbb D}={\mathbb F}$.

Let $W$ be a vector space over ${\mathbb D}$ of finite dimension $n+1\geqslant 2$.  Let $\langle  .,.\rangle  $ be a non degenerate hermitian form on $W$. This form is bi-additive and
\[ \langle  dw,d'w'\rangle  =d\,\,\overline{d'}\langle  w,w'\rangle   ,\quad \langle  w',w\rangle  =\overline{\langle  w,w'\rangle  }\] 
Given a ${\mathbb D}-$linear map $u$ from $W$ into itself, its adjoint $u^*$ is defined by the usual formula
\[ \langle  u(w),w'\rangle  =\langle  w,u^*(w')\rangle  \]

Choose a vector $e$ in $W$ such that $\langle  e,e\rangle  \ne 0$; let $U=e{\mathbb D}$ and $V=U^\perp$, the orthogonal complement. Then $V$ has dimension $n$ and the restriction of the hermitian form to $V$ is non degenerate.

Let $M$ be the unitary group of $W$ that is to say the group of all ${\mathbb D}-$linear maps $m$ of $W$ into itself which preserve the hermitian form or equivalently such that $mm^*=1$. Let $G$ be the unitary group of $V$. With the p-adic topology both groups are of type lctd ( locally compact, totally discontinuous and countable at infinity). They are reductive groups of classical type.

The group $G$ is naturally imbedded into $M$. Our goal is to show that the following conjecture follows from the conjectures of Part I 

\vskip 0.5cm

{\parindent 0pt
{\bf Conjecture $\mathbf 1'$} 
 {\sl If $\pi $ (resp $\rho $) is an irreducible admissible representation of $M$ (resp of $G$) then
\[ \mathrm{dim}\left  ({\mathrm {Hom}}_G(\pi _{| M},\rho )\right  )\leq 1\]
}}

Choose a basis $e_1,\dots e_n$ of $V$ such that $\langle  e_i,e_j\rangle  \in {\mathbb F}$. For 
\[ w=x_0e+\sum _1^nx_ie_i\]
put
\[ \overline w=\overline x_0e+\sum _1^n\overline {x_i}\, e_i\]
If $u$ is a ${\mathbb D}-$linear map from $W$ into itself, let $\overline u$ be defined by
\[ \overline u(w)=\overline{u(\overline w)}\]

Let $\sigma $ be the anti-involution $\sigma (m)=\overline m^{-1}$ of $M$; we will show that conjecture 1' is a consequence of

\vskip 0.5cm

{\parindent 0pt
{\bf Conjecture $\mathbf 2'$} 
{\sl A distribution on $M$ which is invariant under the adjoint action of $G$ is invariant under $\sigma $.
}}
\vskip 0.5cm

The involution $\sigma $ depends on the choice of the basis of $V$. However, changing the choice of the basis will just replace $\sigma $ by $g\sigma g^{-1}$ for some $g\in G$ so that the action on the space of invariant distributions does not depend on this choice.

The proof follows exactly the same path as in Part I. There are two mains differences. On one hand the Levi components of some of the parabolic subgroups of $G$ may not be of the orthogonal or unitary type; some components   of type $GL$ appear and this is why we need the ( conjectural) results of Part I. On the other hand the singular set is much simpler than in Part I. It has a natural stratification stable by the involution and such that on each strata an inductive argument works. This does not seem to be the case for the general linear group.

\section{ \hglue -16pt . Conjecture 2' implies Conjecture 1'}

In Chapter 4 of \cite{6} the following result is proved. Choose $\delta \in GL_{\mathbb F}(W)$ such that $\langle  \delta w,\delta w'\rangle  =\langle  w',w\rangle  $. If $\pi $ is an irreducible admissible representation of $M$, let $\pi ^*$ be its smooth contragredient and define $\pi ^\delta $ by
\[ \pi ^\delta (x)=\pi (\delta x\delta ^{-1})\] 
Then $\pi ^\delta $ and $\pi ^*$ are equivalent. We choose $\delta =1$ in the orthogonal case ${\mathbb D}={\mathbb F}$. In the unitary case, fix an orthogonal basis of $W$, say $e_1,\dots ,e_{n+1}$, such that $e_2,\dots ,e_{n+1}$ is a basis of $V$; put $\langle  e_i,e_i\rangle  =a_i$. Then
\[ \langle  \sum x_ie_i,\sum y_je_j\rangle  =\sum a_ix_i\overline{y_i}\]
Define $\delta $ by 
\[ \delta \left  (\sum x_ie_i\right  )=\sum \overline{x_i}e_i\] 
Note that $\delta ^2=1$.

Let $E_\pi $ be the space of $\pi $. Then, up to equivalence $\pi ^*$, is the representation $m\mapsto \pi (\delta m\delta ^{-1})$. If $\rho $ is an admissible irreducible representation of $G$ in a vector space $E_\rho $ then an element $A$ of ${\mathrm {Hom}}\left  (\pi ^*_{|G},\rho \right  )$ is a linear map from $E_\pi $ into $E_\rho $ such that
\[ A\pi (\delta g\delta ^{-1})=\pi (g)A,\quad g\in G\]
In turn the contragredient $\rho ^*$ of $\rho $  is equivalent to the representation $g\mapsto \rho (\delta g\delta ^{-1})$ in $E_\rho $. Then an element $B$ of ${\mathrm {Hom}}\left  (\pi _{|G},\rho ^*\right  )$ is a linear map from $E_\pi $ into $E_\rho $ such that
\[ B\pi (g)=\rho (\delta g\delta ^{-1})B,\quad g\in G\] 
As $\delta ^2=1$ the conditions on $A$ and $B$ are the same:
\[ {\mathrm {Hom}}\left  (\pi ^*_{|G},\rho \right  )\approx {\mathrm {Hom}}\left  (\pi _{|G},\rho ^*\right  ) \] 
However, assuming Conjecture 2', by Corollary 2-1 we have
\[ \dim\biggl ({\mathrm {Hom}}\left  (\pi ^*_{|G},\rho \right  )\biggr )\times \dim\biggl ({\mathrm {Hom}}\left  (\pi _{|G},\rho ^*\right  )\biggr )\leq 1 
\] 
so that both dimensions are 0 or 1. replacing $\rho $ by $\rho ^*$ we get Conjecture 1'. From now on we forget about Conjecture 1'.

\section{\hglue -14pt . A partial linearization}

The group $G$ acts on itself by the adjoint action and on V by the natural action so it acts on $G\times V$ and we may consider the space of invariant distributions ${\mathcal S}(G\times V)'$ on $G\times V$. Put
\[ \sigma ((g,v))=(\overline g^{-1},-\overline v)\]
It is an involution. 

\vskip 0.5cm
{\parindent 0pt
{\bf Conjecture $\mathbf 3'$} 
{\sl  A distribution on $G\times V$ which is invariant under $G$ is invariant under $\sigma $.
}}
\vskip 0.5cm

We claim that Conjecture 3'
 implies Conjecture 2'. Assume Conjecture 3' and apply it to $M$ acting on $M\times W$. Recall that $W=U\oplus   V$ with $U$ a non isotropic subspace of dimension 1. Let $e$ be a non zero element of $U$. Let
 \[ Y=\{(m,w)\in M\times W\, |\, \langle  w,w\rangle  =\langle  e,e\rangle  \, \}\]
 It is a closed subset, invariant under $M$ and under $\sigma $. Any distribution on $Y$ which is invariant under $M$ is thus invariant under $\sigma $.
 
 Now let $X\subset Y$ be the set of all pairs $(m,e)$. By Witt's theorem $MX=Y$. If $(m,e)\in X$ and if for some $m'\in M$ one has $m'(m,e)\in X$ then $m'e=e$ which is equivalent to $m'\in G$ and then $m'X=X$. We are in position to use a Frobenius type descent as described in the Appendix. In our case both groups $M$ and $H$ are unimodular. Choose a Haar measure on $M/H$. Let $S$ be a distribution on $X$ which is invariant under $G$. Then define a distribution $\theta (S)$ on $Y$ by
 \[ \langle  \theta (S),f\rangle  =\int _{M/G}\langle  S,f(mx)\rangle  \, dm\]
This is a distribution on $Y$, invariant under $M$ hence invariant under $\sigma $. The map $S\mapsto \theta (S)$ is a bijection between the space of $G-$invariant distributions on $X$ and the space of $M-$invariant distributions on $Y$. Let us compute $\theta (\sigma (S))$.
\begin{eqnarray*} \langle  \theta (\sigma (S)),f\rangle  &=&\int _{M/G}\langle  \sigma (S),f(m(m',e))\rangle  dm\\
&=&\int _{M/G}\langle  S,f((m\sigma (m')m^{-1},me))\, dm\\
&=&\int _{M/G}\langle  S,\sigma (f)(\sigma (m^{-1})m'\sigma (m),-\overline m\, e)\, dm
\end{eqnarray*} 
Now the map $m\mapsto \sigma (m^{-1})$ fixes $G$ and defines an involutive automorphism of $M/G$. In the above integral we change $m$ into $\sigma (m^{-1})=\overline m$:
\[ \langle  \theta (\sigma (S),f\rangle  =\int _{M/G}\langle  S,\sigma (f)(mm'm^{-1},- m\, e)\, dm\]
The linear map $u:w\mapsto -w$ belongs to the center of $M$ hence, changing $m$ into $um$ we finally get
\begin{eqnarray*} \langle  \theta (\sigma (S)),f\rangle  &=&\int _{M/G}\langle  S,\sigma (f)(mm'm^{-1},me)\rangle  \, dm\\
&=&\langle  \theta (S),\sigma (f)\rangle  
\end{eqnarray*}  
As $\theta (S)$ is invariant under $\sigma $ this implies that $\theta (S)=\theta (\sigma (S))$ hence $S=\sigma (S)$. But $X\approx M$ so we do get Conjecture 2'.

Let ${\mathfrak {g} }$ be the Lie algebra of $G$.The group $G$ acts on ${\mathfrak {g} }\times V$. On ${\mathfrak {g} }\times V$ consider the involution $\sigma $ given by $\sigma (x,v)=(-\overline x,-\overline v)$. To prove Conjecture 3' we shall first show that, in turn, it is implied by

 \vskip 0.5cm
{\parindent 0pt
{\bf Conjecture $\mathbf 4'$
 {\sl Any distribution on ${\mathfrak {g} }\times V$ which is invariant under $G$ is invariant under $\sigma $.}}
 \vskip 0.5cm

This is done using Harish-Chandra descent method. This means that we have to consider the Levi factors of the parabolic subgroups of $G$. This involves unitary groups of lower ranks and, unfortunately, groups of type $GL(n)$ over finite extensions of the base field. This is why the results of this second part are subject to the validity of the conjectures of the first part.

We assume that Conjecture 4 of the first part has been proved and then proceed by induction on the dimension of $V$.
\section{\hglue -14pt . The case $n=1$}
\
 In the orthogonal case ${\mathbb D}={\mathbb F}$, the group $G$ is $\pm {\rm Id}$. It acts trivially on itself and on its Lie algebra ( which is reduced to $(0)$). On $V={\mathbb F}$, the involution is $-{\rm Id}$ so both Conjecture 3' and 4' are tautological.
 
 In the unitary case ${\mathbb D}$ is quadratic extension of ${\mathbb F}$. The group $G$ is the group of $g\in G$ such that $g\overline g=1$; it acts trivially on itself and on its Lie algebra 
and by multiplication on $W={\mathbb D}$. On $G$ and ${\mathfrak {g} }$ the involution is trivial and on ${\mathbb D}$ it is $d\mapsto \overline d$. The orbits of $G$ in $G\times V$ are stable by $\sigma $ and furthermore if we denote by $\gamma (g)$ the action of $G$ on $G\times V$ then $\gamma (g)\sigma =\sigma \gamma (\overline g)$. A classical result of Bernstein-Zelevinsky  ( see \cite{6}, chapter 4, page 91) asserts that Conjecture 3' is true. The same argument works for the Lie algebra, giving Conjecture 4'.

\section{\hglue -14pt . Harish-Chandra's descent}

Let us go back to the general situation. Let $a\in G$, semi-simple; we want to describe its centralizer in     $G$. View $a$ as a ${\mathbb D}-$linear endomorphism of $V$ and call $P$ its minimal polynomial. Then, as $a$ is semi-simple, $P$ decomposes into irreducible factors $P=P_1\dots P_r$ two by two relatively prime. Let $V_i=\mathrm{Ker}    P_i(a)$ so that $V=\oplus   V_i$. Any element $x$ which commutes with $a$ will satisfy $xV_i\subset V_i$ for each $i$.
For
\[ R(T)=d_0+\cdots +d_mT^m,\quad d_0d_m\ne 0\] 
let
\[ R^*(T)=\overline{d_0}T^m+\cdots +\overline{d_m}\]
Then, from $aa^*=1$ we obtain, if $m$ is the degree of $P$
\[ \langle  P(a)v,v'\rangle  =\langle  v,a^{-m}P^*(a)v'\rangle  \] 
( note that the constant term of $P$ can not be 0 because $a$ is invertible).
It follows that $P^*(a)=0$ so that $P^*$ is proportional to $P$. Now $P^*=P_1^*\dots P_r^*$ hence there exists a bijection $\tau $ from $\{1,2,\dots ,r\}$ onto itself such that $P^*_i$ is proportional to $P_{\tau (i)}$.
Let $m_i$ be the degree of $P_i$. Then, for some non zero constant $c$
\[ 0=\langle  P_i(a)v_i,v_j\rangle  =\langle  v_i,a^{-m_i}P_i^*(a)v_j\rangle  =c\langle  v_i,a^{-m_i}P_{\tau (i)}(a)v_j\rangle  
,\quad v_i\in V_i,\,\, v_j\in V_j\] 
We have two possibilities.

{\bf Case 1:}$\,\, \tau (i)=i$. The space $V_i$ is orthogonal to $V_j$ for $j\ne i$; the restriction of the hermitian form to $V_i$ is non degenerate. Let ${\mathbb D}_i={\mathbb D}[T]/(P_i)$ and consider $V_i$ as a vector space over ${\mathbb D}_i$ through the action $(R(T),v)\mapsto R(a)v$. As $a_{|V_i}$ is invertible, $T$ is invertible modulo $(P_i)$; choose $U$ such that $T U=1$ modulo $(P_i)$. Let $\tau _i$ be the semi-linear involution of ${\mathbb D}_i$, as an algebra over ${\mathbb D}$:
\[ \sum d_jT^j\mapsto \sum \overline{d_j}U^j\quad{\rm modulo}\,\, (P_i)\] 
If $P_i$ is of degree one, that is to say if the restriction of $a$ to $V_i$ is a scalar operator then ${\mathbb D}_i\approx{\mathbb D}$. Any linear map $x_i$ from $V_i$ into itself commutes with $a_i=a_{|V_i}$. In this case we define $G_i$ to be the unitary group for the restriction to $V_i$ of our original hermitian form. We also put ${\mathbb F}_i={\mathbb F}$.

If $P_i$ has degree at least 2, let ${\mathbb F}_i$ be the subfield of fixed points for $\tau _i$. It is a finite extension of ${\mathbb F}$ and ${\mathbb D}_i$ is a quadratic extension of ${\mathbb F}_i$. Let $L$ be a non zero ${\mathbb D}-$linear form on ${\mathbb D}_i$.Then any ${\mathbb D}-$linear form $\ell$ on ${\mathbb D}_i$ may be written as $d\mapsto L(\lambda d)$ for some unique $\lambda \in {\mathbb D}_i$. We claim that we can choose $\lambda \ne 0$ in such a way that $\ell (\tau _i(d))=\tau _i(\ell (d))$ for all $d\in {\mathbb D }_i$. Indeed $d\mapsto \tau _i\left  (L(\tau _i(d))\right  )$ is ${\mathbb D}-$linear form so there exists $\alpha  \in {\mathbb D}_i$ such that, for all $d\in {\mathbb D}_i$
\[ \tau _i\left  (L(\tau _i(d))\right  )=L(\alpha  d)\] 
One gets $L(d)=L(\alpha  \tau _i(\alpha  )d)$, hence $\alpha  \tau _i(\alpha  )=1$. By Hilbert's theorem 90, there exists $\mu \in {\mathbb D}_i$ such that $\alpha  =\mu /\tau _i(\mu )$. Then $\ell : d\mapsto L(\mu d)$ has the required property.

If $v,v'\in V_i$ then $d\mapsto \langle  d(a)v,v'\rangle  $ is ${\mathbb D}-$linear map on ${\mathbb D}_i$ hence there exists $S(v,v')\in {\mathbb D}_i$ such that
\[ \langle  d(a)v,v'\rangle  =\ell (dS(v,v'))\]
One checks that $S$ is a non degenerate hermitian form on $V_i$ as a vector space over ${\mathbb D}_i$, a quadratic extension of ${\mathbb F}_i$.

Also a ${\mathbb D}-$linear map $x_i$ from $V_i$ into itself commutes with $a_i$ if and only if it is ${\mathbb D}_i$-linear and it is unitary with respect to our original hermitian form if and only if it is unitary with respect to $S$. So in this case we call $G_i$ is the unitary group of $S$. It does not depend upon the choice of $\ell$. As no confusion may arise, for $\lambda \in {{\mathbb D}_i}$ we define $\overline \lambda =\tau _i(\lambda )$.

We have an involution $\sigma $ on $G\times V$ well defined as soon we have chosen a basis of $V$over ${\mathbb D}$ such that the hermitian form$\langle  .,.\rangle  $ has, relative to this basis, a matrix with coefficients in ${\mathbb F}$. To build such a basis we are going to pick a suitable basis of $V_i$
 for each $i$. 
 
 We stay with Case 1. Choose a basis $(\xi _k)$ of $V_i$ over ${\mathbb D}_i$ such that $S(\xi _k,\xi _{k'})\in {\mathbb F}_i$.Then, for $v=\sum \lambda _k\xi _k$ we define $\overline v=\sum \overline{\lambda _k}\xi _k$. If $f\in \mathrm{End}  _{{\mathbb D}_i}(V_i)$ we define $\overline f$ by $\overline f(v)=\overline{f(\overline v)}$. Then, on $(G_i,V_i)$ we have the involution 
 \[ \sigma _i:\,\,\, (g_i,v_i)\mapsto (\overline{g_i}^{-1},-\overline {v_i})\] 
 and ${\mathfrak {g} }_i$ being the Lie algebra of $G_i$, on ${\mathfrak {g} }_i\times V_i$
 \[ \sigma _i:\,\,\, (X_i,v_i)\mapsto (-\overline{X_i},-\overline{v_i})\] 
 Now let $(u_r)$ be a basis of ${\mathbb F}_i$ over ${\mathbb F}$; it is also a basis of ${\mathbb D}_i$ over ${\mathbb D}$. Then $(u_r,\xi _k)$ is a basis of $V_i$ over ${\mathbb D}$ and, owing to our choice of $\ell$
 \begin{eqnarray*} \langle  u_r\xi _k,u_{r'}\xi _{k'}\rangle  &=&\ell\left  (S(u_r\xi _k,u_{r'}\xi _{k'}\right  )=\ell\left  (u_ru_{r'}S(\xi _k,\xi _{k'}))\right  )\\
 &=&\ell\left  (u_ru_{r'}S(\xi _{k'},\xi _{k}))\right  )=\langle  u_{r'}\xi _{k'},u_{r}\xi _{k}\rangle  \end{eqnarray*}
Thus to define the global involution $\sigma $, for the $V_i$ part, we may choose this basis. Note also, for $v=\sum y_{r,k}u_r\xi _k\in V_i$ with $y_{r,k}\in {\mathbb D}$ we have $\overline v=\sum \overline{y_{r,k}}u_r\xi _k$. If we consider $G_i\times V_i$ as imbedded, in the obvious way, inside $G\times V$ then on this subspace the involution $\sigma $ coincides with the involution $\sigma _i$.  
 
{\bf Case 2.} Suppose now that $j=\tau (i)\ne i$. Then $V_i\oplus   V_j$ is orthogonal to $V_k$ for $k\ne i,j$ and the restriction of the hermitian form to $V_i\oplus   V_j$ is non degenerate, both $V_i$ and $V_j$ being totally isotropic subspaces. Choose an inverse $U$ of $T$ modulo $P_j$. Then for any $P\in {\mathbb D}[T]$
\[ \langle  P(a)v_i,v_j\rangle  =\langle  v_i,\overline P(U(a))v_j\rangle  ,\quad v_i\in V_i,\,\, v_j\in V_j\]
where $\overline P$ is the polynomial deduced from $P$ by changing its coefficients into their conjugate.This defines a map, which we call $\tau _i$ from ${\mathbb D}_i$ onto ${\mathbb D}_j$. In a similar way we have a map $\tau _j$ which is the inverse of $\tau _i$. Then, for $\lambda \in {\mathbb D}_i$ we have $\langle  \lambda v_i,v_j\rangle  =\langle  v_i,\tau _i(\lambda )v_j\rangle  $.

 View $V_i$ as a vector space over ${\mathbb D}_i$. The action
\[ (\lambda ,v_j)\mapsto 
 \tau _i(\lambda )v_j\] 
defines a structure of ${\mathbb D}_i$ vector space on $V_j$. However note that for $\lambda \in {\mathbb D}$ we have $\tau _i(\lambda )=\overline \lambda $ so that $\tau _i(\lambda )v_j$ may be different from $\lambda v_j$. To avoid confusion we shall write, for $\lambda \in {\mathbb D}_i$
\[ \lambda v_i=\lambda *v_i\quad {\rm and}\quad \tau _i(\lambda )v_j=\lambda *v_j\] 

As in the first case choose a non zero ${\mathbb D}-$linear form $\ell$ on ${\mathbb D}_i$. For $v_i\in V_i$ and $v_j\in V_j$ the map $\lambda \mapsto \langle  \lambda *v_i,v_j\rangle  $ is a ${\mathbb D}-$linear form on ${\mathbb D}_i$, hence there exists a unique element $S(v_i,v_j)\in {\mathbb D}_i$ such that, for all $\lambda $
\[ \langle  \lambda *v_i,v_j\rangle  =\ell(\lambda S(v_i,v_j))\] 
The form $S$ is ${\mathbb D}_i-$ bilinear and non degenerate so that we can view $V_j$ as the dual space over ${\mathbb D}_i$ of the ${\mathbb D}_i$ vector space $V_i$.

Let $(x_i,x_j)\in \mathrm{End}  _{\mathbb D}(V_i)\times \mathrm{End}  _{\mathbb D}(V_j)$. They commute with $(a_i,a_j)$ if and only if they are ${\mathbb D}_i$-linear. The original hermitian form will be preserved, if and only if $S(x_iv_i,x_jv_j)=S(v_i,v_j)$ for all $v_i,v_j$. This means that $x_j$ is the inverse of the transpose of $x_i$. In this situation we define $G_i$ as the linear group of the ${\mathbb D}_i-$vector space $V_i$.

We now have to take care of the involutions. Let $(\xi _k)$ be a basis of $V_i$ over ${\mathbb D}_i$; let $(\eta _k)$ be the dual basis of $V_j$. Let $u$ be the ${\mathbb D}_i-$linear map from $V_i$ onto $V_j$ such that $u(\xi _k)=\eta _k$. Then on $G_i\times V_i\times V_j$ the "local" involution $\sigma _i$ is
\[ \sigma _i(g_i,v_i,v_j)=(u^{-1}\,\,{}^tg_i
\,\,u,u^{-1}(v_j),u(v_i))\] 
Next let $(\alpha  _r)$ be a basis of ${\mathbb D}_i$ over ${\mathbb D}$ and define the dual basis $(\beta _s))$ of ${\mathbb D}_i$ by $\ell (\alpha  _r\beta _s\rangle  =\delta _{r,s}$. Then $(\alpha  _r*\xi _k)$ is a basis of $V_i$ over ${\mathbb D}$ and $(\beta _s*\eta _t\rangle  )$ is a basis of $V_j$ over ${\mathbb D}$. We have $\langle  \alpha  _r*\xi _k,\beta _s*\eta _t\rangle  =\delta _{r,s}\delta _k,t$.   Using these basis define $\overline {v_i}$ for $v_i\in V_i$ and $\overline {v_j}$ for $v_j\in V_j$. Then if $x_i$ is a ${\mathbb D}-$linear map from $V_i$ into itself we define $\overline {x_i}$ by $\overline{x_i}(v_i)=\overline{x_i(\overline{v_i})}$ and similarly for $V_j$. Note that 
\[ u(\overline{v_i})=\overline{u(v_i)},\quad u^{-1}(\overline{v_j})=\overline{u^{-1}(v_j)}\] 
Then we define 
 \[ \gamma _i:\,\, V_i\times V_j\rightarrow V_i\times V_j\] 
 by
 \[ \gamma _i(v_i,v_j)=(-u^{-1}(\overline{v_j}),-u(\overline{v_i}))\] 
This mapping is ${\mathbb D}-$linear. Indeed if $\lambda \in {\mathbb D}$ then $\lambda v_i=\lambda *v_i$ whereas $\lambda v_j=\overline{\lambda }*v_j$ so that
 \begin{eqnarray*} \gamma _i\left  (\lambda v_i,\lambda v_j\right  )&=&\left  (-u^{-1}(\lambda *\overline{v_j}),-u(\overline \lambda *\,\,\overline{v_i})\right  )\\
 &=&\left  (-\lambda *u^{-1}(\overline{v_j}),-\overline \lambda *u(\overline{v_i})\right  )\cr
 &=&\left  (-\lambda u^{-1}(\overline{v_j}),-\lambda u(\overline{v_i})\right  )\\
 &=&\lambda \gamma _i(v_i,v_j)
 \end{eqnarray*}
Furthermore if $v_i,w_i\in V_i$ and $v_j,w_j\in V_j$ then
\begin{eqnarray*} \langle  \gamma _i(v_i,v_j),\gamma _i(w_i,w_j)\rangle  &=&\langle  u^{-1}(\overline{v_j}),u(\overline{w_i})\rangle  +\langle  u(\overline{v_i}),u^{-1}(\overline{w_j})\rangle  \\
&=&\langle  u^{-1}(\overline{v_j}),u(\overline{w_i})\rangle  +\overline{\langle  u^{-1}(\overline{w_j}),u(\overline{v_i})}\rangle  \\
&=&\ell\left  (S(u^{-1}(\overline{v_j}),u(\overline{w_i}))\right  )+\overline{\ell\left  (S(u^{-1}(\overline{w_j}),u(\overline{v_i}))\right  )}
\end{eqnarray*}
    
Over ${\mathbb D}_i$ the duality between $V_i$ and $V_j$ is defined by $S$ and $u$ is self adjoint, hence
\begin{eqnarray*} \langle  \gamma _i(v_i,v_j),\gamma _i(w_i,w_j)\rangle  &=&\ell\left  (S(\overline{w_i},\overline{v_j}\right  )+\overline{\ell\left  (S(\overline{v_i},\overline{w_j})\right  )}\\
&=&\langle  \overline{w_i},\overline{v_j}\rangle  +\overline{\langle  \overline{v_i},\overline{w_j}\rangle  }\\
&=&\langle  v_j,w_i\rangle  +\langle  v_i,w_j\rangle  \\
&=&\langle  (v_i,v_j),(w_i,w_j)\rangle  
\end{eqnarray*} 
Thus if we extend $\gamma _i$ to $V$ by letting it act by the identity on the $V_k$ for $k\ne i,j$ then $\gamma _i\in G$.

Let $(x_i,v_i,v_j)\in G_i\times V_i\times V_j$. On $V_i\times V_j$ $x_i$ acts as $(x_i,\,\, ^tx_i^{-1})$. The transposition is taken relative to the duality over ${\mathbb D}_i$. Let $(w_i,w_j)\in V_i\times V_j$.
\[ \gamma _i\circ (x_i,\,\, ^tx_i^{-1})\circ \gamma _i^{-1}(w_i,w_j)=\left  (u^{-1}\,\,{}^t\overline x^{-1}u(w_i),u\,\, \overline{x_i}u^{-1}(w_j)\right  )\] 
Thus if we imbed $G_i$ in $G$ by letting it act trivially on $V_k$ for $k\ne i,j$  we get
\[ \gamma _ix_i\gamma _i^{-1}=u^{-1}{}^t(\overline{x_i})^{-1}u\] 
Now $\sigma _i(x_i)=u^{-1}\,\, {}^tx_iu$ and $\sigma (x_i)=\overline{x_i}^{-1}$, therefore
\[ ad (\gamma _i)x_i=\sigma _i\sigma (x_i)\] 
The same is true for the action on $V_i\times V_j$. In other words, for the component $V_i\times V_j$ the involution $\sigma $ and $\sigma _i$ differ by the  action of an element $\gamma _i$ of $G$.
Recall that in case 1 they coincide. 

Let $M$ be the centralizer of $a$. Then $(M,V)$ decomposes as a "product", each "factor" being either of type $(G_i,V_i)$ with $G_i$ a unitary group (case 1) or $(G_I,V_i\times V_j)$ with $G_i$ a general linear group (case 2). Gluing together the involutions $\sigma _i$ we obtain an involution $\sigma _a$ on $(M,V)$. We define  $\gamma \in G$ as the product of the $\gamma _i$ ( in case 1 we take $\gamma _i$ trivial). Note that $\gamma ^2=1$ It is easy to check that $\sigma _a(a)=a$. For $m\in M$ we have $\sigma _a(m)=\gamma (\overline m)^{-1}\gamma =(\mathrm {ad}      \gamma )(\sigma (m))$ and for $v\in V$, $\sigma _a(v)=\gamma (\sigma (v))$.

We can now apply Harish-Chandra descent, repeating the arguments of part 1. By Lemma 5-2 and Lemma 5-3
of this part, there exists an open neighborhood $\omega $ of 0 in the Lie algebra ${\mathfrak {m}   }$ of $M$ which is closed, invariant under the adjoint action of $M$ and such that the exponential map is everywhere defined and submersive on $\omega $. Put $\Omega =\mathrm {Exp}   (\omega )$; it is open and invariant under the adjoint action of $M$. Then we may assume that the map $(g,m,v)\mapsto (gmag^{-1},gv)$ of $G\times \Omega \times V$ into $G\times V$ is everywhere submersive. 
The image of this map is a product $U\times V$ with $U$ an open neighborhood of $a$ in $G$ invariant under the adjoint action of $G$.

By Harish-Chandra submersion principle there exists a map $f\mapsto F_f$ from the space ${{\mathcal S}}(G\times \Omega \times V)$ onto ${\mathcal S}(U\times V)$ such that for any $\varphi \in {\mathcal S}(U\times V)$
\[ \int _{G\times \Omega \times V}f(g,m,v)\varphi (gmag^{-1},gv)dg\, dm\, dv=\int _{U\times V}F_f(u,v)\varphi (u,v)du\, dv\]
(the various Haar measures are fixed once for all).
\overfullrule=0pt

Define
\[ H_f(m,v)=\int _Gf(g,m,v)dg\]
Then, by Proposition 5-2 of Part 1 there is a well defined one to one map $\theta $ from the space of $G-$invariant distribution $T$ on $U\times V$ into the space of $M-$invariant distributions on $\Omega \times V$ such that, for all $f$ 
\[ \langle  T,F_f\rangle  =\langle  \theta (T),H_f\rangle  \]
Put
\[ f_1(g,m,v)=f(g,\sigma _a(m),\sigma _a(v))\]
Then
\[ H_{f_1}(m,v)=H_f(\sigma _a(v),\sigma _a(v))\]
Let us compute $F_{f_1}$.
\[ \int _{G\times \Omega \times V}f(g,\sigma _a(m),\sigma _a(v))\varphi (gmag^{-1},gv)dg\, dm\, dv=\int _{U\times V}F_{f_1}(u,v)\varphi (u,v)du\, dv\]
Change $m$ into $\sigma _a(m)$ and $v$ into $\sigma _a(v)$
\[ \int _{G\times \Omega \times V}f(g,m,v)\varphi (g\sigma _a(m)ag^{-1},g\sigma _a(v))dg\, dm\, dv=\int _{U\times V}F_{f_1}(u,v)\varphi (u,v)du\, dv\] 
Now
\[ g\sigma _a(m)ag^{-1}=g\sigma _a(am)g^{-1}=g\sigma _a(ma)g^{-1}=g\gamma \sigma (ma)\gamma ^{-1}g^{-1}=\sigma \left  (\sigma (g\gamma )^{-1}ma\sigma (g\gamma )\right  )\] 
and
\[ g\sigma _a(v)=g\gamma \sigma (v)=-g\gamma \overline v=\sigma \left  (\overline{g\gamma }v\right  )=\sigma \left  (\sigma (g\gamma )^{-1}v\right  )\] 
Therefore, changing $g$ into $\overline g \gamma ^{-1}$
\begin{eqnarray*} \int _{U\times V}F_{f_1}(u,v)&\varphi (u,v)dudv\\
&=&\int _{G\times \Omega \times V}f(g,m,v)\varphi \Bigl (\sigma \left  (\sigma (g\gamma )^{-1}ma\sigma (g\gamma )\right  ),\sigma \left  (\sigma (g\gamma )^{-1}v\right  )\Bigr )dgdmdv\\
&=&\int _{G\times \Omega \times V}f(\overline g \gamma ^{-1},m,v)\varphi ^\sigma (gmag^{-1},gv)dgdmdv\\
&=&\int _{U\times V}F_{f_2}(u,v)\varphi ^\sigma (u,v)du\, dv
\end{eqnarray*}
where
\[ \varphi ^\sigma (u,v)=\varphi (\sigma (u),\sigma (v))\] 
and
\[ f_2(g,m,v)=f(\overline g\gamma ^{-1},m,v)\] 
Note that $H_{f_2}=H_f$. We just found that $F_{f_1}=F^\sigma _{f_2}$. We also have $H_{f_1}=H^{\sigma _a}_f$. Therefore if $T$ is an invariant distribution on $U\times V$ then
\[ \langle  \theta (T)^{\sigma _a},H_f\rangle  =\langle  \theta (T),H_{f_1}\rangle  =\langle  T,F_{f_1}\rangle  =\langle  T,F_{f_2}^\sigma \rangle  =\langle  T^\sigma ,F_{f_2}\rangle  
=\langle  \theta (T^\sigma ),H_{f}\rangle  
\]
We conclude that $\theta (T^\sigma )=\theta (T)^{\sigma _a}$

On the other hand we also have a map from $\omega \times V$ onto $\Omega \times V$ everywhere submersive and which commute with the involution $\sigma _a$. Applying again the submersion principle we get a one to one map $\theta '$ from the space of $M-$invariant distributions on $\Omega \times V$ into the space of $M-$invariant distributions on $\omega \times V$ and this map $\theta '$ also commutes with $\sigma _a$.

Assume Conjecture 4'. As $\omega \times V$ is closed inside ${\mathfrak {g} }\times V$ any distribution on $\omega \times V$ extends to a distribution on ${\mathfrak {m}   }\times V$ with support in $\omega \times V$. Theorem D implies that such a distribution , if invariant, is invariant under $\sigma $. Using $\theta '$ we deduce that any invariant distribution on $\Omega \times V$ is invariant under $\sigma _a$. Using $\theta $ we now get that any $G-$invariant distribution on $U\times V$ is invariant under $\sigma $.

Start with an invariant distribution $T$ on $G\times V$ and suppose that $\sigma (T)=-T$. Let $g\in G$ and let $a$ be the semi-simple part of $g$. With notations as above we have that the restriction of $T$ to $U\times V$ must be 0. However $a$ belongs to the closure of the orbit of $g$ so that for some $x\in G$ we have $xgx^{-1}\in U$. Therefore $T$ is 0 in a neighborhood of $\{xgx^{-1}\}\times V$ but $T$is invariant hence is 0 in a neighborhood of $\{g\}\times V$ so $T$ must be 0. Any invariant distribution is symmetric with respect to $\sigma $.

This leaves us with the task of proving Conjecture 4'. The proof will be by induction on the dimension $n$ of $V$. The case $n=1$ has been checked. We shall assume that the theorem is true for $\dim V<n$ and also for the general linear case ( Part 1). We take $\dim V=n$.

The first step is to apply Harish-Chandra descent, this time on the Lie algebra. This is very similar to what we have done on the group.

We start with a non central semi-simple element $a\in {\mathfrak {g} }$. It is a ${\mathbb D}-$linear map from $V$ into itself such that $a+a^*=0$ where the adjoint is relative to the hermitian form. Let $P$ be the minimal polynomial of $a$ and $P=P_1\dots P_r$ its decomposition into irreducible factors. Put $V_i=\mathrm{Ker}    P_i(a)$ so that $V=\oplus   V_i$. If $M$ (resp. ${\mathfrak {m}   }$) is the cogmmutant of $a$ in $G$ (resp. in ${\mathfrak {g} }$) then an element $m\in M$ (resp. $X\in {\mathfrak {m}   }$) is diagonal with respect to the above decomposition of $V$.

Let ${\mathbb D}_i={\mathbb D}[T]/(P_i)$. Then $V_i$ has a structure of vector space over ${\mathbb D}_i$. Let $\tau $ be the semi-linear involution on ${\mathbb D}[T]$ such that $\tau (T)=-T$ and, for $\lambda \in {\mathbb D}$, $\tau (\lambda )=\overline \lambda $. Then $\tau (P)$ is proportional to $P$ so that for each $i$ the polymomial $\tau (P_i)$ is proportional to one of the $P_j$. This gives an involutive bijection, again called $\tau $, of $\{1,2,\dots ,r\}$ onto itself. As before we have two cases, Case 1: $\tau (i)=i$ and Case 2 $\tau (i)=j$ with $i\ne j$. With some trivial modifications we analyze each case as before.

In Case 1 we get an hermitian form on $V_i$ relative to some extension ${\mathbb F}_i$ of ${\mathbb F}$. The group $G_i$ is the corresponding unitary group. In Case 2 the group $G_i$ is the general linear group of $V_i$ as a vector space over ${\mathbb D}_i$ and $V_j$ is identified to the dual of $V_i$ over ${\mathbb D}_i$. The group $M$ is isomorphic to the product of the $G_i$. The situation $(M,V)$ is a product of $(G_i,V_i)$ ( case 1)
 and $(G_i,V_i,V_j)$ (Case 2).
 
 For a suitable choice of the basis, with respect to which the conjugations are defined, we obtain an involution $\sigma _a$ on $({\mathfrak {m}   }\times V)$, product of local involutions $\sigma _i$. The involution $\sigma $ on ${\mathfrak {g} }\times V$ is such that, for some $\gamma \in G$
 \[ \mathrm {Ad}  (\gamma )\sigma (X,v)=\sigma _a(X,v),\quad X\in {\mathfrak {m}   },\,\,\, v\in V\]

Because $a$ is not central,  Conjecture 4 or 4' is assumed for each component $G_i$. It follows that on ${\mathfrak {m}   }\times V$ any $M-$invariant distribution is symmetric with respect to $\sigma _a$.
Using Lemma 5-1 and 5-2 of Part 1 and proceeding exactly as in the group case we get an $M-$invariant open and closed neighborhood of $a$ in ${\mathfrak {m}   }$ and, if $U=\mathrm {Ad}  (G)\Omega $ a one to one linear map $\theta $ from the space ${\mathcal S}'(U\times V)^G$ of $G-$invariant distributions on $U\times V$ into the space ${\mathcal S}'(\Omega \times V)^M$ of $M-$invariant distributions on $\Omega \times V$. Furthermore $\theta (\sigma (T))=\sigma _a(\theta (T))$. By induction we know that $\sigma _a(\theta (T))=\theta (T)$ therefore $\sigma (T)=T$.

As in the group case we conclude that if $T$ is an invariant distribution on ${\mathfrak {g} }\times V$ such that $\sigma (T)=-T$ and if $X$ is an element of ${\mathfrak {g} }$ with a semi-simple part which is not central then $T=0$ in a neighborhood of $X$. 

Let ${\mathfrak z }$ be the center of ${\mathfrak {g} }$ and let ${\mathcal N}$ be the cone of nilpotent elements in $[{\mathfrak {g} },{\mathfrak {g} }]$. We have to prove that an invariant distribution $T$ on ${\mathfrak {g} }\times V$ with support contained in $({\mathfrak z }\times {\mathcal N})\times V$ is symmetric with respect to $\sigma $. The involution is trivial on ${\mathfrak z }$ and the group $G$ acts trivially on the center so it is enough to consider the case of a distribution with support contained in ${\mathcal N}\times V$. We will now do some reduction on the $V$ side.

\section{\hglue -12pt . Regular orbits}
We keep our notations. In particular $n=\dim V$. A pair $(X,v)\in {\mathfrak {g} }\times V$ is called {\bf regular} if $\{ v,Xv,\dots ,X^{n-1}v\}$ is a basis of $V$.

For any $(X,v)\in {\mathfrak {g} }\times V$ we have $\langle  X^rv,v\rangle  =(-1)^r\langle  v,X^rv\rangle  $. If we define
\[ q_i(X,v)=\langle  X^{i}v|v\rangle  \]
then
\[ \langle  X^iv,X^jv\rangle  =
(-1)^jq_{(i+j)}(X,v) \] 
When no confusion may arise we write $q_i$  for $q_i(X,v)$. Then define the
$n\times n$ matrix $A(X,v)$ ( or simply $A$) by
\[ A(X,v)=(a_{i,j})\quad {\rm with}\quad a_{i,j}=\langle  X^{i-1}v|X^{j-1}v\rangle  \] 
Let $D(X,v)$ be the determinant of $A(X,v)$.
\begin{proposition}
 $(X,v)$ is regular if and only if $D(X,v)\ne 0$. The set of regular
elements is a non empty Zariski open subset of ${\mathfrak {g} }\times V$.
\end{proposition}
If $(X,v)$ is regular then $A(X,v)$ is the matrix of $\langle  .,.\rangle  $ in the basis $v,Xv,\dots
X^{n-1}v $ and as the hermitian form is non degenerate, the matrix is non
singular. Conversely if there exists a non trivial linear relation among the
vectors $v,Xv,\dots X^{n-1}v$ then the same relation will be true for the
rows of $A$ which is thus singular.

The second assertion simply means that $D\ne 0$. Fix a basis $e_1,\dots e_n$ of $V$ with
coordinates $z_1,z_2,\dots z_n$, such that
\[ \langle  v,v\rangle  =\sum _1^n\lambda _jz_j\overline{z_j}\] 
Here the $\lambda _j$ are non zero elements of ${\mathbb F}$. Define $X$ by
\begin{eqnarray*} &Xe_1=e_2,\,\, Xe_2=-(\lambda _1/\lambda _2)\, e_1+e_3,\dots ,\\
&Xe_{n-1}=-(\lambda _{n-2}/\lambda _{n-1})\, e_{n-2}+e_n,\, Xe_n=-(\lambda _{n-1}/\lambda _n)\, e_{n-1}
\end{eqnarray*} 
Then $X\in {\mathfrak {g} }$ and $(X,e_1)$ is regular.

For $X\in {\mathfrak {g} }$ define the $D_j(X)$ by
\[ \mathrm {Det}  (T\,{\rm Id} -X)=T^n-\sum _0^{n-1}D_j(X)T^j\]
Like the $q_j$ the $D_j$ are invariant polynomial functions, with values in ${\mathbb D}$.
Note that the $q_j$ and $D_j$ are not algebraically independent.
\begin{proposition}
 Two regular elements are conjugate under $G$ if and only if they
give the same values to the invariants $q_j$ and $D_j$.
\end{proposition}
The necessity is clear. Conversely let $(X,v)$ and $(Y,w)$ be two regular
elements such that
\[ q_j(X,v)=q_j(Y,w),\quad D_j(X,v)=D_j(Y,w)\quad j=1,\dots n-1\] 
In particular 
\[ A(X,v)=A(Y,w)\]
so the linear map $g$ from $V$ to $V$ defined by $g(X^pv)=Y^pw$ for $p=0,\dots
n-1$ belongs to $G$. We claim that $gXg^{-1}=Y$. It is enough to check that
\[ g^{-1}YgX^pv=X^{p+1}v,\quad p=0,\dots n-1\] 
Now, if $p\leq n-2$
\[ g^{-1}YgX^pv=g^{-1}Y^{p+1}w=X^{p+1}v\] 
For $p=n-1$
\[ g^{-1}YgX^{n-1}v=g^{-1}Y^nw=g^{-1}\sum _0^{n-1}D_j(Y)Y^jw=\sum _{0}^{n-1}D_j(X)X^jv
=X^nv\] 
We define the {\bf regular orbits} as the orbits of regular elements. Each
such orbit is defined by the values of the $q_j$ and the $D_j$. These values
must be such that $D$, which is a polynomial in the $q_j$ should not take the
value 0 and also the  relations between the $q_j$ and
$D_j$ must be satisfied by the chosen values. Each such orbit is (Zariski)
closed (the invariants are constant  on the closure of any orbit !). If
$(X,v)$ is regular and if $g\in G$ is such that $g(X,v)=(X,v)$ which means  that
$gv=v$ and $gXg^{-1}=X$, then $gX^pv=(gXg^{-1})^pgv=X^pv$ for all $p$. By
definition of regular $\{v,Xv,\dots ,X^{n-1}v\}$ is a basis of $V$ so we
conclude that $g={\rm Id}$: the isotropy subgroup of a regular element is
trivial. In fact
\begin{theorem}
  An orbit  is regular if and only if it is closed and if the centralizer in $G$ of
an element of the orbit is trivial..
\end{theorem}
It is enough to prove that if the space $E$ generated by the $X^pv$ for
$p=0,\dots n-1$ is a proper subspace of $V$ then either $(X,v)$ has a non
trivial isotropy subgroup in $G$ or the orbit of $(X,v)$ is not closed
(or both \dots ).

Suppose first that the restriction of $\langle  .,.\rangle  $ to $E$ is non degenerate and let $F$ be
the orthogonal of $E$. Then $V=E\oplus   F$, an orthogonal decomposition. We know
that $X(E)\subset E$, hence $X(F)\subset F$. Relative to the above decomposition $X$ may be
writen as
\[ X=\pmatrix{X_1&0\cr 0&X_2
\cr}\]
Put
\[ Y=\pmatrix{0&0\cr 0&X_3
\cr}\] 
where $X_3$ belongs to the Lie algebra of the unitary group of the restriction of $\langle  .,.\rangle  $
to $F$ and  commutes with $X_2$. We may choose $X_3\ne 0$  Then
$Y\in {\mathfrak {g} }$ and $[Y,X]=0$ and also $Yv=0$. The orbit of $(X,v)$ is not of maximal
dimension, hence is not regular. 

In general the restriction of $\langle  .,.\rangle  $ to $E$ is degenerate. Let $N$ be its kernel and
choose a subspace $E_1$ of $E$ such that $E=N\oplus   E_1$. Note that $X(N)\subset N$ and
$X(E_1)\subset N\oplus   E_1$. Let $E_2$ be the orthogonal complement of $E_1$. Then
$V=E_1\oplus   E_2$ and $N\subset E_2$. Fix a decomposition
\[ E_2=E_2^0\oplus   E_2^+\oplus   E_2^-\] 
such that the restriction of $\langle  .,.\rangle  $ to $E_2^0$ is anisotropic, the orthogonal
complement of $E_2^0$ in $E_2$ is $E_2^+\oplus   E_2^-$, the subspaces $E_2^+$ and
$E_2^-$ are  maximal totally isotropic subspaces of $E_2$, the  form
$\langle  .,.\rangle  $ is a non degenerate  (hermitian) duality between $E_2^+$ and $E_2^-$ and finally
$N\subset E_2^+$. Let $M\subset E_2^+$ be a subspace such that $E_2^+=N\oplus   M$. Viewing
$E_2^-$ as the (anti)dual of $E_2^+$ we get a corresponding decomposition
$E_2^-=N^-\oplus   M^-$. Then
\[ V=E_1\oplus   (N\oplus   M)\oplus   (N^-\oplus   M^-)\oplus   E_2^0\]
For $t\in {\mathbb F}^*$ let $a_t$ be the endomorphism of $V$ defined by
\[ (a_t)_ {|N}=t{\rm Id}, \,\, (a_t)_{|N^-}=t^{-1}{\rm Id}\]
and $a_t={\rm Id}$ on the other components. One checks easily that $a_t\in G$.

If $v=v_0+v_1$ with $v_0\in N$ and $v_1\in E_1$ then $a_t(v)=tv_0+v_1$ so that
\[ \lim_{t\rightarrow 0}a_tv=v_1.\] 
Let $u\in V$ and suppose that $Xu$ has a non zero component  in $N^-$; then for
some $n\in N$ we have $\langle  Xu,n\rangle  \ne 0$ hence $\langle  u,Xn\rangle  \ne 0$.  However $Xn\in N$ hence this
implies that $u$ has a non zero component in $N^-$.

Let us compute $a_tXa_t^{{-1}}$. First if $n^-\in N^-$ then
$a_tXa_t^{-1}n^-=ta_tXn^-$ so that, when $t$ goes to 0, the limit is the $N^-$
component of $Xn^-$. Next if $n\in N$ then $Xn\in N$ and $a_tXa_t^{-1}n=Xn$. Finally
if $u\in V$ has a zero component relative to $N\oplus   N^-$ then $a_t^{-1}u=u$ and
we saw that $Xu$ has a zero component relative to $N^-$. It follows that the
limit of $a_tXa_t^{-1}u$ when $t$ goes to 0 is the projection of $Xu$ onto the
subspace
\[ E_1\oplus   M\oplus   M^-\oplus   E_2^0\] 
 We have checked that $a_tXa_t^{{-1}}$ has a limit
$X_1$. The point $(X_1,v_1)$ belongs to the closure of the orbit of $(X,v)$. 
Relatively to the decomposition $E=N\oplus   E_1$, the restriction of $X$ to $E$ may
be written in matrix form
\[ X_ {|E}=\pmatrix{Y&Y'\cr 0&Z
\cr}\] 
and the restriction of $X_1$ to $E$ is then
\[ (X_1)_{|E}=\pmatrix{Y&0\cr 0&Z
\cr}\] 
It follows that $X_1^pv_1$ is the projection onto $E_1$ of $X^pv$. Thus the
subspace generated by $X_1^pv_1,\,\, p=0,\dots n-1$ is $E_1$. If $N\ne (0)$,
(which we tacitly assumed) then $(X,v)$ and $(X_1,v_1)$ are not conjugate. Indeed if,
for some $g\in G$ we have $gXg^{{-1}}=X_1$ and $gv=v_1$ then, for all $p$ we get
$gX^pv=X_1^pv_1$ so that $g(E)=E_1$ which is impossible because $\dim (E)>\dim (E_1)$.
 
We have a more precise result.

Let $p$ be the dimension of $E$. Let $q$ be the largest integer such that the vectors
$v,Xv,\dots ,X^{q-1}v$ are linearly independent. Then $X^qv$ is a linear combination
\[ X^qv=\sum _0^{q-1}\alpha  _rX^rv\]
which implies that the subspace generated by the $X^rv,\,\, r=0,\dots ,q-1$ is stable
by $X$ hence must be equal to $E$. We have $p=q$.

If the restriction of $\langle  .,.\rangle  $ to $E$ is non degenerate, let $F$ be the orthogonal of $E$. Then
$F$ is stable by $X$. Relatively to the decomposition $V=E\oplus   F$
\[ X=\pmatrix{X_ E& 0\cr  0&X_ F\cr}\] 
Note that $X_ E$ (resp. $X_ F$ ) belongs to the Lie algebra ${\mathfrak {g} }_ E$ (resp.
${\mathfrak {g} }_ F$) of the unitary group $G_ E$ of $E$ (resp. $G_
F$ of $F$).
\begin{theorem}
For the orbit of $(X,v)$ to be closed it is necessary for the restriction of $\langle  .,.\rangle  $
to the subspace $E$ to be non degenerate. If this is true then the orbit is closed if and
only if $X_ F$ is a semi-simple element of ${\mathfrak {g} }_ F$.
\end{theorem}
The first assertion has just been proved; let us prove the second one. Suppose that the
orbit is closed and let $Y_ F\in {\mathfrak {g} }_ F$ be an element of the closure of the
orbit $\mathrm {Ad}  G_ F(X_ F)$. There exists a sequence $u_n$ of elements of $G_ F$
such that $g_nX_ Fg_n^{{-1}}$ converges to $Y_ F$. Define
\[ g_n=\pmatrix{1&0\cr 0&u_n\cr}\] 
Then $g_n\in G$ and $g_nv=v$ so that $g_n(X,v)$ converges to 
\[ \Biggl (\pmatrix{X_ E&0\cr 0&Y_ F\cr}\, ,\, v\Biggr )\] 
The orbit is closed so there exists $g\in G$ such that $gv=v$ and
\[ g\pmatrix{X_ E&0\cr 0&X_ F\cr}g^{-1}=\pmatrix{X_ E&0\cr 0& Y_ F\cr}\] 
For all $q$ we have
 \[ g(X^qv)=g(X_ E^qv)=\pmatrix{X_ E&0\cr 0& Y_ F\cr}^qv=X_ E^qv\] 
so that $g$ is the identity on $E$. Thus we may write
\[ g=\pmatrix{1&0\cr 0&u\cr},\quad u\in G_ F\] 
So $Y_ F=u(X_ F)$. The orbit of $X_ F$ is closed which is equivalent to the
semi-simplicity of $X_ F$.

Conversely assume that $X_ F$ is semi simple and let $(Y,w)=\lim g_n((X,v)$ be some
element of the closure of the orbit of $(X,v)$. We know that $\lim g_n(X^qv)=Y^qw$. In
particular the linear relations among the $X^qv$ remains valid for the $Y^qw$. The
subspace $E'$ generated by the $Y^qw$ is thus generated by $w,Yw,\dots ,Y^{p-1}w$. We
claim that $\dim E'$ is exactly $p$. Indeed we have $\langle  X^rv,X^sv\rangle  =\langle  Y^rw,Y^sw\rangle  $. Now
$v,Xv,\dots X^{p-1}v$ is a basis of $E$ and $\langle  .,.\rangle  $ restricted to $E$ is non degenerate so that the matrix
\[ \Bigl (\langle  X^{{i-1}}v,X^{j-1}v\rangle  \Bigr )_{0\leq i,j\leq p-1}\] 
is non singular. Hence the matrix
\[ \Bigl (\langle Y^{{i-1}}w,Y^{j-1}w\rangle  \Bigr )_{0\leq i,j\leq p-1}\] 
is also non singular which implies that the vectors $w,Yw,\dots ,Y^{p-1}w$ are
linearly independent and that the restriction of $\langle  .,.\rangle  $ to $E'$ is non degenerate. There
exists $u\in G$ such that $u(X^qv)=Y^qw$ for $q=0,1,\dots p-1$. To prove that $(Y,w)$
belongs to the orbit of $(X,v)$ it is enough to prove that $u^{-1}((Y,w))$ belongs to this
orbit. In other words we may assume that, for all $e\in E$ we have $\lim g_n(e)=e$; in
particular $w=v$.

Consider the decomposition $V=E\oplus   F$. The hermitian form $\langle  .,.\rangle  $ defines a semi-linear map $s$
from $V$ to its dual. Put
\[ s=\pmatrix{s_ E&0\cr
0&s_ F
\cr}\] 
For any linear map from $V$ to $V$ let $f^*=s^{-1}\, ({}^tf)s$ be the adjoint map; we use
a similar notation for $E$ and $F$. 

Write $g_n$ in matrix form
\[ g_n=\pmatrix{\alpha  _n&\beta _n\cr \gamma _n&\delta _n
\cr}\] 
and put
\[\beta _n^*=s_ F^{-1}(^t\beta _n)s_ E,\quad \gamma _n^*=s_ E^{-1}(^t\gamma _n)s_ F\] 
Then $g_n\in G$ is equivalent to the set of relations
 \begin{eqnarray*}\alpha  _n^*\alpha  _n+\gamma _n^*\gamma _n&=&1\\
\beta _n^*\beta _n+\delta _n^*\delta _n&=&1\\
\alpha  _n^*\beta _n+\gamma _n^*\delta _n&=&0
\end{eqnarray*}
The condition $\lim g_n(e)=e$ for $e\in E$ means that 
\[ \lim \gamma _n=0,\quad \lim \alpha  _n=1\]
In particular $\alpha  _n$ is invertible for $n$ large enough and from the third equation above
we get
\[\beta _n=-(\alpha  _n^*)^{-1}\gamma _n^*\delta _n\]
and
\[\beta _n^*\beta _n=\delta _n^*\gamma _n\alpha  _n^{-1}(\alpha  _n^{*})^{-1}\gamma _n^*\delta _n\] 
Thus
\[ \delta _n^*\gamma _n\alpha  _n^{-1}(\alpha  _n^{*})^{-1}\gamma _n^*\delta _n+\delta _n^*\delta _n=1\]
which we rewrite as
\[ (\delta _n^*)^{-1}(\delta _n)^{-1}=1+\varepsilon _n\] 
with
\[ \varepsilon _n=-\gamma _n\alpha  _n^{-1}(\alpha  _n^{*})^{-1}\gamma _n^*\] 
Note that $\varepsilon _n=\varepsilon _n^*$ and has limit 0. Consider the map $f\mapsto f^*f$ from the
space of self adjoint maps from $F$ to $F$ into itself. The differential at the origin is
$h\mapsto h^*+h$, a bijection. Hence we have around the identity a local
diffeomorphism. Forgetting about a finite number of values of $n$, we see that
there exists a sequence $f_n$ of  maps from $F$ to $F$, converging to 1
and such that $f^*_nf_n=1+\varepsilon _n$. Going back to $\delta _n$ we have
\[ (\delta _n^*)^{-1}(\delta _n)^{-1}=f_n^*f_n\]
or
\[ (f_n\delta _n)^*(f_n\delta _n)=1\] 
Then $u_n=f_n\delta _n\in G_ F$ and 
\[ \delta _n=f_n^{-1}u_n\]

 Now \[ g_nXg_n^{-1}=g_nXg_n^*\]  is explicited in
matrix form as \[ \pmatrix{\alpha  _n&\beta _n\cr \gamma _n&\delta _n
\cr}\pmatrix{X_ E&0\cr
0&X_ F
\\}\pmatrix{\alpha  _n^*&\gamma _n^*\cr
\beta _n^*&\delta _n^*
\\}=\pmatrix{\alpha  _nX_ E\alpha  _n^*+\beta _nX_ F\beta _n^*&\alpha  _nX_ E\gamma _n^*+\beta _nX_ F\delta _n^*\cr
\gamma _nX_ E\alpha  _n^*+\delta _nX_ F\beta _n^*&\gamma _nX_ E\gamma _n^*+\delta _nX_ F\delta _n^*
\cr}\]
The limit of this matrix is
\[ Y=\pmatrix{X_ E&0\cr
0&Y_ F
\cr}\] 
so, in particular
\[ Y_ F=\lim (\gamma _nX_ E\gamma _n^*+\delta _nX_ F\delta _n^*) \] 
But $\lim \gamma _n=0$ so, using the above formula for $\delta _n$ we get
\[ Y_ F=\lim (f_n^{-1}u_nX_ Fu_n^{-1}(f_n^*)^{{-1}})\] 
Define $\eta _n$ by
\[ Y_ F=f_n^{-1}u_nX_ Fu_n^{-1}(f_n^*)^{{-1}}+\eta _n\] 
so that $\lim \eta _n=0$ and rewrite the equality as
\[ f_n(Y_ F-\eta _n)f_n^*=u_nX_ Fu_n^{-1}\] 
The left side belongs to the orbit $G_ F X_ F$ so its limit $Y_ F$ belongs to the
closure of this orbit. However $X_ F$ is supposed to be semi simple so the orbit is
closed and we get that $Y_ F=uX_ F u^{-1}$ for some $u\in G_ F$. Finally put
\[ g=\pmatrix{1&0\cr
0&u
\cr}\] 
Then $g\in G$ and $gv=v$ and $gXg^{-1}=Y$.

Finally let us note that it follows from Theorem 17-1 that a regular orbit always carries a $G-$invariant measure. Also a regular orbit is stable by $\sigma $ (Proposition 17-2) so that this invariant measure is  stable by $\sigma $. it follows  that any invariant distribution on the regular set is symmetric with respect to $\sigma $.

\section{\hglue -12pt . Reduction to the singular set}
 Let $X\in {\mathfrak {g} }$ and $v\in V$. Consider the sequence of vectors $(X^jv)$. As $X^n$ is a linear combination
of the $X^j$ for $j<n$, for $s\geqslant n$, $X^sv$ is a linear combination of the $X^jv$ for $j<n$. Let $r$ be an integer,
$1\leq r\leq n$ and choose a class $\gamma $ of hermitian forms, non degenerate and of rank $r$.  We shall denote by
$\Sigma (\gamma )$ the set of all $(X,v)$ such that the subspace
\[ E_r(X,v)=\sum _0^{r-1}{\mathbb D}X^jv\] 
is of dimension $r$ and such that the restriction of $\langle  .,.\rangle  $ to this subspace belongs to the class $\gamma $. Note that ,
for $r=n$, if $\langle  .,.\rangle  $ does not belong to the class $\gamma $, then $\Sigma (\gamma )=\emptyset  $ and if $\langle  .,.\rangle  $ belongs to $\gamma $ then $\Sigma
(\gamma )$ is the subset of regular elements.
\begin{lemma}
 Each $\Sigma (\gamma )$ is a, possibly empty, open subset. The union of all $\Sigma (\gamma )$ is the
set of all $(X,v)$ such that $E_{n}(X,v)$ is not a totally isotropic subspace
\end{lemma}
Indeed,the linear independance of the $X^jv.\, j=0,\dots ,r-1$ and the non degeneracy of the restriction of $\langle  .,.\rangle  $
to $E_r(X,v)$ are equivalent to the non vanishing of the determinant of the matrix
\[ A_r(X,v)=\left  (\langle  X^{i-1}v,X^{j-1}v\rangle  \right  ),\quad i,j=0,\dots ,r-1\] 
Once $r$ is fixed, there is only a finite number of classes $\gamma $ characterized by the discrimant, the determinant
of the above matrix, modulo the squares (in the base field ${\mathbb F}$) and the Hasse symbol in the case of quadratic forms (${\mathbb D}={\mathbb F}$).

To prove the second assertion let $q_j=\langle  X^jv,v\rangle  $. If $q_j=0$ for $j=0,\dots ,n-1$ then all the $q_j$ are 0 and the
restriction of $\langle  .,.\rangle  $ to $E_{n}(X,v)$ is 0. Otherwise let $r$ be the smallest integer such that $q_{r-1}\ne 0$ The
matrix
\[ A_r(X,v)=\pmatrix{0&0&0&\dots&\dots&\dots&q_{r-1}\cr
0&0&0&\dots&\dots&-q_{r-1}&*\cr
\dots&\dots&\dots&\dots&q_{r-1}&*&\dots\cr
\dots&\dots&\dots&\dots&\dots&\dots&\dots\cr
0&-q_{r-1}&*&\dots&\dots&\dots&\dots\cr
q_{r-1}&*&\dots&\dots&\dots&\dots&*
\cr}\] 
is non singular. Hence $(X,v)$ belongs to some $\Sigma (\gamma )$.

Assume that there exists a subspace
$E$ of $V$, of dimension $r$, such that the restriction of $\langle  .,.\rangle  $ to $E$ belongs to a given  class $\gamma $. Let $\Xi (\gamma )$ be
the
 closed subset of $\Sigma (\gamma )$ defined by $E_r(X,v)=E$. By Witt's theorem $G\Xi (\gamma )=\Sigma(\gamma )$.

Suppose that $(X,v)\in \Xi (\gamma )$, that $g\in G$ and that $g(X,v)\in \Xi (\gamma )$. Put $w=gv$ and $Y=gXg^{-1}$. Then $\{ w,Yw,\dots
Y^{q-1}w\}$ is a basis of $E$. But $g(X^jv)=Y^jw$ so $g(E)=E$. This implies that $g(\Xi (\gamma ))=\Xi (\gamma )$. Let $H$ be the
stabilizer of $\Xi (\gamma )$ in $G$ or what amounts to the same the stabilizer of $E$ in $G$. For any subspace $F$ of $V$
call $G_ F$  the unitary group of the restriction of $\langle  .,.\rangle  $ to $F$ and ${\mathfrak {g} }_ F$ its Lie
algebra. Then $H$ is isomorphic to $G_ E\times G_{E^\perp}$ where $E^\perp$ is the orthogonal complement to $E$. We
are in a situation to apply the results of the Appendix . 

Use the decomposition $V=E^\perp \oplus   E$ to write $X$ in matrix form
\[ X=\pmatrix{x_{1,1}&x_{1,2}\cr x_{2.1}&x_{2,2}
\cr}\] 
Also the hermitian form $\langle  .,.\rangle  $ defines a semi linear map from $V$ to $V^*$ which we write as
\[ \pmatrix{s^\perp&0\cr 0&s
\cr}\] 
Then $X\in {\mathfrak {g} }$ is equivalent to
\begin{eqnarray*} {}^tx_{1,1}s^\perp+s^\perp x_{1,1}&=&0\cr
{}^tx_{2,2}s+sx_{2,2}&=&0\cr
sx_{2,1}+{}^tx_{1,2}s^\perp&=&0
\end{eqnarray*} 

If $(X,v)\in \Xi (\gamma )$, for $j\leq r-2$ we have $X^jv=x_{2,2}^jv$ and $x_{1,2}X^jv=0$ The linear map $x_{1,2} $ from $E$ to
$E^\perp$ is of rank at most  1. More precisely there is a unique $w\in E$ such that
\[ \langle  \sum _0^{r-1}\lambda _jX^jv,w\rangle  =\lambda _{r-1}\] 
and if we define $u\in E^\perp$ by $u=x_{1,2}X^{r-1}v$ then, for $e\in E$, we have $x_{1,2}(e)=\langle  e,w\rangle  u$.
Thus starting from $(X,v)\in \Xi (\gamma )$ we obtain, on one hand  $(x_{1,1},u)$ which belongs to ${\mathfrak {g} }_
{E^\perp}\times E^\perp$, with no particular condition and on the other hand $(x_{2,2},v)$ which is a regular
element of ${\mathfrak {g} }_ E\times E$.

Conversely if we start with two such elements $(x_{1,1},u)$ and a regular $(x_{2,2},v)$ then we define $w$
using $v$ and $x_{2,2}$ and then $x_{1,2}$ using $w$ and $u$ which allows us to recover $X$. Thus if we call
$({\mathfrak {g} }_ E\times E)_{\rm reg}$ the open set of regular elements in ${\mathfrak {g} }_ E\times E$ then we have an
homeomorphism
\[ ({\mathfrak {g} }_{E^\perp}\times E^\perp)\times ({\mathfrak {g} }_ E\times E)_{\rm reg}\rightarrow \Xi (\gamma )\] 
Write an element $h$ of $H$ as
\[ h=\pmatrix{h_ {E^\perp}&0\cr 0&h_ E
\cr}\] 
Then, for $(X,v)\in \Xi (\gamma )$ we get $hv=h_ Ev$ and
\[ hXh^{-1}=\pmatrix{h_{E^\perp}x_{1,1}h^{-1}_{E^\perp}&h_{E^\perp}x_{1,2}h_ E^{-1}\cr
h_ Ex_{2,1}h_{E^\perp}^{-1}&h_ Ex_{2,2}h_ E^{-1}
\cr}\] 
For $e\in E$
\[ h_{E^\perp}x_{1,2}h_ E^{-1}(e)=\langle  h_ E^{-1}(e))h_{E^\perp}(u),w\rangle  =\langle  e,h_
E(w)\rangle  h_{E^\perp}(u)\]
But $h_ E(w)$ satisfies
\[ \left  \langle  h_ Ex_{2,2}^jh_ E^{-1}(h_ Ev),h_ E(w)\right  \rangle  =\langle  x_{2,2}^jv)|w\rangle  \quad j=0,\dots ,q-1\] 
so that, under the action of $h$, the vector $w$ is replaced by $h_ E(w)$ and we get that
the action of $H$ on $\Xi (\gamma )$ is the product of the actions of the two smaller unitary groups.

We will now draw several consequences of the above remarks. First the $G-$orbits in $\Sigma (\gamma )$ are in one
to one correspondance with the $H-$orbits in $\Xi (\gamma )$ which themselves are obtained by taking the product of a
regular orbit of $G_ E$ in ${\mathfrak {g} }_ E\times E$ by an orbit of $G_{E^\perp}$ in ${\mathfrak y
g}_{E^\perp}\times E^\perp$. Let $(X,v)\in \Xi (\gamma )$.
\begin{proposition}
 $(X,v)$ is regular if and only if $(x_{1,1},u)$ is regular
\end{proposition}
Indeed by induction one checks that there exists constants $c_{q,j}$ such that
\[ X^{q+r}v-\Biggl (x_{1,1}^qu-\sum _0^{q-1}c_{q,j}x_{1,1}^ju\Biggr )\in E,\quad q=0,\dots ,n-1-r\] 
We know that $(X,v)$ is regular if and only if $X^jv,\, j=0,\dots ,n-1$ is a basis of $V$. In all cases $X^jv,\,
j=0,\dots ,r-1$ is a basis of $E$ thus $(X,v)$ is regular if and only if $x_{1,1}^qu,\, q=0,\dots ,n-1-r$ is a basis
of $E^\perp$ which means that $(x_{1,1},u)$ is regular.

\begin{proposition} 
 Let $(X,v)\in \Xi (\gamma )$. The orbit $G(X,v)$ is fixed by $\sigma $ if and only if the orbit $G_{E^\perp}(x_{1,1},u)$ is
fixed. If any $G_ {E^\perp}-$invariant distribution on $E^{\perp}$ is symmetric with respect to the involution $\sigma _ {E^\perp}$  then any invariant distribution on
$\Sigma (\gamma )$ is symmetric with respect to the involution $\sigma $.
\end{proposition}
To define $\sigma $, we may choose a basis $e_1,\dots ,e_n$ of $V$ over ${\mathbb D}$ relative to which the hermitian
form is diagonalized:
\[ \left  \langle  \sum z_ie_i,\sum z'_ie_i\right  \rangle  =\sum \alpha  _iz_i\overline{z'_i},\quad \alpha  _i=\overline{\alpha  _i}\] 
Using the coordinates relative to this basis we define $\overline v, \overline X,\dots$ and $\sigma (X,v)=(-\overline
X,-\overline v)$. A different choice of basis will simply replace $\sigma $ by $g\circ \sigma $ for some $g\in G$ so that the
above Proposition is independent of this choice. We choose a basis adapted to the decomposition
$V=E\oplus   E^\perp$.

 With the same notations as above, if we replace $(X.v)$ by $\sigma (X,v)=(-\overline X,-\overline v)$ then $u$ is
changed into $(-1)^{r+1}\overline u$ and $x_{1.1}$ into $-\overline{x_{1,1}}$,so that $(x_{1,1},u)$ is replaced by
$(-\overline {x_{1,1}},(-1)^{r+1}\overline u)=\sigma (x_{1,1},(-1)^ru)$. However in $E^\perp$ the map $-{\rm Id}$
belongs to the center of the unitary group so
 \[ \sigma (x_{1,1},(-1)^ru)=\sigma \circ (-{\rm Id})^r(x_{1,1},u)\] 
On the other hand $(x_{2,2},v)$ is replaced by $(-\overline x_{2,2},-\overline v)=\sigma (x_{2,2},v)$. As any regular
orbit is stable we get the first assertion.Also  for invariant distributions our map commutes with $\sigma $ which
implies the second assertion because on the regular set an y distributionn is symmetric.

Going back to our general situation, as we proceed by induction we may assume the theorem for $E^{\perp}$. Therefore it remains to consider the case  of invariant distributions with support contained into the set of all $(X,v)$ such that  $\langle  X^jv,v\rangle  =0$ for all $j$. On the other hand we also know that it is enough to consider distributions with support contianed into the set of all $(X,v)$ with $X$ nilpotent.

We shall say that $(X,v)$ is {\bf singular} if $X$ is nilpotent and if, for all $j$, $\langle  BX^jv,v\rangle  =0$. Let $\Sigma$ be the singular set. We have to prove that any invariant distribution with support contained into the singular set is stable by $\sigma $.

\section{\hglue -12pt . The singular set}
We keep our notations. An invariant distribution $T$ on $V$ is called symmetric (resp. skew symmetric) if $\sigma (T)=T$ (resp. $\sigma (T)=-T$).  
An element $(X,v)\in {\mathfrak {g} }\times V$ belongs to $\Sigma$ if and only if the subspace $E(X,v)$
generated by the vectors $X^jv$ is totally isotropic and if $X$ is nilpotent. Let
$\Sigma_r$
 be the subset of $\Sigma$ defined by $\dim E(X,v)=r$; then 
\[ \Sigma=\bigcup _0^{r_0}\Sigma_r\] 
where $r_0$ is the largest possible dimension for a totally isotropic subspace of
$V$.

The subset $\Sigma_r$ is defined by the following conditions: the invariants $D_j(X)$
and $q_j(X,v)$ are all 0, the vectors $v,Xv,\dots ,X^{r-1}v$ are linearly
independant and $X^rv=0$. Indeed as $E(X,v)$ is stable under $X$ which is nilpotent
the restriction of $X^r$ to $E(X,v)$ is 0. In particular $\Sigma_r$ is locally closed, hence
is an lctd space and we may introduce the spaces ${\mathcal S}(\Sigma_r)$ and ${\mathcal
S}'(\Sigma_r)$. Also $\Sigma_r$ is stable by $\sigma $. We will show in a moment that $\Sigma_r\ne \emptyset  $ for $r\leq r_0$.

We claim that in order to prove our theorem it is enough to prove that,
for each $r$, an element of ${\mathcal S}'(\Sigma_r)$ which is invariant under $G$ is
symmetric . Indeed, let $T$ be an invariant distribution
supported on $\Sigma$ and suppose that it is skew symmetric. The restriction of $T$ to
the open subspace $\Sigma_{r_0}$ of $\Sigma$ is invariant and  skew
symmetric. If we have proved that an invariant distribution on $\Sigma_{r_0}$ is
always symmetric we conclude that this restriction is 0. Next we restrict to
$\Sigma_{r_0-1}$ an open subspace of the complement of $\Sigma_{r_0}$ in $\Sigma$ and so on,
obtaining finally that $T=0$. As the involution commutes with the action of $G$,
any invariant distribution is the sum of an invariant symmetric distribution and a
skew symmetric one,  and this is enough to conclude.

From now on we fix $r\leq r_0$, $r\ne 0$ .We shall deal later with the case $r=0$. Two totally isotropic subspaces of dimension $r$ are
conjugate under $G$. Fix one of them $E$. Up to conjugation by $G$ we may assume
that $\{ v,Xv,\dots ,X^{r-1}v\}$ is a basis of $E$. Choose a second totally isotropic
subspace $E^*$, of dimension $r$ such that the sum $E+E^*$ is direct and that $\langle  .,.\rangle  $
is a non degenerate (semi) duality between $E$ and $E^*$. Let $F=(E\oplus   E^*)^\perp$. We
write $V$ as
\[ V=E\oplus   F\oplus   E^*\] 
Fix a basis $\{ e_0,\dots ,e_{r-1}\}$ of $E$ and let $\{ e_0^*,\dots ,e^*_{r-1}\}$ be the dual basis of $E^*$ so that $\langle  e_i,e^*_j\rangle  =\delta _{i,j}$. Choose a basis $(u_j)$ of $F$ such that $\langle  u_j,u_k\rangle  \in {\mathbb F}$. Then, relative to these choices, the hermitian form $\langle  .,.\rangle  $ has a matrix $s_ V$ which is equal to its adjoint ( transpose of the conjugate) and in fact has its coefficients in ${\mathbb F}$ .In matrix notations the hermitian form is then $v^*s_ Vv'=\langle  v',v\rangle  $, the unitary group is defined by $g^*s_ Vg=s_ V$ and the Lie algebra by $X^*s_ V+s_ VX=0$ Relatively to the decomposition $V=E\oplus   F\oplus   E^*$  in matrix form 
\[ s_ V=\pmatrix{0&0&{\rm Id}\cr
0&s&0\cr
{\rm Id}&0&0
\cr},\quad s=s^*\] 

If $h\in {\rm GL}(E)$ then
\[ g=\pmatrix{h&0&0\cr
0&{\rm Id}_ F&0\cr
0&0& (h^{*})^{-1}
\cr}\]
belongs to $G$. Up to conjugation by $G$
we may assume that $X^jv=e_{r-1-j}$ for $j=0,\dots ,r-1$. . Let $\Xi _r$ be the set of all $(X,v)$ satisfying
the above conditions. Note that it is a closed subspace of $\Sigma_r$ and we know that
$G\Xi _r=\Sigma _r$.

If $(X,v)\in \Xi _r$ then $v=e_{r-1}$ and, making explicit the condition $X\in {\mathfrak {g} }$ we
see that
\[ X=\pmatrix{n_r&-b^*s&a\cr
0&Y&b\cr
0&0&-n_r^*
\cr}\]
where $n_r$ is defined by $n_r(e_j)=e_{j-1}$
 for $j=1,\dots ,r-1$ and $n_r(e_0)=0$. Also $a+a^*=0$ 
 and $b$ is arbitrary. Finally $Y\in {\mathfrak y
g}_ F:\,\, Y^*s+sY=0$ and $X$ is nilpotent if and only if $Y$ is nilpotent. 

Let $(X,v)\in \Xi _r$ and let $g\in G$ be such that $g(X,v)\in \Xi _r$. Put $g(X,v)=(X',v')$ so that
$v=v'=e_{r-1}$ and $X'=gXg^{-1}$. Then $g(e_{r-1})=e_{r-1}$ and
\[ g(e_{r-1-j})=(gX^j)(e_{r-1})=(X^{'j}g)(e_{r-1})=X^{'j}(e_{r-1})=e_{r-1-j}\]
Thus $g$ is the identity on $E$. Conversely if $g\in G$ is the identity on $E$ then
$g\Xi _r=\Xi _r$. Let $H$ be the subgroup of $g\in G$ which are the identity on $E$.

If $h\in H$, then, expliciting the condition $h\in G$ we get
\[ h=\pmatrix{1&-x^*sy&z\cr
0&y&x\cr
0&0& 1
\cr}\] 
where $x,y,z$ are such that
\begin{eqnarray*} y^*sy&=&s\cr
z+ z^*+x^*sx&=&0
\end{eqnarray*}
The reductive part of $H$ is isomorphic to $G_ F$, the unitary group of the
restriction of the form to $F$. 
We can now use the results of the Appendix. There is a one to one
correspondence $\theta :S\mapsto T$ from the space of $H-$invariant distributions on
$\Xi _r$ onto the space of $G-$invariant distributions on $\Sigma_r$ such that
\[ \langle  T,f\rangle  =\int _{G/H}\langle  S_\xi ,f(g\xi )\rangle  dg\] 
Now consider the involution $(X,v)\mapsto (-\overline X,-\overline v)$ ( we use the above basis of $V$). The strata $\Sigma_r$ is stable but
not $\Xi _r$. To correct the situation, define $\tau _ E : E\rightarrow E$ by $\tau _
E(e_{r-1-j})=(-1)^{j+1}e_{r-1-j}$ and $\tau \in G$ by
\[ \tau =\pmatrix{\tau _ E&0&0\cr
0&-{\rm Id}_ F&0\cr
0&0& (\tau  _ E^*)^{-1}
\cr}\] 
If $(X,v)\in \Xi _r$ then $\tau (-\overline X,-\overline v)$ belongs to $\Xi _r$.  Indeed  $v=e_{r-1}=\overline v=-\tau (v)$ and, using the
above notations for $X$ we get
\[ \tau (-\overline X)\tau ^{-1}=\pmatrix{n_r&-\tau _ E\, \overline{b^*}s&-\tau _ E\overline a\,
\tau ^*_ E\cr
0&-\overline Y&\overline b\,\tau ^*_ E\cr
0&0&-n_r^*
\cr}\] 
If $S\in {\mathcal S}'(\Xi _r)^H$
\begin{eqnarray*} \langle  \theta (S),f(-\overline X,-\overline v)\rangle  &=&\int _{G/H}\langle  S,f(g(-\overline X,-\overline v))\rangle  dg\\
&=&\int _{G/H}\langle  S,f\left  (g\tau ^{-1}(\tau (-\overline X)\tau ^{-1},v)\right  )\rangle  dg\\
&=&\int _{G/H}\langle  S,f\left  (g(\tau (-\overline X)\tau ^{-1},v)\right  )\rangle  dg
\end{eqnarray*}
For the last equality note that $\tau $ normalizes $H$, thus operates on $G/H$ by right
multiplication and being an involution it leaves the measure on $G/H$ invariant.
The distribution $T=\theta (S)$ is symmetric if and only if the distribution $S$ is
symmetric relative to the involution $(X,v)\mapsto (-\tau \overline X\tau ^{-1},v)$ of $\Xi _r$.

Our next step is to use the invariance under the unipotent radical $U$ of $H$. First
note that $v$ does not play any role so we may as well consider $\Xi _r$ as a
subset of ${\mathfrak {g} }$.If
\[ X=\pmatrix{n_r&- b^*s&a\cr
0&Y&b\cr
0&0&- n_r^*
\cr},\quad{\rm and}\quad u=\pmatrix{1&- x^*s&z\cr
0&1&x\cr
0&0& 1
\cr}\] 
then
\[ uXu^{-1}=\pmatrix{n_r&-(b')^*s&a'\cr
0&Y&b'\cr
0&0&- n_r^*
\cr}\]
with

\begin{eqnarray*}b'&=&b-(Yx+xn_r^*)\\
a'&=&a-(n_rz+z\, n_r^*)-n_r\, x^*sx+x^*sYx+
 b^*sx-x^*sb
\end{eqnarray*}
 Also in the sequel the fact that $Y$ is nilpotent will play no role whatsoever. So
we will just forget about this assumption and take for $Y$ any element of ${\mathfrak y
g}_ F$.
We study the map $L:x\mapsto Yx+x n_r^*$ from ${\mathrm {Hom}}(E^*,F)$
into itself. Let $c_j=x(e^*_j)$ and $\gamma _j= L(x)(e^*_j)$. Then
\[ \gamma _j=Yc_j+c_{j+1},\quad j=0,\dots r-2,\quad \gamma _{r-1}=Yc_{r-1}\] 
In particular $L(x)=0$ if and only if $c_j=(-1)^jY^jc_0$ with $Y^rc_0=0$. Hence the
kernel of $L$ is isomorphic to the kernel of $Y^r$. Next we look for the $x$ such
that $\gamma _0=\gamma _1=\dots =\gamma _{r-2}=0$. The conditions are $c_j=(-1)^jY^jc_0$ for
$j=1,\dots r-1$ and then $\gamma _{r-1}=(-1)^{r-1}Y^rc_0$. In other words the
intersection of the image of $L$ with the subspace $\gamma _0=\dots =\gamma _{r-2}=0$ is the
subspace $\gamma _0=\dots =\gamma _{r-2}=0, \gamma _{r-1}\in {\rm Im}Y_r$. Its dimension is equal to
the rank of $Y^r$, hence if $F_0\subset F$ is any subspace supplementary to this image
then the subspace $\tilde F_0:\gamma _0=\dots =\gamma _{r-2}=0,\gamma _{r-1}\in F_0$ has the same
dimension as the kernel of $L$ and intersects trivially the image of $L$ hence is a
supplementary subspace of the image of $L$. It follows that if $X\in \Xi _r$ then up to
conjugation we can assume that $b\in \tilde F_0$. At least for the moment we will
just remember that, up to conjugation we may assume that $b(e^*_j)=0$ for $j\leq r-2$
(in matrix terms, only the last column of $b$ is non zero). let $\tilde F$ be the
space of all such $b$.

Let us take $x=0$; only $a$ is affected and now we have to study the map
$N:z\mapsto n_rz+z\,  n^*_r$ with $z+z^*=0$ . We use matrix notations with $z=(z_{i,j})$ and $N(z)=(\zeta _{i,j})$. Recall
that $0\leq i,j\leq r-1$; we make the convention that $z_{r,j}=z_{i,r}=0$. Then
\[ \zeta _{i,j}=z_{i+1,j}+z_{i,j+1}\] 
Here we deal with the structure of vector space over the base field ${\mathbb F}$.
The kernel of $N$ is the subspace of anti-hemitian matrices $z$ such that
\[ z_{i,j}=\cases{(-1)^jz_{i+j,0}&if $i+j\leq r$ and $i+j$ odd\cr
0& if $i+j\geqslant r$ or $i+j$ even
\cr}\]
Its dimension is $r/2$ if $r$ is even and $(r-1)/2$ if $r$ is odd. 

A straightforward computation shows   that the elements of the image of $N$
such that $\zeta _{i,j}=0$ for $j\leq r-2$ and $i\leq r-2$ (only the last column and the last row
are non zero) are the anti-hermitian matrices of size $r$ such that the last row is
\[ (\dots\dots \alpha  _3\,\, 0\,\,\alpha  _2\,\, 0\,\, \alpha  _1\,\,0\,\,0)\] 
with all the $\alpha  _i\in {\mathbb F}$.
Then denote by $E_0$ the subspace of $E$ of all vectors 
\[ \dots +\dots +\beta _3e_{r-6}+\beta _2e_{r-4}+\beta _1e_{r-2}\]
with again $\beta _i\in {\mathbb F}$
and let $\tilde E_0$ be the space of all antihermitian matrices such that only the
last column and last row are non zero and such that the last column belongs to
$E_0$. Then the dimension of $\tilde E_0$ is equal to the dimension of the kernel
of $N$, the intersection of $\tilde E_0$ with the image of $N$ is reduced to $(0)$,
so that $\tilde E_0$ is a supplementary subspace of the image of $N$.

Let $\Theta _r\subset\Xi _r$ be the set of all $X\in \Xi _r$ such that $Y\in {\mathfrak {g} }_
F\,\,\,$,$a\in \tilde E_0$ and $b\in \tilde F$. we have just checked that $U\Theta _r=\Xi _r$. We
cannot use the Appendix  but Harish-Chandra 's
submersion principle does apply. Indeed both $\Xi _r$ and $\Theta _r$ are vector spaces,
hence analytic manifolds ans we have a surjective map
\[ \pi : U\times\Theta _r\rightarrow \Xi _r\] 
given by $\pi (u,\xi )=u\xi u^{-1}$. To compute the differential of $\pi $ we evaluate
\[ \pi (u\mathrm {Exp}  (tA),\xi +\eta ),\quad A\in {\mathfrak  u}={\rm Lie}(U),\,\,\, \eta \in\Theta _r\]
The differential is
\[ (A,\eta )\mapsto \mathrm {Ad}  (u)(\eta +[A,\xi ])\] 
If
\[ A=\pmatrix{0&-\beta ^*s&\alpha  \cr
0&0&\beta \cr
0&0&0&
\cr}\quad \xi =\pmatrix{n_r&-b^*s&a\cr
0&Y&b\cr
0&0&-n_r^*
\cr}\] 
then
\[ [A,\xi ]=\pmatrix{0&*&-(\alpha  \, n_r^*+n_r\alpha  )+(b^*s\beta -
\beta ^*sb)\cr
0&0&-(Y\beta +\beta ^*n_r)\cr
0&0&0
\cr}\]
By definition of $\Theta _r$ this formula shows that $\Xi _r=\Theta _r+\mathrm {Ad} (\xi ){\mathfrak  u}$ and the
map is everywhere submersive.

Thus there exists a surjective map $f\mapsto F_f$ of ${\mathcal S}(U\times\Theta _r)$ onto
${\mathcal S }(\Xi _r)$ such that, for any $\varphi \in {\mathcal S}(\Xi _r)$
\[ \int _{U\times\Theta _r}f(u,\xi )\varphi \left  (\pi (u,\xi )\right  )dud\xi =\int _{\Xi _r}F_f(X)\varphi (X)dx\] 
Here $du, d\xi $ and $dX$ are Haar measures fixed, a priori, once for all. By
transposition we obtain a one to one map $\theta :S\mapsto T$ from ${\mathcal S}'(\Xi _r)$
into ${\mathcal S}'(U\times\Theta _r)$ such that $\langle  T,f\rangle  =\langle  S,F_f\rangle  $.

Let $g\in U$; we have $d(gu)=du$ and $d(\mathrm {Ad}  g(X))=dX$. Put $f^g(u,\xi )=f(g^{-1}u,\xi )$.
\begin{eqnarray*} \int _{\Xi _r}F_{f^g}(X)\varphi (X)dX&=&\int _{U\times\Theta _r}f(g^{-1}u,\xi )\varphi \left  (\pi (u,\xi )\right 
)dud\xi \\
&=&\int _{U\times\Theta _r}f(u,\xi )\varphi \left  (\mathrm {Ad}  g \,\,\pi (u,\xi )\right  )dud\xi \\
&=&\int _{\Xi _r}F_f(X)\varphi (\mathrm {Ad}  g\,\, X)dX\\
&=&\int _{\Xi _r}(F_f)^g(X)\varphi (X)dX
\end{eqnarray*} 
This implies that $(F_f)^g=F_{f^g}$ and that $S$ is $U-$invariant if and only if
$T=\theta (S)$ is $U-$invariant. If this is the case then we may decompose $T$ as
$T=duæR$ where $R$ is a distribution on $\Theta _r$.

Now consider the action of the reductive part of $H$ that is to say of $G_
F$. If $y\in G_ F$, imbedded in $H$, it acts on $\Xi _r$ simply by changing
$b$ into $yb$ and on $Y$ by the adjoint action. In particular $\Theta _r$ is fixed. Put
$^yf(u,\xi )=f(u,\mathrm {Ad}  y (\xi ))$ and $f_y(u,\xi )=f(y^{-1}uy,\xi )$.
\begin{eqnarray*} \int _{\Xi _r}F_ {^y\hskip -2 pt f}(X)\varphi (X)dX&=&\int _{U\times\Theta _r}f(u,\mathrm {Ad} 
y(\xi ))\varphi \left  (\pi (u,\xi )\right  )dud\xi \\
 &=&\int _{U\times\Theta _r}f(u,\xi )\varphi \left  (\pi (u,\mathrm {Ad}  y^{-1}(\xi ))\right  )dud\xi \\
&=&\int _{U\times\Theta _r}f(u,\xi )\varphi \left  (\mathrm {Ad}  (uy^{-1})\xi \right  )dud\xi \\
&=&\int _{U\times\Theta _r}f_y(u,\xi )\varphi \left  (\mathrm {Ad}  (y^{-1})\pi (u,\xi )\right  )dud\xi \\
&=&\int _{\Xi _r}F_{f_y}(X)\varphi (\mathrm {Ad}  y^{-1}(X))dX\\
&=&\int _{\Xi _r}F_ {f_y}\left  (\mathrm {Ad}  y(X)\right  )\varphi (X)dX
\end{eqnarray*} 
This implies that
\[ F_{^y\hskip -2pt f}(X)=F_{f_y}\left  (\mathrm {Ad}  y(X)\right  )\]
Suppose that the distribution $S$ is invariant under $H$. Then it is invariant under
$U$ and $y$. As before put $T=\theta (S)=duæR$. Then, by the biinvariance of $du$ we have
\[ \langle  T,f\rangle  =\langle  T, f_y\rangle  =\langle  S,F_{f_y}\rangle  \] 
By the invariance of $S$ under $y$
\[ \langle  S,F_{f_y}(X)\rangle  =\langle  S,F_{f_y}\left  (\mathrm {Ad}  y(X)\right  )\rangle  =\langle  S,F_{^y\hskip
-2pt f}\rangle  =\langle  T,^y\hskip -2pt f\rangle  \]
The conclusion is that $T$ is invariant under $y$ and $S\mapsto R$ is a one to one
map from the space of $H-$invariant distributions on $\Xi _r$ into the space of
$G_ F-$invariant distributions on $\Theta _r$.

Let us look at the involution $X\mapsto \tau (-\overline X)$ on $\Xi _r$. Take $X\in\Theta _r$. Then $a\in {\mathbb F}$ and  $\tau _
E\,a\, \tau ^{*}_ E=-a$ and $- b\tau _ E^*=b$. As $\tau $ leaves the $Y$
component invariant we get that, on $\Theta _r$ the involution is given by $Y\mapsto
-\overline Y,\, b\mapsto -\overline b,\, a\mapsto a$.

Put $^\tau \hskip -2pt f(u,\xi )=f(u,\mathrm {Ad}  \tau (-\overline \xi ))$ and $f_\tau (u,\xi )=f(\tau ^{-1}\overline u\tau ,\xi )$. This makes
sense because $\mathrm {ad}      \tau $ leaves $U$ invariant and note that $du$ is also invariant by
$\mathrm {ad}      \tau $.
 \begin{eqnarray*} \int _{\Xi _r}F_ {^\tau \hskip -2 pt f}(X)\varphi (X)dX&=&\int _{U\times\Theta _r}f(u,\mathrm {Ad} 
\tau (-\overline \xi ))\varphi \left  (\pi (u,\xi )\right  )dud\xi \\
 &=&\int _{U\times\Theta _r}f(u,\xi )\varphi \left  (\pi (u,\mathrm {Ad}  \tau ^{-1}(-\overline \xi ))\right  )dud\xi \\
&=&\int _{U\times\Theta _r}f(u,\xi )\varphi \left  (-\mathrm {Ad}  (u\tau ^{-1})\overline \xi \right  )dud\xi \\
&=&\int _{U\times\Theta _r}f_\tau (u,\xi )\varphi \left  (\mathrm {Ad}  (\tau ^{-1})(-\overline{\pi (u,\xi )})\right  )dud\xi \\
&=&\int _{\Xi _r}F_{f_\tau }(X)\varphi (\mathrm {Ad}  \tau ^{-1}(-\overline X))dX\\
&=&\int _{\Xi _r}F_ {f_\tau }\left  (\mathrm {Ad}  \tau (-\overline X)\right  )\varphi (X)dX
\end{eqnarray*} 
This implies that
\[ F_{^\tau \hskip -2pt f}(X)=F_ {f_\tau }\left  (\mathrm {Ad}  \tau (-\overline X)\right  )\] 
or, replacing $f$ by $f_\tau $
\[ F_{^\tau \hskip -2pt (f_\tau )}(X)=F_f\left  (\mathrm {Ad}  \tau (-\overline X)\right  )\] 
If $S$ is $H-$invariant and $\theta (S)=T=du\otimes R$ then
\[ \langle  S,F_f\left  (\mathrm {Ad}  \tau (-\overline X)\right  )\rangle  =\langle  S,F_{^\tau \hskip -2pt (f_\tau )}\rangle  =\langle  T,^\tau \hskip
-2pt (f_\tau )\rangle  \] 
Assume, for a moment that $T$ is symmetric : $\langle  T,^\tau \hskip -2pt f\rangle  =\langle  T,f\rangle  $, then
\[ \langle  T,^\tau \hskip -2pt (f_\tau )\rangle  =\langle  T,f_\tau \rangle  =\langle  du\otimes R,f(\tau ^{-1}u\tau ,\xi )\rangle  =\langle  T,f\rangle  =\langle  S,F_f\rangle  \] 
and $S$ is symmetric.

Therefore we must prove that $T=du\otimes R$ is symmetric. First suppose that $r>0$ 
 so that $F$ is of dimension strictly positive and strictly smaller than the
dimension of $V$. identify $\Theta _r$ with $({\mathfrak {g} }_ F\oplus   F)\oplus   E_0$. Let $\alpha  \in {\mathcal
S}(U)$ and $\beta \in {\mathcal S}(E_0)$. On ${\mathfrak {g} }_ F\oplus   F$ the distribution
\[ \psi (Y,b)\mapsto \langle  T,\alpha  æ\beta æ\psi \rangle  =\int _U\alpha  (u)du\, \langle  R,\psi (Y,b)\beta (a)\rangle  \] 
is invariant hence symmetric by the induction hypothesis and this is enough to
prove that $T$ is symmetric.

Finally if $r=0$ then $v=0$ so that $\Sigma_0$ is essentially ${\mathfrak {g} }$ and we have to
prove that any invariant distribution on ${\mathfrak {g} }$ which is supported on the
nilpotent set is symmetric. We know that the invariant measures carried by the
nilpotent orbits are a basis of this space of distributions so it is enough to
prove that if $X$ is nilpotent then $X$ and $-\overline X$ are conjugate. 

Proposition 1-2, Chapter 4 of [6] asserts that there exists an ${\mathbb F}-$ linear map $g$ from $V$ to $V$ such that $gXg-^{1}=-X$ and $\langle  g(v),g(v')\rangle  =\langle  v',v\rangle  $; in particular $g(\lambda v)=\overline \lambda g(v)$. As usual let us define the conjugation with some basis of $V$ such that the hermitian form has a matrix with coefficients in ${\mathbb F}$ and put $\delta (v)=\overline v$. Then $\delta g$ is ${\mathbb D}-$linear, unitary and $(\delta g)X(\delta g)^{-1}=-\overline X$.

This completes the proof but we assumed the case of the general linear group\dots

\vfill\eject
\centerline{\normalfont \Large \bfseries Appendix: A Frobenius type descent}
\vskip 1cm

If $X$ is a Hausdorff totally disconnected  locally compact topological space (lctd
space in short) we denote by ${\mathcal S}(X)$ the vector space of locally constant
applications with compact support of $X$ into the field of complex numbers ${\mathbb C}$ .
The dual space ${\mathcal S}'(X)$ of ${\mathcal S}(X)$ is the space of distributions on $X$. Once
for all we assume that all the lctd spaces we introduce are countable at infinity. 

If an lctd topological group  $G$ acts continuously on a lctd space $X$ then it acts on
${{\mathcal S}}(X)$ by \[ (gf)(x)=f(g^{-1}x)\] 
and on distributions by
\[ (gT)(f)=T(g^{-1}f)\]
The space of invariant distributions is denoted by ${\mathcal S}'(X)^G$.
Let $G $  be a lctd topological group and $H$ a closed subgroup. Suppose that $H$ acts
continuously on a lctd space $X$. Let $H$ acts on $G\times X$ by
\[ (h,(g,x))\mapsto (gh^{-1},hx)\] 
and let $Y$ be the quotient space. The equivalence relation is open and its graph is closed,
hence $Y$ is Hausdorff and a lctd space. Let $\pi $ be the projection map from $G\times X$ onto $Y$.
The group $G$ acts on $G\times X$ and on $Y$ on the left: $g(g',x)=(gg',x)$ and $g(\pi (g',x))=\pi (gg',x)$

Let $e$ be the neutral element of $G$ and consider the subspace $\{e\}\times X$. We have
\[ \pi ^{-1}\left  (\pi (\{e\}\times X)\right  )=H\times X\]
a closed subspace. Hence $\pi (\{e\}\times X$ is closed in $Y$ and as $\pi $ is open it follows that $X$ is
homeomorphic to $\pi (\{e\}\times X)$. We shall identitfy $X$ with this image. So $X$ is a closed
subspace of $Y$, we have $GX=Y$ and if, for $g\in G$ there exists $x\in X$ such that $gX\in X$ then
$g\in H$. Finally $HX=X$. Note however that $\pi (h((e,x))=\pi ((h,x)=\pi ((e,h^{-1}x))$

 Let $\Delta _H$ be
the module function of $H$ and $\Delta _G$ the module function of $G$ Let ${\mathcal S}_\Delta (G)$ be the
space of functions $\varphi $ defined on $G$, locally constant, such that  
\[ \varphi (gh)=\frac{\Delta _H(h)}{\Delta _G(h)}\varphi (g),\quad g\in G,\, h\in H\] 
and with support compact modulo $H$. Then there exist a positive linear form $\mu $ on this
space, invariant by right translations by elements of $G$; it is unique up to multiplication
by a strictly positive number. We use the notation
\[ \oint _{G/H}\varphi (g)d\mu (g)\] 
If $d_{\ell}g$ is a left Haar measure on $G$ and $d_rh$ a right Haar measure on $H$ we may
normalize in such a way that
\[ \int _Gf(g)=dg=\oint_{G/H}\int _Hf(gh)\Delta _G(h)d_rh\, d\mu (g)\] 
Let ${\mathcal S}'(X)_\Delta ^H$ be the space of distributions $S$ on $X$ such that
\[ \langle  S,\psi (hx)\rangle  =\frac{\Delta _H(h)}{\Delta _G(h)}\langle  S,\psi \rangle  \]

 Let $S\in {\mathcal S}'(X)_\Delta ^H$ and define a distribution $T=\theta (S)$ on $Y$ by
\[ \langle  T,f\rangle  =\oint _{G/H}\langle  S,f(gx)\rangle  d\mu (g)\]
To check that this makes sense consider the continuous map $(g,x)\mapsto gH$ of $G\times X$ onto
$G/H$. It defines a continuous map $\nu$ from $Y$ to $G/H$ . If $y\in Y$ and if $g\in G$ and $x\in X$ are
such that $y=gx$ then $\nu (y)=gH$.

Let $U$ be the support of $f$; it is compact and open. Then $gx\in U$ implies that $gH\in \nu (U)$ so
that the support of $\langle  S,f(gx)\rangle  $ is compact modulo $H$. Furthermore $f$ is fixed by some
open compact subgroup of $G$ and the same is true for $\langle  S,f(gx)\rangle  $; we can apply $\mu $. The
distribution $T$ is invariant under $G$.
\begin{proposition}
 The map $\theta $ is bijective. Furthermore for any $f\in {\mathcal S}(Y)$ there exists
$\varphi \in {\mathcal S}(X)$ such that, for any $S\in {\mathcal S}'(X)^H_\Delta $
\[ \langle  S,\varphi \rangle  =\langle  \theta (S),f\rangle  \]
\end{proposition}
For $F\in {\mathcal S}(G\times X)$ put
\[ f(gx)=\int _HF(gh^{-1},hx)d_rh\]
Then $f$ is well defined and $f\in {\mathcal S}(Y)$. The map $F\mapsto f$ is onto. Also note that
the support of $f$ is contained in the projection of the support of $F$. The transpose map
$T\mapsto U$ is one to one linear from ${\mathcal S}'(Y)$ into ${\mathcal S}'(G\times X)$ and its image is
exactly the subspace of distributions on $G\times X$ which under the action of $H$ satisfy
\[ \langle  U,F(gh^{-1},hx)\rangle  =\Delta _H(h)\langle  U,F(g,x)\rangle  \] 
 The group $G$ operates on $Y$ and also on $G\times X$, by left
multiplication on the first variable ($(g',(g,x))\mapsto (g'g,x)$. Then $T\mapsto U$ is a
bijection from ${\mathcal S}'(Y)^G$ onto the space ${\mathcal S}'(G\times X)^{G,H}$ of distributions on
$G\times X$ invariants by $G$ and having the above invariance property with respect to $H$.

Such a distribution $U$ is uniquely written as a tensor product $d_\ell g\otimes S$ where $S\in {\mathcal
S}'(X)^H_\Delta $. Then $\theta (S)=T$.

Indeed, for $f\in {\mathcal S}(Y)$ choose a ``lift'' $F$. Then
\begin{eqnarray*} \langle  \theta (S),f\rangle  &=&\oint_{G/H}\langle  S,f(gx)\rangle  d\mu (g)\\
&=&\oint
_{G/H}\langle  S,\int _HF(gh^{-1},hx)d_rh\rangle  d\mu (g)\\
&=&\oint_{G/H}\int _H\langle  S,F(gh^{-1},hx)\rangle  d_rh\, d\mu (g)\\
&=&\oint _{G/H}\int _H\langle  S,F(gh^{-1},x)\frac{\Delta _H(h)}{\Delta _G(h)}d_rh\, d\mu (g)\rangle  \cr
&=&\langle  S,\oint_{G/H}\int _HF(gh,x)\Delta _G(h)d_rh\, d\mu (g)\rangle  \\
&=&\langle  S,\int _GF(g,x) d_\ell g\rangle  \\
&=&\langle  d_\ell g\otimes S,F\rangle  \\
&=&\langle  U,F\rangle  \\
&=&\langle  T,f\rangle  
\end{eqnarray*}
Finally for $f$  as above, choose $F$ and define
\[ \varphi (x)=\int _GF(g,x)dg\]
Then, for all $S$ we have $\langle  S,\varphi \rangle  =\langle  \theta (S),f\rangle  $.
\vskip 0.3cm

  For our peace of mind here is a detailed proof.
\begin{lemma}
Let $Y$ be a locally compact totally disconnected space which is countable at
infinity. If \[ Y=\bigcup_{i\in I}U_i\] is an open covering of $Y$, then there exists a finer
covering 
\[ Y=\bigcup_{j\in J}V_j\] 
locally finite and such that each $V_j$ is open and compact. Furthermore we may assume $J$
to be countable.
\end{lemma}
 Let $Y_n$ an increasing sequence of
open compact subsets of $Y$ such that $Y=\cup Y_n$. Let $W_n=K_n-K_{n-1},\,\, n\geqslant 1$ and
$W_0=K_0$. Then each $W_n$ is open and compact. For each $y\in W_n$ choose $i\in I$ such that
$x\in U_i$ and let $V_x$ be an open and compact neighbourhood of $x$ contained in $W_nßU_i$.
We can cover $W_n$ by a finite number of such $V_x$. Letting $n$ vary we get a covering of
$Y$ with the required properties.

With the notations of the lemma, let $g_j$ be the characteristic function of $V_j$. Because
the covering is locally finite, the sum $\sum g_j$
is well defined and strictly positive. Put $\kappa _j=g_j/\sum g_r$. Then each $\kappa _j\in {\mathcal S}(Y)$, the
support of $\kappa _j$ is $V_j$ and $\sum \kappa _j=1$. We shall call the familly of $\kappa _j$ the partition of
unity associated to the covering $V_j$. 

 \begin{lemma}
 Let $H$ be a lctd group,  and $Z$ a lctd space countable at
infinity. Suppose that $H$ acts continuously and properly on $Z$, on the left.There exists a
locally constant function $u$ defined on $Z$, strictly positive whose support has a compact
intersection with the inverse image of any compact subset of $H\setminus Z$ and such that
\[ \int _Hu(hz)d_rh=1\] 
\end{lemma}
On $H\setminus Z$ we put the quotient topology. The equivalence relation is open, hence the
projection map $\pi $ from $Z$ to $H\setminus Z$ is open. Because the action is proper, the quotient is
Hausdorff. Hence this quotient is also a lctd, countable at infinity. Call $\pi $ the projection
map from $Z$ onto $H\setminus Z$.

For any $y\in H\setminus Z$ choose $z$ such that $\pi (z)=y$, an open compact neighbourhood $K_y$ of $z$ and
denote by $u_y$ the characteristic function of $K_y$. Consider the covering
\[ H\setminus Z=\bigcup_{y\in H\setminus Z}\pi (K_y)\]
Apply the first lemma. We get a covering
\[ H\setminus Z=\bigcup_{j\in J}V_j\] 
by open compact subsets, locally finite and for each $j\in J$ we can choose $y_j\in H\setminus Z$ such
that $V_j\subset\pi (K_y)$. Let $(\kappa _j)_{j\in J}$ be the partition of unity associated to the $V_j$. Put
\[ F=\sum _j(\kappa _j\circ \pi )u_{y_j}\] 
This sum is locally finite and $F$ is locally constant. The support condition is satisfied. For
all $z$
\[ \int _HF(hz)d_rh\] 
is strictly positive and locally constant. Thus we can take
\[ u(z)=\frac{F(z)}{\int _HF(hz)d_rh}\]

Keep the notations of this last lemma. For $f\in {\mathcal S}(Z)$ put
\[ \varphi (\pi (z))=\int _Hf(hz)d_rh\]
\begin{lemma}
 The map $f\mapsto \varphi $ is a surjective map from ${\mathcal S}(Z)$ onto ${\mathcal
S(H\setminus Z)}$. The transpose map is a linear bijection from ${\mathcal S}'(H\setminus Z)$
onto the space of distributions $U$ on $Z$ such that
\[ \langle  U,f(hz)\rangle  =\Delta _H(h)\langle  U,f(z)\rangle  \]
\end{lemma}
The support of $\varphi  $ is contained into the projection of the support of $F$, hence is compact.
To prove that $\varphi $ is locally constant it is enough to consider the case where $f$ is the
characteristic function of some open compact suset $K$. Fix $z_0\in Z$ and an open compact
neighbourhood $U$ of $z_0$. The action of $H$ being proper, there exists a compact subset $L$
of $H$ such that, for $h\in H$ the condition $hUK\ne \emptyset  $ implies $h\in L$. Thus, for $z\in U$
\[ \varphi (\pi (z))=\int _Lf(hz)d_rh\] 
For $z\in U$ let
\[ L_z=\{h\in L\, |\, hz\in K\},\quad M_z=\{h\in L\, |\, hz\notin  K\}\] 
They are both open and closed subsets of $L$. For each $h_0\in L_{z_0}$ there exists a
neighbourhood $V_{h_0}$ of $z_0$ in $U$ and a neighbourhood $W_{z_0}$ of $h$ in $L_{z_0}$ such
that $hz\in K$ for $h\in W_{z_0}$ and $z\in V_{z_0}$. But $L_{z_0}$  is compact so we can find a
neighbourhood $V$ of $z_0$  in $U$ such that $hz\in K$ for $z\in V$ and any $h\in L_{z_0}$. This means
that $L_{z_0}\subset L_z$ for $z\in V$. But we can argue in exactly the same way with $M_z$ which,
for $z$ close enough to $z_0$ will give the inclusion $L_z\subset L_{z_0}$. Hence $L_z$ is ``locally
constant'' and so is $\varphi $.
Conversely let us start with $\varphi \in {\mathcal S}(H\setminus Z). $ Fix $u$ as in lemma B and put
$f(z)=\varphi (\pi (z))u(z)$; then $f$  maps to $\varphi $.

Now let $S\in {\mathcal S}'(H\setminus Z)$. We lift it to a distribution $U$ on $z$ by
\[ \langle  U,f\rangle  =\langle  S,\int _Hf(h'z)d_rh'\rangle  \] 
If we replace $f(z)$ by $f(hz)$
\[ \int _Hf(hh'z)d_r(h')=\Delta _H(h)\int _Hf(h'z)d_rh'\] 
so that
\[ \langle  U,f(hz)\rangle  =\Delta _H(h)\langle  U,f(z)\rangle  \] 
Conversely if $U$ is a distribution on $Z$ satisfying the above condition we have to show
that it is 0 on the kernel of the map $f\mapsto \varphi $. So let $f$ be such that
\[ \int _Hf(hz)d_rh=0\] 
We may assume that $f$ is not identically 0.
Let $u$ be as in Lemma A-3 and consider
\[ \int _Hf(z)u(hz)d_rh\] 
Assume that $z$ remains in the support of $f$, a compact subset of $Z$. Then because of the
property of the support of $u$ and the hypothesis that the action of $H$ is proper
the set $K$ of all $h$ such that, for some $z$, $f(z)u(hz)\ne 0$ is open and compact. Note that it
is a neighbourhood of the neutral element of $H$. There exists an open compact subgroup
$K_1$, contained in $K$ and such that 
\[ u(k_1hz)f(z)=u(hz)f(z),\quad k_1\in K_1,\, h\in K,\, f(z)\ne 0\] 
The subset $K$ is a finite union of $K_1-$cosets:
\[ K=\bigcup_{h\in K_1\setminus K}K_1h\] 
Then
 \[ \int _Hf(z)u(hz)d_rh=\int _Kf(z)u(hz)d_rh\]

and
\[ \int _Hf(z)u(hz)d_rh=\sum _{h\in K_1\setminus K}f(z)u(hz){\rm Vol}(K_1)\] 
Then
\begin{eqnarray*} \langle  U,f\rangle  &=&\langle  U,\int _Hf(z)u(hz)d_rh\rangle  \\
&=&\langle  U,\sum _{h\in K_1\setminus K}f(z)u(hz){\rm Vol}(K_1)\rangle  \cr
&=&\sum _{h\in K_1\setminus K}{\rm Vol}(K_1)\langle  U,f(z)u(hz)\rangle  \\
&=&\int _K\langle  U,f(z)u(hz)\rangle  d_rh\\
&=&\int _H\langle  U,f(z)u(hz)\rangle  d_rh\\
&=&\int _H\langle  U,f(h^{-1}z)u(z)\rangle  \Delta _H(h)d_rh\\
&=&\int _H\langle  U,f(hz)u(z)\rangle  d_rh
\end{eqnarray*} 
Again if $K'$ is the set of $h$ such that, for some $z$, $f(hz)u(z)\ne 0$ then $K'$ is open and
compact and contains an open compact subgroup $K'_1$ such that $f(k'hz)u(z)=f(hz)u(z)$ for
$k_1\in K'_1,\, h\in K',\, f(hz)u(z)\ne 0$. replacing the integral by a finite sum we conclude that
\[ \int _H\langle  U,f(hz)u(z)\rangle  d_rh=\langle  U,u(z)\int _Hf(hz)d_rh\rangle  =\langle  U,u\varphi \rangle  \]
which is 0 if $\varphi =0$.

Now let us go back to the situation of Proposition A-1. Apply Lemma A-4 with $X=G\times X$ and $H$
acting by $(g,x)\mapsto (gh^{-1},hx)$. Lemma A-4 provides a justification for the assertions
concerning the map $F\mapsto f$ and its transpose. That $U$ must be a tensor product
$d_\ell \otimes S$ follows from the equality ${\mathcal S}(G\times X)={\mathcal S}(G)\otimes {\mathcal S}(X)$ and the fact
that the left invariant distributions on $G$ are the multiples of the Haar measure. The last
computation in the proof is just a simple aplication of Fubini's theorem.
\vskip 0.5cm
Let $x\in X$ and let $H_x$ be its centralizer in $H$. A relatively invariant measure on $H/H_x$ of
multiplier  $\chi $, a character of $H$, is a non zero measure  $\nu_x$ such that
\[ \int _{H/H_x} \Phi  (sh)d\nu_x(h)=\chi (s)^{-1}\int _{H/H_x}\Phi  (h)d\nu_x(h)\] 
Such a measure exists, and is essentially unique, if and only if
\[ \chi (h)=\frac{\Delta _{H_x}(h)}{\Delta _H(h)},\quad h\in H_x\] 
We are interested in the case where $\chi =\Delta _G/\Delta _H$ so that we must assume that
$\Delta _G=\Delta _{H_x}$ on restriction to $H_x$. This is exactly the condition for the existence of
an invariant measure on $G/H_x$. Also note that the stabilizer of $x$ in $G$ is also $H_x$.

We say that the orbit $Hx$ of $x$ in $X$ carries a relatively invariant measure of
multiplier  $\chi $ if, for $\nu_x$ as above
\[ \int _{H/H_x}|\varphi (hx)|d\nu_x<+\infty ,\quad \varphi \in {\mathcal S}(X)\]
Similarly if $\mu _x$ is an invariant measure on $G/H_x$ we say that  the orbit $Gx$ carries
an invariant measure if and only if
\[ \int _{G/H_x}|f(gx)|d\mu _x(g)<+\infty ,\quad f\in {\mathcal S}(Y)\] 
\begin{proposition}
 The orbit $Hx$ carries a relatively invariant measure  $\nu_x$ of multiplier
$\Delta _G/\Delta _H$ if and only if the orbit $Gx$ carries an invariant measure $\mu _x$. Furthermore, if
this is the case, we can normalize $\mu _x$ in such a way that $\mu _x=\theta (\nu_x)$
\end{proposition}
Suppose that $\nu_x$ and $\mu _x$ exist as measures on $H/H_x$ and $G/H_x$ respectively. Then
\[ \Phi  \mapsto \oint_{G/H}\int _{H/H_x}\Phi  (gh)d\nu_x(h)d\mu (g),\quad \Phi  \in {\mathcal S}(G/H_x)\] 
is well defined and is in fact a left invariant measure. So we can normalize in such a way
that
\[ \int _{G/H_x}\Phi  (g)d\mu _x(g)= \oint_{G/H}\int _{H/H_x}\Phi  (gh)d\nu_x(h)d\mu (g)\] 
Now if we start with $\varphi \in {\mathcal S}(X)$ we can choose $F\in {\mathcal S}(G\times X)$ such that
\[ \varphi (x)=\int _GF(g,x)d_\ell g\]  and then define $f\in {\mathcal S}(Y)$ by
\[ f(gx)=\int _HF(gh^{-1},hx)d_rh\] 
If $\varphi \geqslant 0$ we may choose $F\geqslant 0$, hence $f\geqslant 0$. Conversely if we start with $f$ we choose $F$ and
then compute $\varphi $, both positive if $f$ is positive. In this situation we have
$\langle  \theta (S),f\rangle  =\langle  S,\varphi \rangle  $.

Suppose $f\geqslant 0$. then
\begin{eqnarray*} \int _{G/H_x}f(gx)d\mu _x(g)&=&\int _{G/H_x}\int _HF(gh^{-1},hx)d_rhd\mu _x(g)\\
&=&\oint_{G/H}\int _{H/H_x}\int _HF(gsh^{-1},hx)d_rh\,d\nu_x(s)d\mu (g)\\
&=&\oint_{G/H}\int _{H/H_x}\int _HF(gh^{-1},hsx)d_rh\,d\nu_x(s)d\mu (g)\\
&=&\oint_{G/H}\int _H\int _{H/H_x}F(gh^{-1},hsx)d_rh\,d\nu_x(s)d\mu (g)\\
&=&\oint_{G/H}\int _H\int _{H/H_x}F(gh^{-1},sx)\frac{\Delta _H(h)}{\Delta _G(h)}d_rh\,d\nu_x(s)d\mu (g)\\
&=&\oint_{G/H}\int _H\int _{H/H_x}F(gh,sx)\Delta _G(g)d_rh\,d\nu_x(s)d\mu (g)\\
&=&\int _G\int _{H/H_x}F(g,sx)d\nu_x(s)d_\ell g\\
&=&\int _{H/H_x}\varphi (sx)d\nu_x(s)
\end{eqnarray*} 
This shows that the first and last integrals are both finite or infinite, proving the first
assertion of the proposition. The second one follows from the same computation with
complex valued functions.}

 \vskip 1cm
 S.Rallis\hfill\break
 Department of mathematics Ohio State Unibversity,\hfill\break
 COLUMBUS, OH 43210\hfill\break
 E-mail: haar@math.ohio-state.edu
 \vskip 1cm
 G.Schiffmann\hfill\break
  Institut de Recherche Math\'ematique Avanc\'ee,\hfill\break
  Universit\'e Louis-Pasteur et CNRS, 7, rue Ren\'e-Descartes,\hfill\break
   67084 Strasbourg Cedex, France.\hfill\break
  e-mail: schiffma@math.u-strasbg.fr  
\end{document}